\documentclass[a4,11pt]{article}
\usepackage[subrefformat=parens]{subcaption}
\usepackage[dvips]{graphicx}
\usepackage{amssymb}
\usepackage{amstext}
\usepackage{amsmath}
\usepackage{amscd}
\usepackage{amsthm}
\usepackage{amsfonts}
\usepackage{enumerate}
\usepackage{graphicx}
\usepackage{latexsym}
\usepackage{mathrsfs}
\usepackage{mathtools}
\usepackage{cases}
\usepackage[svgnames]{xcolor}
\usepackage{verbatim}
\usepackage{bm}
\usepackage{tikz}
\usepackage{caption}
\usepackage{wrapfig}
\usepackage[top=3.5cm, bottom=3.5cm, left=3.5cm, right=3.5cm]{geometry}
\newtheorem{thm}{Theorem}[section]

\newtheorem{lem}[thm]{Lemma}

\theoremstyle{definition}
\newtheorem{rem}[thm]{Remark}
\newtheorem{dfn}[thm]{Definition}

\newtheorem*{claim*}{Claim}

\numberwithin{equation}{section}

\title{\textbf{Remarks on the factorization and monotonicity method for inverse acoustic scatterings}}

\author{Takashi FURUYA}

\date{}
\begin{document}
\maketitle
\begin{abstract}
We study the factorization and monotonicity method for inverse acoustic scattering problems. Firstly, we give a new general functional analysis theorem for the monotonicity method. Comparing with the factorization method, the general theorem of the monotonicity generates reconstruction schemes under weaker {\it a priori} assumptions for unknown targets, and can directly deal with mixed problems so that the unknown targets have several different boundary conditions. Using the general theorem, we give the reconstruction scheme for the mixed crack that the Dirichlet boundary condition is imposed on one side of the crack and the Neumann boundary condition on the other side, which is a new extension of monotonicity method.
\end{abstract}
\section{Introduction}
In this paper, we study the factorization and monotonicity method for inverse acoustic scattering problems. The factorization method has first been introduced by Kirsch (\cite{Kirsch1}) for the inverse acoustic  obstacle scattering. It has been studied so far by many authors (see e.g., \cite{Anagnostopoulos and Charalambopoulos and Kleefeld, Bondarenko and Kirsch and Liu, Boukari and Haddar, Furuya1, Grinberg and Kirsch, Kirsch2, Kirsch and Grinberg, Kirsch and Liu, Kirsch and Ritter, Lechleiter, Wu and Yan}), and it is well-known as one of classical qualitative methods, which includes the linear sampling method of Colton and Kirsch (\cite{Colton and Kirsch}), the singular sources method of Potthast (\cite{Potthast}), the probe method of Ikehata (\cite{Ikehata}), etc. The monotonicity method, on the other hand, has been recently introduced by Harrach in \cite{Harrach and Ullrich} for the electrical impedance tomography. Very recently, it was extended to the Helmholtz equation in the bounded domain \cite{Harrach and Pohjola and Salo1, Harrach and Pohjola and Salo2}, the inverse acoustic obstacle \cite{Albicker and Griesmaier}, Dirichlet crack \cite{Furuya3}, and medium scatterings \cite{Griesmaier and Harrach, Lakshtanov and Lechleiter}. It was found from these works that the monotonicity method has advantage over the factorization method that we can give the reconstruction scheme under weaker {\it a priori} assumptions for unknown targets. 
\par
The contributions of this paper are the following.
\begin{itemize}
  \item[\bf{[A]}] We give a general functional analysis theorem for the monotonicity method (Theorem \ref{General theorem for MM}) including previous works \cite{Albicker and Griesmaier, Furuya3, Griesmaier and Harrach} for inverse acoustic scatterings.
  \item[\bf{[B]}] Using the general theorem, we give the reconstruction scheme for the inverse mixed crack scattering, which is a new  extension of the monotonicity method (Theorems \ref{Mixed crack part1} and \ref{Mixed crack part2}).
\end{itemize}
\par
A characteristic of the factorization method is to prepare the general functional analysis theorem which generates reconstruction schemes by spectrums of the far-field operator (see Theorem 2.15 of \cite{Kirsch and Grinberg}). Based on this idea, the contribution {\bf[A]} will be discussed as the version of the monotonicity method. 
\par
The general theorem of the factorization method assumes that the real part of the middle operator of the far-field operator has a decomposition into a positive coercive operator and a compact operator, while the imaginary part of the middle operator becomes strictly positive (see the assumptions {\bf(b)} and {\bf(c)} in Theorem \ref{General theorem for FM}). This two assumptions cause {\it a priori} assumptions for the unknown target, in particular, the positivity of the imaginary part corresponds to restrictions for the wave number. For example, it is necessary for the inverse medium scattering that the wave number is not a transmission eigenvalue with respect to the unknown medium (see e.g., Theorem 4.10 of \cite{Kirsch and Grinberg}), and for the inverse obstacle scattering that the wave number is not an eigenvalue with respect to the boundary condition of the unknown obstacle (see e.g., Corollary 2.16 of \cite{Kirsch and Grinberg}). However, the general theorem of the monotonicity method does not assume the positivity of the imaginary part (see Theorem \ref{General theorem for MM}), which means that the monotonicity can essentially avoid restrictions for the wave number. In fact, monotonicity reconstructions for inverse obstacle (Theorem 5.3 of \cite{Albicker and Griesmaier}) and medium (Theorems 5.1--5.3 of \cite{Griesmaier and Harrach}) have been successful without restrictions for the wave number.
\par
The advantage of the monotonicity over the factorization is not only weaker {\it a priori} assumptions, but also to directly deal with mixed problems that the unknown target consists of two separate components with different boundary conditions. The real part of middle operator of the far-field operator for mixed problems is not decomposed into a positive coercive operator and a compact operator (see e.g., Theorem 3.4 of \cite{Kirsch and Grinberg}, Theorem 3.2 of \cite{Kirsch and Liu}). In order to make a decomposition of the coercivity and the compactness, the factorization method for mixed problems needs masking approaches that one component is covered by an artificial domain disjoint with the other component we want to reconstruct (see e.g., Lemma 3.5 of \cite{Kirsch and Grinberg}, (3.26) of \cite{Kirsch and Liu}). However, the general theorem of the monotonicity method (see {\bf(3)} of Theorem \ref{General theorem for MM}) does not assume such a real part decomposition, which means that the monotonicity essentially do not need masking approaches. In fact, monotonicity reconstructions for inverse mixed obstacle (Theorem 5.5 of \cite{Albicker and Griesmaier}) have been successful without the masking approach. 
\par
This paper studies not only previous works from the viewpoint of the general theorem, but also a new extension of the monotonicity to the inverse acoustic mixed crack scattering, which corresponds to the contribution {\bf[B]}. The mixed crack consists of only one component, but imposes the Dirichlet boundary condition on one side of the unknown crack and the Neumann boundary condition on the other side (see the beginning of Section 5.5). The factorization method for the mixed crack has been studied in \cite{Wu and Yan}, but an extensive closed curve of the unknown crack should be known, which is a very restrictive assumption (see Theorem 3.3 of \cite{Wu and Yan}). Using the general theorem for the monotonicity method, we give a reconstruction scheme without assuming such an extensive curve (see Theorems \ref{Mixed crack part1} and \ref{Mixed crack part2}).
\par
This paper is organized as follows. In Section 2, we define the inverse acoustic scattering problem. In Section 3, we recall the functional analysis theorem for the factorization method. In Section 4, we give a new general functional analysis theorem (Theorem \ref{General theorem for MM}) for the monotonicity method. In Section 5, we study several applications of the general theorem. The results discussed in Sections 5.1--5.4 are the same in previous works \cite{Albicker and Griesmaier, Griesmaier and Harrach, Furuya3}, while the one in Section 5.5 is new in the monotonicity for the inverse mixed crack scattering. Finally in Section 6, we give numerical examples for our theoretical results.
\section{Inverse acoustic scattering}
First of all, we define the inverse acoustic scattering problem. Let $k>0$ be the wave number, and let $\theta \in \mathbb{S}^{d-1}$ ($d=2,3$) be the incident direction. We denote the incident field $u^{inc}(x, \theta)$ with incident direction $\theta$ by the plane wave of the form 
\begin{equation}
u^{inc}(x, \theta):=\mathrm{e}^{ikx \cdot \theta}, \ x \in \mathbb{R}^d.
\end{equation}
Let $\Omega \subset \mathbb{R}^{d}$ be a bounded open set with the smooth boundary $\partial \Omega$ such that the exterior $\mathbb{R}^d \setminus  \overline{\Omega}$ is connected. Here, we denote by $\overline{\Omega}=\Omega \cup \partial \Omega$. In particular, we discuss the following two cases. The first case is that the scatterer $\Omega$ is a penetrable medium, and determine the total field $u=u^{sca}+u^{inc}$ such that
\begin{equation}
\Delta u+k^2(1+q)u=0 \ \mathrm{in} \ \mathbb{R}^d, \label{Medium}
\end{equation}
\begin{equation}
\lim_{r:=|x| \to \infty} r^{\frac{d-1}{2}} \biggl( \frac{\partial u^{sca}}{\partial  r}-iku^{sca} \biggr)=0, \label{SRC}
\end{equation}
where $q \in L^{\infty}(\mathbb{R}^d)$ has a compact support such that $\Omega = \mathrm{supp}\ q$, and $\Delta$ is the {\it Laplace operator}. The {\it Sommerfeld radiation condition} (\ref{SRC}) holds uniformly in all directions $\hat{x}:=\frac{x}{|x|}$. The second case is that $\Omega$ is an impenetrable obstacle, and determine the total field $u=u^{sca}+u^{inc}$ such that
\begin{equation}
\Delta u+k^2 u=0 \ \mathrm{in} \ \mathbb{R}^d \setminus  \overline{\Omega}, \label{Obstacle}
\end{equation}
\begin{equation}
\mathcal{B}u=0 \ \mathrm{on} \ \partial \Omega, \label{BC}
\end{equation}
and $u^{sca}$ satisfies the Sommerfeld radiation condition (\ref{SRC}) where (\ref{BC}) means the boundary conditions, for example, the Dirichlet boundary condition $\mathcal{B}u=u$, the Neumann boundary condition $\mathcal{B}u=\frac{\partial u}{\partial \nu}$, etc. In both cases, it is well-known that there exists a unique solution $u^{sca}$ that has the following asymptotic behaviour (see e.g., \cite{Colton and Kress})
\begin{equation}
u^{sca}(x)=\frac{\mathrm{e}^{ikr}}{r^{\frac{d-1}{2}}}\Bigl\{ u^{\infty}(\hat{x},\theta)+O\bigl(1/r \bigr) \Bigr\} , \ r \to \infty.
\end{equation}
The function $u^{\infty}$ is called the {\it far-field pattern} of the scattered field $u^{sca}$. With the far-field pattern $u^{\infty}$, we define the far-field operator $F :L^{2}(\mathbb{S}^{d-1}) \to L^{2}(\mathbb{S}^{d-1})$ by
\begin{equation}
Fg(\hat{x}):=\int_{\mathbb{S}^{d-1}}u^{\infty}(\hat{x},\theta)g(\theta)ds(\theta), \ \hat{x} \in \mathbb{S}^{d-1}. \label{FFO}
\end{equation}
In the inverse acoustic scattering problem, we reconstruct $\Omega$ from the far-field pattern $u^{\infty}(\hat{x},\theta)$ for all $\hat{x},\theta \in \mathbb{S}^{d-1}$, and fixed $k>0$. In other words, given the far-field operator $F$, we reconstruct $\Omega$.
\section{The factorization method} 
Here, we recall the general functional analysis theorem for the factorization method. The following functional analytic theorem is proved by Theorem 2.15 of \cite{Kirsch and Grinberg} and Theorem 3.1 of \cite{Furuya2}.
\begin{thm}[Theorem 2.15 of \cite{Kirsch and Grinberg} and Theorem 3.1 of \cite{Furuya2}]\label{General theorem for FM}
Let $X \subset U\subset X^{*}$ be a Gelfand triple with a Hilbert space $U$ and a reflexive Banach space $X$ such that the embedding is dense. Furthermore, let Y be a Hilbert space and let $F:Y \to Y$, $G:X \to Y$, $T:X^{*} \to X$ be linear bounded operators such that 
\begin{equation}
F=GTG^{*}.
\end{equation}
We make the following assumptions: 
\begin{description}

\item[(a)] $G$ is compact with dense range in $Y$. 

\item[(b)] $\mathrm{Re}T$ has the form $\mathrm{Re}T=C+K$ where $K:X^{*} \to X$ is some self-adjoint compact operator and $C: X^{*} \to X$ is some positive coercive operator, i.e., there exists a constant $c>0$ such that
\begin{equation}
\langle \varphi,  C \varphi \rangle \geq c \left\| \varphi \right\|_{X^{*}}^2 \ for \ all \ \varphi \in X^{*},
\end{equation}
where $\langle \cdot , \cdot \rangle$ denotes the duality pairing between $X^{*}$ and $X$. \label{2}

\item[(c)]
$\mathrm{Im} \langle \varphi, T \varphi \rangle > 0$ for all $\varphi \in \overline{\mathrm{Ran}(G^{*})}$ with $\varphi \neq 0$.
\end{description}
Then, the operator $F_{\#}:=\bigl|\mathrm{Re}F\bigr|+\mathrm{Im}F$ is non-negative, and the ranges of $G:X \to Y$ and $F_{\#}^{1/2}:Y \to Y$ coincide with each other, that is, we have the following range identity;
\begin{equation}
\mathrm{Ran}(G)=\mathrm{Ran}(F_{\#}^{1/2}). \label{range identity}
\end{equation}
\end{thm}
Here, the real part and the imaginary part of an operator $A$ are self-adjoint operators given by
\begin{equation}
\mathrm{Re}A=\displaystyle \frac{A+A^{*}}{2} \ \ \ \mathrm{and} \ \ \ \mathrm{Im}A=\displaystyle \frac{A-A^{*}}{2i}.
\end{equation}
\begin{rem}
It has been well-known that in Theorem 2.1 of \cite{Lechleiter} the assumption {\bf(c)} can be replaced by the injectivity of $T$, and it has been mainly used especially for the relaxation of the assumption that the wave number $k>0$ is not a transmission eigenvalue in inverse medium scatterings (see e.g., \cite{Furuya1, Kirsch and Grinberg, Kirsch and Liu, Lechleiter}). However, it was found that this replacement is not correct (see Remark 3.2 of \cite{Furuya2}), thus the factorization method for inverse medium scatterings essentially needs the assumption of transmission eigenvalues.
\end{rem}
\section{General theorems for the monotonicity method}
In this section, we give a new general functional analysis theorem for the monotonicity method.
\begin{dfn}
Let $A, B:H \to H$ be self-adjoint compact linear operators on a Hilbert space $H$. We write
\begin{equation}
 A\leq_{\mathrm{fin}} B,\label{2.1}
\end{equation}
if $B-A$ has only finitely many negative eigenvalues.
\end{dfn}

\begin{thm}\label{General theorem for MM}
Let $X \subset U\subset X^{*}$, $\tilde{X} \subset \tilde{U}\subset \tilde{X}^{*}$ be Gelfand triples with Hilbert spaces $U$, $\tilde{U}$ and reflexive Banach spaces $X$, $\tilde{X}$ such that the embeddings are dense. Furthermore, let Y be a Hilbert space and let $F:Y \to Y$,  $\tilde{F}:Y \to Y$, $G:X \to Y$, $\tilde{G}:\tilde{X} \to Y$, $T:X^{*} \to X$, $\tilde{T}:\tilde{X}^{*} \to \tilde{X}$ be linear bounded operators such that 
\begin{equation}
F=GTG^{*},\ \ \  \tilde{F}=\tilde{G}\tilde{T}\tilde{G}^{*}. \label{Factorization}
\end{equation}

\begin{description}

\item[(1)] Assume that
\begin{description}
\item[(1a)] $\mathrm{Re}T$ has the form $\mathrm{Re}T=C+K$ where $C:X^{*} \to X$ is some positive coercive operator and $K:X^{*} \to X$ is some self-adjoint compact operator.

\item[(1b)]There exists a compact operator $R:\tilde{X} \to X$ such that $\tilde{G}=GR$.
\end{description}
Then, 
\begin{equation}
\mathrm{Re}\tilde{F} \leq_{\mathrm{fin}} \mathrm{Re}F. \label{conclusion(1)}
\end{equation}

\item[(2)] Assume that 
\begin{description}
\item[(2a)] $\mathrm{Re}\tilde{T}$ has the form $\mathrm{Re}\tilde{T}=\tilde{C}+\tilde{K}$ where $\tilde{C}:\tilde{X}^{*} \to \tilde{X}$ is some positive coercive operator and $\tilde{K}:\tilde{X}^{*} \to \tilde{X}$ is some self-adjoint compact operator.

\item[(2b)]There exists an infinite dimensional subspace $W$ in $\mathrm{Ran}(\tilde{G})$ such that $W \cap \mathrm{Ran}(G) = \{ 0 \}$. 
\end{description}
Then, 
\begin{equation}
\mathrm{Re}\tilde{F} \not\leq_{\mathrm{fin}} \mathrm{Re}F.\label{conclusion(2)}
\end{equation}

\item[(3)] Let $X_{j} \subset U_{j}\subset X^{*}_{j}$ ($j=1,2$) be a Gelfand triple with a Hilbert space $U_{j}$ and a reflexive Banach space $X_j$ such that the embedding is dense. Let $F^{Mix}:Y \to Y$, $G^{Mix}:X_{1} \times X_{2}^{*} \to Y$, $T^{Mix}:X_{1}^{*} \times X_{2} \to X_{1} \times X_{2}^{*}$ be linear bounded operators such that 
\begin{equation}
F^{Mix}=G^{Mix}T^{Mix}G^{Mix \ *}.
\end{equation}
Assume that
\begin{description}
\item[(3a)] $\mathrm{Re}T^{Mix}$ has the form $\mathrm{Re}T^{Mix}=\left(\begin{array}{cc}
      C_{11} &C_{12} \\
      C_{21} &C_{22} 
    \end{array}\right) + K^{Mix}$ where $C_{11}:X_{1}^{*} \to X_{1}$ is some positive coercive operator, $C_{12}:X_{2} \to X_{1}$ and $C_{21}:X_{1}^{*} \to X_{2}^{*}$ are some linear bounded operators, $C_{22}:X_{2} \to X_{2}^{*}$ is some negative coercive operator (that is, $-C_{22}$ is positive coercive), and $K^{Mix}:X_{1}^{*} \times X_{2} \to X_{1} \times X_{2}^{*}$ is some self-adjoint compact operator. 

\item[(3b)] There exists an infinite dimensional subspace $W_{1}$ in  $\mathrm{Ran}(G^{Mix}R_{1}^{*})$ such that $W_{1} \cap \mathrm{Ran}\left( \left[ \tilde{G}, G^{Mix}R_{2}^{*} \right] \right) = \{ 0 \}$ where $R_{1}:X_{1}^{*}\times X_{2} \to X_{1}^{*}$ is defined by $R_{1}\left(
    \begin{array}{ccc}
      f \\
      g \\
    \end{array}
\right):=f$, and its adjoint operator $R^{*}_{1}:X_{1} \to X_{1} \times X_{2}^{*}$ is given by $R_{1}^{*}\varphi:=\left(
    \begin{array}{ccc}
      \varphi \\
      0 \\
    \end{array}
\right)$.
\end{description}
\end{description}
Then, 
\begin{equation}
\mathrm{Re}F^{Mix} \not\leq_{\mathrm{fin}} \mathrm{Re}\tilde{F}. \label{conclusion(3)}
\end{equation}
\end{thm}
\begin{rem}
If the assumption {\bf(3b)} is replaced by 
\begin{description}
\item[(3b)'] {\it There exists a finite dimensional subspace $W_{2}$ in  $\mathrm{Ran}(G^{Mix}R_{2}^{*})$ such that $W_{2} \cap \mathrm{Ran}\left( \left[ \tilde{G}, G^{Mix}R_{1}^{*} \right] \right) = \{ 0 \}$ where $R_{2}:X_{1}^{*}\times X_{2} \to X_{2}$ is defined by $R_{2}\left(
    \begin{array}{ccc}
      f \\
      g \\
    \end{array}
\right):=g$, and its adjoint operator $R^{*}_{2}:X^{*}_{2} \to X_{1} \times X_{2}^{*}$ is given by $R_{2}^{*}\psi:=\left(
    \begin{array}{ccc}
      0 \\
      \psi \\
    \end{array}
\right)$,
}
\end{description}
then, we can show that $-\mathrm{Re}F^{Mix} \not\leq_{\mathrm{fin}} \mathrm{Re}\tilde{F}$ by the same argument.
\end{rem}
We recall the following technical lemmas which will be useful to prove Theorem \ref{General theorem for MM}.
\begin{lem}[Corollary 3.3 of \cite{Harrach and Pohjola and Salo2}]\label{equivalence of definition}
Let $A, B:H \to H$ be self-adjoint compact linear operators on a Hilbert space $H$ with an inner product $( \cdot, \cdot )_H$. Then, the following statements are equivalent:
\begin{description}
\item[(1)]
$A\leq_{\mathrm{fin}} B$
\item[(2)]
There exists a finite dimensional subspace $V$ in $H$ such that
\begin{equation}
\left( (B-A)v, v \right)_H \geq0,
\end{equation}
for all $v \in V^{\bot}$.
\end{description}
\end{lem}

\begin{lem}[Lemma 4.6 in \cite{Harrach and Pohjola and Salo2}]\label{Range inclusion}
Let $X$, $Y$, and $Z$ be Hilbert spaces, and let $A:X \to Y$ and $B:X \to Z$ be bounded linear operators. Then,
\begin{equation}
\exists C>0: \ \left\| Ax \right\|_{Y}^2 \leq  C\left\| Bx \right\|_{Z}^2 \ for \ all \ x \in X \ \ \ \ \Longleftrightarrow \ \ \ \ \mathrm{Ran}(A^{*})\subset \mathrm{Ran}(B^{*}).
\end{equation}
\end{lem}

\begin{lem}[Lemma 4.7 in \cite{Harrach and Pohjola and Salo2}]\label{Vector space lemma}
Let $X$, $Y$, $V \subset Z$ be subspaces of a vector space $Z$. If 
\begin{equation}
X\cap Y = \{ 0 \}, \ \ \ \ and \ \ \ \ X \subset Y+V,
\end{equation}
then, $\mathrm{dim}(X) \leq \mathrm{dim}(V)$.
\end{lem}
\vspace{0.3cm}
\begin{proof}[\bf{Proof of Theorem \ref{General theorem for MM}}]
{\bf (1)} Since the restriction $C\bigl|_{U}:U \to U$ is positive, there exists a positive square root $\hat{W}:U \to U$, i.e., $C\bigl|_{U}=\hat{W}^2$. Since we have,
\begin{equation}
\| \hat{W}\varphi \|^{2}_{U}=(\varphi, \hat{W}^2 \varphi )_{U}=( \varphi, C\bigl|_{U} \varphi )_{U}=\left< \varphi, C \varphi \right> \leq \| C \| \| \varphi \|_{X^{*}}^{2}, \ \ for \ all\ \varphi \in U,
\end{equation}
and the embedding $U \subset X^{*}$ is dense, $\hat{W}$ has a bounded extension $W:X^{*} \to U$ of $\hat{W}$. By the positive coercivity of $C$, there exists a constant $c>0$ such that for all $\varphi \in U$,
\begin{equation}
c\| \varphi \|^{2}_{X^{*}} \leq \langle \varphi,  C \varphi \rangle =( \varphi, C\bigl|_{U} \varphi )_{U} = \| \hat{W} \varphi \|_{U}^{2}.
\end{equation}
Hence, by the dense embedding $U \subset X^{*}$, we have $c\| \varphi \|^{2}_{X^{*}} \leq \| W \varphi \|_{U}^{2}$ for all $\varphi \in X^{*}$, which implies that the extension $W: X^{*} \to U$ of $\hat{W}$ is bounded invertible. It is easy to check that $C=W^{*}W$. By this, the factorization (\ref{Factorization}) of operators $F$ and $\tilde{F}$, assumptions {\bf(1a)} and {\bf(1b)}, we have
\begin{eqnarray}
&&\mathrm{Re}F - \mathrm{Re}\tilde{F}
=G \left[ C + \mathrm{Re}K - R (\mathrm{Re}\tilde{T})R^{*} \right]G^{*}
\nonumber\\
&=& \left[GW^{*}\right] \left[ W^{*^{-1}} C W^{-1} +  W^{*^{-1}}\left\{ \mathrm{Re}K - R (\mathrm{Re}\tilde{T})R^{*} \right\} W^{-1} \right]\left[GW^{*}\right]^{*}
\nonumber\\
&=:& \hat{G} [ I_{U}+\hat{K} ] \hat{G}^{*}.
\end{eqnarray}
\par
Let $\left\{ \mu_j, \phi_j \right\}$ be an eigensystem of the self-adjoint compact operator $\hat{K}:U \to U$, and let 
\begin{equation}
V:=\mathrm{span}\left\{ \phi_{j}: \mu_j \leq -\displaystyle \frac{1}{2} \right\}. \label{spectrum1}
\end{equation}
Then, $V$ is a finite dimensional subspace of $U$, and for all $\varphi \in \left[\mathrm{Ran}\left( \hat{G}\bigl|_{V} \right) \right]^{\bot}$, which is equivalent to $\hat{G}^{*}\varphi \in V^{\bot}=\mathrm{span}\left\{ \phi_{j}: \mu_j > -\displaystyle \frac{1}{2} \right\}$,
\begin{equation}
( ( \mathrm{Re}F - \mathrm{Re}\tilde{F} ) \varphi, \varphi )_{Y} = ( ( I_{U} + \hat{K} ) \hat{G}^{*} \varphi, \hat{G}^{*} \varphi )_{U} \geq \displaystyle \frac{1}{2} \| \hat{G}^{*} \varphi \|_{U}^{2} \geq 0.\label{spectrum2}
\end{equation}
From $\mathrm{dim}\left[\mathrm{Ran}\left( \hat{G}\bigl|_{V} \right)\right] < \infty$ and Lemma \ref{equivalence of definition}, we conclude (\ref{conclusion(1)}).
\par
{\bf(2)} Assume on the contrary that $\mathrm{Re}\tilde{F} \leq_{\mathrm{fin}} \mathrm{Re}F.$ Then by Lemma \ref{equivalence of definition}, there exists a finite dimensional subspace $V_{1}$ in $Y$ such that
\begin{equation}
( \mathrm{Re}\tilde{F}\varphi, \varphi )_{Y} \leq \left( \mathrm{Re}F \varphi, \varphi \right)_{Y}, \label{(2)-1}
\end{equation}
for all $\varphi \in V^{\bot}_{1}$. By the same argument in the beginning of {\bf (1)}, there exists a bounded invertible operator $\tilde{W}:\tilde{X}^{*} \to \tilde{U}$ and a self-adjoint compact operator $K$ such that 
\begin{equation}
\mathrm{Re}\tilde{F}= [\tilde{G}\tilde{W}^{*}] [ I_{\tilde{U}}+K] [\tilde{G}\tilde{W}^{*}]^{*},
\end{equation}
which implies that by the same argument in (\ref{spectrum1})--(\ref{spectrum2}) there exists a finite dimensional subspace $V_{2}$ in $Y$ and a constant $c>0$ such that
\begin{equation}
( \mathrm{Re}\tilde{F} \varphi, \varphi)_{Y} \geq \displaystyle \frac{1}{2} \| \tilde{W} \tilde{G}^{*} \varphi \|_{\tilde{U}}^{2} \geq c \| \tilde{G}^{*} \varphi \|_{\tilde{X}^{*}}^{2}, \label{(2)-2}
\end{equation}
for all $\varphi \in V^{\bot}_{2}$. Setting $V:=V_{1}\cup V_{2}$, $V$ is a finite dimensional subspace in $Y$. Then by (\ref{(2)-1}) and (\ref{(2)-2}) we have 
\begin{equation}
 c \| \tilde{G}^{*} \varphi \|_{\tilde{X}^{*}}^{2} \leq ( \mathrm{Re}F \varphi, \varphi)_{Y} \leq \| \mathrm{Re}T \|  \| G^{*} \varphi \|_{X^{*}}^{2}, \label{(2)-3}
\end{equation}
for all $\varphi \in V^{\bot}$.
\par
On the other hand, by the assumption {\bf (2b)} and Lemma \ref{Vector space lemma}, we have $W \not \subset \mathrm{Ran}(G) + V$, which implies that by $W \subset \mathrm{Ran}(\tilde{G})$
\begin{equation}
\mathrm{Ran}(\tilde{G}) \not \subset \mathrm{Ran}(G) + V = \mathrm{Ran}\left([G, P_{V}] \right),
\end{equation}
where $P_{V}$ denotes the orthogonal projection on $V$. By this and Lemma \ref{Range inclusion}, there exists a sequence $\displaystyle (\varphi_{n})_{n \in \mathbb{N} } \subset Y$ such that 
\begin{equation}
\| \tilde{G}^{*} \varphi_{n} \|_{\tilde{X}^{*}}^{2} > n^{2} \left( \| G^{*} \varphi_{n} \|_{X^{*}}^{2} + \| P_{V} \varphi_{n} \|_{Y}^{2} \right),
\end{equation}
for all $n \in \mathbb{N}$. Setting $\psi_{n}:=\displaystyle \frac{ n^{\frac{1}{2}} \varphi_n}{\| \tilde{G}^{*} \varphi_{n} \|_{\tilde{X}^{*}}}$, we obtain
\begin{equation}
 \| G^{*} \psi_{n} \|_{X^{*}}^{2} + \| P_{V} \psi_{n} \|_{Y}^{2} = n\displaystyle \frac{ \| G^{*} \varphi_{n} \|_{X^{*}}^{2} + \| P_{V} \varphi_{n} \|_{Y}^{2}}{\| \tilde{G}^{*} \varphi_{n} \|^{2}_{\tilde{X}^{*}}} \leq \frac{1}{n}, \ \ \ \  \| \tilde{G}^{*} \psi_{n} \|_{X^{*}}^{2} = n.
\end{equation}
Setting $\tilde{\psi}_{n}:=(I-P_{V})\psi_{n} \in V^{\bot}$, we finally obtain
\begin{equation}
\| \tilde{G}^{*} \tilde{\psi}_{n} \|_{\tilde{X}^{*}} \geq \| \tilde{G}^{*} \psi_{n} \|_{\tilde{X}^{*}} - \|\tilde{G}^{*} \| \| P_{V} \psi_{n} \|_{Y} \to \infty, \ \ \ \mathrm{as} \ n \to \infty,
\end{equation}
\begin{equation}
\| G^{*} \tilde{\psi}_{n} \|_{X^{*}} \leq \| G^{*} \psi_{n} \|_{X^{*}} + \|G^{*} \| \| P_{V} \psi_{n} \|_{Y}^{2} \to 0, \ \ \ \mathrm{as} \ n \to \infty,
\end{equation}
which contradicts (\ref{(2)-3}). Therefore, we conclude (\ref{conclusion(2)}).
\par 
{\bf(3)} Assume on the contrary that $\mathrm{Re}F^{Mix} \leq_{\mathrm{fin}} \mathrm{Re}\tilde{F}$. Then by Lemma \ref{equivalence of definition}, there exists a finite dimensional subspace $V_{1}$ in $Y$ such that
\begin{equation}
( \mathrm{Re}F^{Mix}\varphi, \varphi )_{Y} \leq ( \mathrm{Re}\tilde{F} \varphi, \varphi )_{Y}, \label{(3)-1}
\end{equation}
for all $\varphi \in V^{\bot}_{1}$. Since $C_{11}: X_{1}^{*} \to X_{1}$ and $-C_{22}^{-1}: X_{2}^{*} \to X_2$ are positive coercive, by the same argument in the beginning of {\bf (1)} there exists a bounded invertible operator $W_{jj}:X^{*}_{j} \to U_{j}$ such that 
\begin{equation}
C_{11}=W_{11}^{*}W_{11}, \ \ \ \ \ -C_{22}^{-1}=W_{22}^{*}W_{22},
\end{equation}
We denote by the operator $W:=\left(\begin{array}{cc}
      W_{11}^{*} &0 \\
      0 &W_{22}^{-1} 
    \end{array}\right):U_{1} \times U_{2} \to X_{1} \times X_{2}^{*}$, hence, we have
\begin{eqnarray}
\mathrm{Re}F^{Mix} &=& [G^{Mix}W] \left[ W^{-1} \left\{ \left(\begin{array}{cc}
      C_{11} &C_{12} \\
      C_{21} &C_{22} 
    \end{array}\right) + K^{Mix} \right\} W^{*^{-1}} \right][G^{Mix}W]^{*}
\nonumber\\
&=:& \hat{G}^{Mix} \left[ \left(\begin{array}{cc}
      I_{U_{1}} & \hat{C}_{12} \\
      \hat{C}_{21} & -I_{U_{2}} 
    \end{array}\right) +  \hat{K}^{Mix} \right]\hat{G}^{Mix \ *},
\end{eqnarray}
where $\hat{C}_{12}:= W_{11}^{*^{-1}}C_{12}W_{22}^{*}:U_{2} \to U_{1}$, $\hat{C}_{21}:= W_{22}C_{21}W_{11}^{*}:U_{1} \to U_{2}$, and $\hat{K}^{Mix}=W^{-1}K^{Mix}W^{*^{-1}}:U_{1} \times U_{2} \to U_{1} \times U_{2}$. Since $\hat{K}^{Mix}$ is a self-adjoint compact operator, by the same argument in (\ref{spectrum1})--(\ref{spectrum2}), there exists a finite dimensional subspace $V_2$ in Y, and constants $c_1, c_2, c_3 >0$ such that
\begin{eqnarray}
&&(\mathrm{Re}F^{Mix}\varphi, \varphi)_{Y} 
\nonumber\\
&=& (\left[ \left(\begin{array}{cc}
      I_{U_{1}} & \hat{C}_{12} \\
      \hat{C}_{21} & -I_{U_{2}} 
    \end{array}\right) +  \hat{K}^{Mix} \right] \hat{G}^{Mix \ *} \varphi, \hat{G}^{Mix \ *} \varphi)_{U_{1}\times U_{2}}
\nonumber\\
&\geq& \| W_{11} R_{1} G^{Mix\ *} \varphi \|^{2}_{U_{1}} - \| W_{22}^{*^{-1}} R_{2} G^{Mix\ *} \varphi \|^{2}_{U_{2}} 
\nonumber\\
&&-\frac{1}{2} \left\{  \| W_{11} R_{1} G^{Mix\ *} \varphi \|^{2}_{U_{1}} + \| W_{22}^{*^{-1}} R_{2} G^{Mix\ *} \varphi \|^{2}_{U_{2}}  \right\}
\nonumber\\
&&
+(\left(\begin{array}{cc}
      0 & C_{12} \\
      C_{21} & 0 
    \end{array}\right) G^{Mix \ *} \varphi, G^{Mix \ *} \varphi)_{U_{1}\times U_{2}}
\nonumber\\
&\geq& c_{1} \| R_{1} G^{Mix\ *} \varphi \|^{2}_{X_{1}^{*}} - c_{2} \| R_{2} G^{Mix\ *} \varphi \|^{2}_{X_{2}}
\nonumber\\
&& - c_{3} \| R_{1} G^{Mix\ *} \varphi \|_{X_{1}^{*}} \| R_{2} G^{Mix\ *} \varphi \|_{X_{2}},
\label{(3)-2}
\end{eqnarray}
for all $\varphi \in V_{2}^{\bot}$. Setting $V:=V_{1}\cup V_{2}$, $V$ is finite dimensional subspace in $Y$. Then by (\ref{(3)-1}) and (\ref{(3)-2}) we have 
\begin{eqnarray}
c_{1} \| R_{1} G^{Mix \ *} \varphi \|_{\tilde{X}^{*}_{1}}^{2} &\leq& c_{2} \| R_{2} G^{Mix\ *} \varphi \|^{2}_{X_{2}}
\nonumber\\
&+& c_{3} \| R_{1} G^{Mix\ *} \varphi \|_{X_{1}^{*}} \| R_{2} G^{Mix\ *} \varphi \|_{X_{2}} + \| \mathrm{Re}\tilde{T} \|  \| \tilde{G}^{*} \varphi \|_{\tilde{X}^{*}}^{2},
\nonumber \\ 
 \label{(3)-3}
\end{eqnarray}
for all $\varphi \in V^{\bot}$.
\par
On the other hand, by the assumption {\bf (3b)} and Lemma \ref{Vector space lemma}, we have $W_1 \not \subset \mathrm{Ran}(\left[\tilde{G}, G^{Mix}R_{2}^{*} \right]) + V$, which implies that by $W_1 \subset \mathrm{Ran}(G^{Mix}R_{1}^{*})$
\begin{equation}
\mathrm{Ran}(G^{Mix}R_{1}^{*}) \not \subset \mathrm{Ran}\left( [\tilde{G}, G^{Mix}R_{2}^{*}] \right) + V = \mathrm{Ran}\left([\tilde{G}, G^{Mix}R_{2}^{*}, P_{V}] \right).
\end{equation}
By this and Lemma \ref{Range inclusion}, there exists sequence $\displaystyle (\varphi_{n})_{n \in \mathbb{N} } \subset Y$ such that 
\begin{equation}
\| R_{1} G^{Mix \ *} \varphi_{n} \|_{X^{*}_{1}}^{2} > n^{2} \left( \| \tilde{G}^{*} \varphi_{n} \|_{\tilde{X}^{*}}^{2} + \| R_{2} G^{Mix \ *} \varphi_{n} \|_{X_{2}}^{2}  + \| P_{V} \varphi_{n} \|_{Y}^{2} \right),
\end{equation}
for all $n \in \mathbb{N}$. Setting $\psi_{n}:=\displaystyle \frac{ n^{\frac{1}{2}} \varphi_n}{\| R_{1} G^{Mix \ *} \varphi_{n} \|_{X^{*}_{1}}}$, we obtain
\begin{equation}
 \| \tilde{G}^{*} \psi_{n} \|_{\tilde{X}^{*}}^{2} + \| R_{2} G^{Mix \ *} \psi_{n} \|_{X_{2}}^{2}  + \| P_{V} \psi_{n} \|_{Y}^{2} \leq \frac{1}{n}, \ \ \ \  \| R_{1} G^{Mix \ *} \varphi_{n} \|_{X^{*}_{1}}^{2} = n.
\end{equation}
Setting $\tilde{\psi}_{n}:=(I-P_{V})\psi_{n} \in V^{\bot}$, we finally obtain as $\ n \to \infty$
\begin{equation}
\| R_{1}G^{Mix \ *} \tilde{\psi}_{n} \|_{X^{*}_{1}} \geq \| R_{1}G^{Mix \ *} \psi_{n} \|_{X^{*}_{1}} - \|R_{1}G^{Mix \ *} \| \| P_{V} \psi_{n} \|_{Y},
\end{equation}
\begin{equation}
\hspace{-2cm} \| \tilde{G}^{*} \tilde{\psi}_{n} \|_{\tilde{X}^{*}} + \| R_{2}G^{Mix \ *} \tilde{\psi}_{n} \|_{X_{2}} \leq \| \tilde{G}^{*} \psi_{n} \|_{\tilde{X}^{*}} + \|\tilde{G}^{*} \| \| P_{V} \psi_{n} \|_{Y}
\nonumber 
\end{equation}
\begin{equation}
\hspace{2cm} + \| R_{2}G^{Mix \ *} \psi_{n} \|_{X_{2}} + \|\tilde{G}^{*} \| \| P_{V} \psi_{n} \|_{Y} \to 0,
\end{equation}
\begin{eqnarray}
&&\| R_{1}G^{Mix \ *} \tilde{\psi}_{n} \|_{X^{*}_{1}}\| R_{2}G^{Mix \ *} \tilde{\psi}_{n} \|_{X_{2}}
\nonumber\\
&\leq& \left( \| R_{1}G^{Mix \ *} \psi_{n} \| + \|R_{1}G^{Mix \ *} \| \| P_{V} \psi_{n} \| \right) 
\nonumber \\
&&
\times \left( \| R_{2}G^{Mix \ *} \psi_{n} \| + \|R_{2}G^{Mix \ *} \| \| P_{V} \psi_{n} \| \right) \to  1 +  \|R_{2}G^{Mix \ *} \|, 
\nonumber \\ 
\end{eqnarray}
which contradicts (\ref{(3)-3}). Therefore, we conclude (\ref{conclusion(3)}).

\end{proof}
\section{Applications of the general theorem}
In the following, we study many applications of Theorem \ref{General theorem for MM} to inverse acoustic scatterings. 
\subsection{Dirichlet obstacle}
Let $F^{Dir}_{\Omega}$ be the far-field operator for a Dirichlet obstacle $\Omega$, that is, $F^{Dir}_{\Omega}$ is the far-field operator defined by (\ref{FFO}) corresponding to the solution of (\ref{Obstacle})--(\ref{BC}) where $\mathcal{B}u=u$. $F^{Dir}_{\Omega}$ has the factorization (see Theorem 1.15 of \cite{Kirsch and Grinberg})
\begin{equation}
F^{Dir}_{\Omega}=-G^{Dir}_{\Omega}S^{*}_{\Omega}G^{Dir\ *}_{\Omega},
\end{equation}
where $G^{Dir}_{\Omega}:H^{1/2}(\partial \Omega) \to L^{2}(\mathbb{S}^{d-1})$ is the data-to-pattern operator defined by $G^{Dir}_{\Omega}f:=v^{\infty}$ where $v^{\infty}$ is the far-field pattern of a radiating solution $v$ (that is, $v$ satisfies the Sommerfeld radiation condition (\ref{SRC})) such that 
\begin{equation}
\Delta v+k^2v=0 \ \mathrm{in} \ \mathbb{R}^d \setminus \overline{\Omega}, \ \ \ \ v=f \ \mathrm{on} \ \partial \Omega, \label{definition of G}
\end{equation}
and $S_{\Omega}:H^{-1/2}(\partial \Omega) \to H^{1/2}(\partial \Omega)$ is the single layer boundary operator defined by
\begin{equation}
S_{\Omega}\varphi(x):=\int_{\partial \Omega} \varphi(y)\Phi(x,y)ds(y), \ x \in \partial \Omega,
\end{equation}
where $\Phi(x,y)$ denotes the fundamental solution for the Helmholtz equation in $\mathbb{R}^d$, i.e., 
\begin{equation}
\Phi(x,y):=\left\{ \begin{array}{ll}
\displaystyle \frac{i}{4}H^{(1)}_0(k|x-y|), & \quad \mbox{$d=2$},  \\
\displaystyle \frac{e^{ik|x-y|}}{4\pi |x-y|}, & \quad \mbox{$d=3$}, \\
\end{array}\right.
\end{equation}
where $H^{(1)}_0$ is the Hankel function of the first kind of order one. The single layer boundary operator $S_{\Omega}$ is of the form (see Lemma 1.14 of \cite{Kirsch and Grinberg})
\begin{equation}
S_{\Omega}=C_{\Omega}+K_{\Omega},
\end{equation}
where $C_{\Omega}:H^{-1/2}(\partial \Omega) \to H^{1/2}(\partial \Omega)$ is some positive coercive operator and $K_{\Omega}:H^{-1/2}(\partial \Omega) \to H^{1/2}(\partial \Omega)$ is some compact operator. 
\par
For a bounded domain $B \subset \mathbb{R}^{d}$ with the smooth boundary, we define the Herglotz operator $H_{\partial B}:L^{2}(\mathbb{S}^{d-1}) \to L^{2}(\partial B)$ by
\begin{equation}
H_{\partial B}g(x):=\int_{\mathbb{S}^{1}}\mathrm{e}^{ik\theta \cdot x}g(\theta)ds(\theta), \ x \in \partial B, 
\end{equation}
If the range is restricted to the space $H^{1/2}(\partial B)$, we denote its Herglotz operator by $\hat{H}_{\partial B}:L^{2}(\mathbb{S}^{d-1}) \to H^{1/2}(\partial B)$. From the definition, we have $\hat{H}_{\partial B}^{*}= G^{Dir}_{B}S_{B}$ (see e.g, Theorem 1.15 of \cite{Kirsch and Grinberg}). Let $J_{\partial B}:H^{1/2}(\partial B) \to H^{-1/2}(\partial B)$ be a self-adjoint compact embedding. Then, we have
\begin{eqnarray}
H^{*}_{\partial B}H_{\partial B} &=& \hat{H}^{*}_{\partial B}J_{\partial B}\hat{H}_{\partial B}= G^{Dir}_{B}S_{B}J_{\partial B} S_{B}^{*} G^{Dir\ *}_{B}
\nonumber\\
&=&G^{Dir}_{B}\left[C_{B}J_{\partial B}C_{B} + \hat{K}_{B} \right] G^{Dir\ *}_{B}, \label{HBHB}
\end{eqnarray}
where $C_{B}J_{\partial B}C_{B}$ is a positive coercive operator and $\hat{K}_{B}$ is some self-adjoint compact operator. 
\par
Assume that $B \subset \Omega$. Then, we can define $R:H^{1/2}(\partial B) \to H^{1/2}(\partial \Omega)$ by $Rf:=v\bigl|_{\partial \Omega}$ where $v$ is a radiating solution $v$ of (\ref{definition of G}) replacing $\Omega$ with $B$. Since $v\bigl|_{\partial \Omega} \in C^{\infty}(\partial \Omega)$, $R$ is a compact operator, and by the definition we have 
\begin{equation}
G^{Dir}_{B}=G^{Dir}_{\Omega}R, \label{(4.1)-1}
\end{equation}
which corresponds to the assumption {\bf(1a)} of Theorem \ref{General theorem for MM}. Applying {\bf (1)} of Theorem \ref{General theorem for MM} as
\begin{equation*}
F=-F^{Dir}_{\Omega}=G^{Dir}_{\Omega}(C_{\Omega}+K^{*}_{\Omega})G^{Dir\ *}_{\Omega},
\end{equation*}
\begin{equation*}
\tilde{F}=H^{*}_{\partial B}H_{\partial B}=G^{Dir}_{B}\left[C_{B}J_{\partial B}C_{B} + \hat{K}_{B} \right] G^{Dir\ *}_{B},
\end{equation*}
we have
\begin{equation}
H^{*}_{\partial B}H_{\partial B} \leq_{\mathrm{fin}} -\mathrm{Re}F^{Dir}_{\Omega}.
\end{equation}
\par
Assume that $B \not \subset \Omega$. Then, there exists a bounded domain $B_0 \Subset B$ such that $B_{0} \cap \Omega = \emptyset$. We set $W:=\mathrm{Ran}(G^{Dir}_{B_0}) \subset \mathrm{Ran}(G^{Dir}_{B})$, then, $W$ is an infinite dimensional subspace of $L^{2}(\mathbb{S}^{d-1})$ because $G^{Dir}_{B_0}$ is injective (see e.g., Lemma 1.13 of \cite{Kirsch and Grinberg}). From $B_{0} \cap \Omega = \emptyset$, we obtain 
\begin{equation}
W \cap \mathrm{Ran}(G^{Dir}_{\Omega})=\{ 0 \},
\end{equation}
(see e.g., Lemma 4.2 of \cite{Albicker and Griesmaier}) which corresponds to the assumption {\bf(2b)} of Theorem \ref{General theorem for MM}. Applying {\bf(2)} of Theorem \ref{General theorem for MM} as 
\begin{equation*}
F=-F^{Dir}_{\Omega}=G^{Dir}_{\Omega}(C_{\Omega}+K^{*}_{\Omega})G^{Dir\ *}_{\Omega},
\end{equation*}
\begin{equation*}
\tilde{F}=H^{*}_{\partial B}H_{\partial B}=G^{Dir}_{B}\left[C_{B}J_{\partial B}C_{B} + \hat{K}_{B} \right] G^{Dir\ *}_{B},
\end{equation*}
we have
\begin{equation}
H^{*}_{\partial B}H_{\partial B} \not \leq_{\mathrm{fin}} -\mathrm{Re}F^{Dir}_{\Omega}.
\end{equation}
From the above discussion, we conclude the following theorem, which is the same result as Theorem 5.3 of \cite{Albicker and Griesmaier}.
\begin{thm}[Theorem 5.3 of \cite{Albicker and Griesmaier}]\label{Dirichlet obstacle inside}
Let $B \subset \mathbb{R}^d$ be a bounded domain with the smooth boundary. Then,
\begin{equation}
B \subset \Omega \ \ \ \  \Longleftrightarrow \ \ \ \  H^{*}_{\partial B}H_{\partial B} \leq_{\mathrm{fin}} -\mathrm{Re}F^{Dir}_{\Omega},
\end{equation}
\end{thm}
By the same argument in Theorem \ref{Dirichlet obstacle inside}, one can  apply {\bf(1)} and {\bf(2)} of Theorem \ref{General theorem for MM} as 
\begin{equation*}
F=H^{*}_{\partial B}H_{\partial B}=G^{Dir}_{B}\left[C_{B}J_{\partial B}C_{B} + \hat{K}_{B} \right] G^{Dir\ *}_{B},
\end{equation*}
\begin{equation*}
\tilde{F}=-F^{Dir}_{\Omega}=G^{Dir}_{\Omega}(C_{\Omega}+K^{*}_{\Omega})G^{Dir\ *}_{\Omega}.
\end{equation*}
Then, we also conclude the following theorem.
\begin{thm}\label{Dirichlet obstacle outside}
Let $B \subset \mathbb{R}^d$ be a bounded domain with the smooth boundary. Then,
\begin{equation}
\Omega  \subset B \ \ \ \  \Longleftrightarrow \ \ \ \  -\mathrm{Re}F^{Dir}_{\Omega}  \leq_{\mathrm{fin}} H^{*}_{\partial B}H_{\partial B},
\end{equation}
\end{thm}
\begin{rem}
We remark that the factorization reconstruction for the inverse obstacle scattering needs to assume that $k^2$ is not a Dirichlet eigenvalue of $-\Delta$ in $\Omega$ (see e.g., Corollary 2.16 of \cite{Kirsch and Grinberg}). While, the monotonicity reconstruction (Theorems \ref{Dirichlet obstacle inside} and \ref{Dirichlet obstacle outside}) does not require the assumption of Dirichlet eigenvalues.
\end{rem}
\subsection{Inhomogeneous medium}
Let $F^{Med}_{\Omega}$ be the far-field operator for an inhomogeneous medium $\Omega$ with the function $q \in L^{\infty}(\Omega)$, that is, $F^{Med}_{\Omega}$ is the far-field operator defined by (\ref{FFO}) corresponding to the solution of (\ref{Medium})--(\ref{SRC}). Throughout this section, we assume that there exists a constant $q_0>0$ such that $q \geq q_0$ in $\Omega$. $F^{Med}_{\Omega}$ has the following factorization by the same argument in Theorem 4.5 of \cite{Kirsch and Grinberg}
\begin{equation}
F^{Med}_{\Omega}=H^{*}_{\Omega}T_{\Omega}H_{\Omega},
\end{equation}
where $H_{\Omega}:L^{2}(\mathbb{S}^{d-1}) \to L^{2}(\Omega)$ is the Herglotz operator defined by
\begin{equation}
H_{\Omega}g(x):=\int_{\mathbb{S}^{1}}\mathrm{e}^{ik\theta \cdot x}g(\theta)ds(\theta), \ x \in \Omega, 
\end{equation}
and some operator $T_{\Omega}:L^{2}(\Omega) \to L^{2}(\Omega)$ is of the form 
\begin{equation}
T_{\Omega}=k^{2}qI_{L^{2}(\Omega)} + K_{\Omega},
\end{equation}
where $K_{\Omega}$ is some compact operator. 
\par
Let $B \subset \mathbb{R}^{d}$ be a bounded domain with the smooth boundary. Assume that $B \subset \Omega$. We define the restriction operator $R:L^{2}(\Omega) \to L^{2}(B)$ by $Rf:=f\bigl|_{B}$. Then by the definition we have $H_{B}=RH_{\Omega}$, and for $\alpha \in (0, k^2 q_{0})$ we have
\begin{equation}
F^{Med}_{\Omega} - \alpha H_{B}^{*} H_{B} = H^{*}_{\Omega}\left[k^{2}qI_{L^{2}(\Omega)} - \alpha R^{*}R + K_{\Omega} \right]H_{\Omega},
\end{equation}
where the operator $k^{2}qI_{L^{2}(\Omega)} - \alpha R^{*}R$ is positive coercive when $\alpha \in (0, k^2 q_{0})$. Applying {\bf(1)} of Theorem \ref{General theorem for MM} as
\begin{equation*}
F=F^{Med}_{\Omega} - \alpha H_{B}^{*} H_{B} = H^{*}_{\Omega}\left[k^{2}qI_{L^{2}(\Omega)} - \alpha R^{*}R + K_{\Omega} \right]H_{\Omega},
\end{equation*}
\begin{equation*}
\tilde{F}=0,
\end{equation*}
we have for $\alpha \in (0, k^2 q_{0})$
\begin{equation}
\alpha H^{*}_{B}H_{B} \leq_{\mathrm{fin}} \mathrm{Re}F^{Med}_{\Omega}.
\end{equation}
\par
Assume that $B \not \subset \Omega$. Then, there exists a bounded domain $B_0 \Subset B$ such that $B_{0} \cap \Omega = \emptyset$. We set $W:=\mathrm{Ran}(H^{*}_{B_0}) \subset \mathrm{Ran}(H^{*}_{B})$, then, $W$ is an infinite dimensional subspace of $L^{2}(\mathbb{S}^1)$ because $H^{*}_{B_0}$ is injective. From $B_{0} \cap \Omega = \emptyset$, we obtain 
\begin{equation}
W \cap \mathrm{Ran}(H^{*}_{\Omega})=\{ 0 \}.
\end{equation}
(see e.g., Lemma 4.3 of \cite{Griesmaier and Harrach}). Applying {\bf(2)} of Theorem \ref{General theorem for MM} as
\begin{equation*}
F=F^{Med}_{\Omega}=H^{*}_{\Omega}T_{\Omega}H_{\Omega},
\end{equation*}
\begin{equation*}
\tilde{F}=\alpha H^{*}_{B}H_{B},
\end{equation*}
we have for $\alpha \in (0, k^2 q_{0})$,
\begin{equation}
\alpha H^{*}_{B}H_{B} \not \leq_{\mathrm{fin}} \mathrm{Re}F^{Med}_{\Omega}.
\end{equation}
From the above discussion, we conclude the following theorem, which is the same result as Theorem 5.1 of \cite{Griesmaier and Harrach}.
\begin{thm}[Theorem 5.1 of \cite{Griesmaier and Harrach}]\label{Medium inside}
Let $B \subset \mathbb{R}^d$ be a bounded domain with the smooth boundary. Then, for $\alpha \in (0, k^2 q_{0})$
\begin{equation}
B \subset \Omega \ \ \ \  \Longleftrightarrow \ \ \ \  \alpha H^{*}_{B}H_{B} \leq_{\mathrm{fin}} \mathrm{Re}F^{Med}_{\Omega}.
\end{equation}
\end{thm}
\par
Assume that $\Omega \subset B$. Then, we can define the compact operator $R:L^{2}(\Omega) \to H^{1/2}(\partial B)$ by $Rg:=w\bigl|_{\partial B}$ where $w$ is a radiating solution $w$ of 
\begin{equation}
\Delta w+k^2(1+q)w=-k^2qg \ \mathrm{in} \ \mathbb{R}^d, \label{definition of medium G}
\end{equation}
From the definition we obtain $G^{Med}_{\Omega}=G^{Dir}_{B}R$ where the data-to-pattern operator $G^{Med}_{\Omega}$ is defined by $G^{Med}_{\Omega}g:=w^{\infty}$ and $G^{Dir}_{B}$ is defined by (\ref{definition of G}) replacing $\Omega$ by $B$. Since $G^{Med}_{\Omega}=H_{\Omega}^{*}T_{\Omega}$ and $T_{\Omega}$ is bounded invertible (see e.g., the arguments of Theorem 4.5 of \cite{Kirsch and Grinberg}), we have 
\begin{equation}
H_{\Omega}^{*}=G^{Med}_{\Omega}T^{-1}_{\Omega}=G^{Dir}_{B}RT^{-1}_{\Omega}.
\end{equation}
Applying {\bf (1)} of Theorem \ref{General theorem for MM} as 
\begin{equation*}
F=H^{*}_{\partial B}H_{\partial B}=G^{Dir}_{B}\left[C_{B}J_{\partial B}C_{B} + \hat{K}_{\Omega} \right] G^{Dir\ *}_{B},
\end{equation*}
\begin{equation*}
\tilde{F}=F^{Med}_{\Omega}=H^{*}_{\Omega}T_{\Omega}H_{\Omega}=H^{*}_{\Omega}\left[k^{2}qI_{L^{2}(\Omega)} + K_{\Omega}\right]H_{\Omega},
\end{equation*}
we have
\begin{equation}
\mathrm{Re}F^{Med}_{\Omega} \leq_{\mathrm{fin}} H^{*}_{\partial B}H_{\partial B}.
\end{equation}
\par
Assume that $\Omega \not \subset B$. Then, there exists a bounded domain $\Omega_0 \Subset \Omega$ such that $\Omega_{0} \cap B = \emptyset$. We set $W:=\mathrm{Ran}(G^{Med}_{\Omega_0}) \subset \mathrm{Ran}(G^{Med}_{\Omega})=\mathrm{Ran}(H^{*}_{\Omega})$, then, $W$ is an infinite dimensional subspace of $L^{2}(\mathbb{S}^1)$ because $G^{Med}_{\Omega_0}$ is injective. From $\Omega_{0} \cap B = \emptyset$, we obtain 
\begin{equation}
W \cap \mathrm{Ran}(G^{Dir}_{B})=\{ 0 \},
\end{equation}
(see e.g., Lemma 4.3 of \cite{Griesmaier and Harrach}). Applying {\bf (2)} of Theorem \ref{General theorem for MM} as 
\begin{equation*}
F=H^{*}_{\partial B}H_{\partial B}=G^{Dir}_{B}\left[C_{B}J_{\partial B}C_{B} + \hat{K}_{\Omega} \right] G^{Dir\ *}_{B},
\end{equation*}
\begin{equation*}
\tilde{F}=F^{Med}_{\Omega}=H^{*}_{\Omega}T_{\Omega}H_{\Omega}=H^{*}_{\Omega}\left[k^{2}qI_{L^{2}(\Omega)} + K_{\Omega}\right]H_{\Omega},
\end{equation*}
we have
\begin{equation}
\mathrm{Re}F^{Med}_{\Omega} \not\leq_{\mathrm{fin}} H^{*}_{\partial B}H_{\partial B}.
\end{equation}
From the above discussion, we conclude the following theorem.
\begin{thm}\label{Medium outside}
Let $B \subset \mathbb{R}^d$ be a bounded domain with the smooth boundary. Then,
\begin{equation}
\Omega \subset B  \ \ \ \  \Longleftrightarrow \ \ \ \  \mathrm{Re}F^{Med}_{\Omega}  \leq_{\mathrm{fin}} H^{*}_{\partial B}H_{\partial B},
\end{equation}
\end{thm}
\begin{rem}
We remark that the factorization reconstruction for the inverse medium scattering needs to assume that $k^2$ is not a transmission eigenvalue in $\Omega$ (see e.g., Theorem 4.10 of \cite{Kirsch and Grinberg}). While, the monotonicity reconstruction (Theorems \ref{Medium inside} and \ref{Medium outside}) does not require the assumption of transmission eigenvalues.
\end{rem}
\subsection{Dirichlet crack}
\begin{figure}[h]
  \centering
  \includegraphics[scale=0.3]{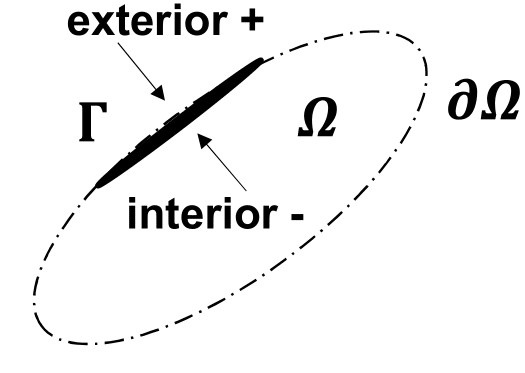}
  \caption{The assumption for $\Gamma$ and the boundary condition.}\label{crack BC}
\end{figure}

Let $F^{Dir}_{\Gamma}$ be the far-field operator for a Dirichlet crack $\Gamma$ where $\Gamma \subset \mathbb{R}^d$ is a smooth non-intersecting open arc ($d=2$) or surface ($d=3$), and we assume that $\Gamma$ can be extended to some smooth, connected, closed curve ($d=2$) or surface ($d=3$) $\partial \Omega$ enclosing a bounded domain $\Omega$ in $\mathbb{R}^d$. The corresponding far-field pattern is defined by solving the scattering problem (\ref{Obstacle})--(\ref{BC}) where $\Omega$ in (\ref{Obstacle}) is replaced by $\Gamma$ and the boundary condition (\ref{BC}) is replaced by
\begin{equation}
u_{-}=0 \ \mathrm{on} \ \Gamma, \ \ \ \ \ \ \ \ u_{+}=0 \ \mathrm{on} \ \Gamma,
\end{equation}
where we denote by $u_{\pm}$ the limit of $u$ approaching the boundary from exterior (+) and interior (-) of an extensive domain $\Omega$ (see Figure \ref{crack BC}). $F^{Dir}_{\Gamma}$ has the factorization (see Lemma 3.4 of \cite{Kirsch and Ritter})
\begin{equation}
F^{Dir}_{\Gamma}=-G^{Dir}_{\Gamma}S^{*}_{\Gamma}G^{Dir\ *}_{\Gamma},
\end{equation}
where $G^{Dir}_{\Gamma}:H^{1/2}(\Gamma) \to L^{2}(\mathbb{S}^{d-1})$ is the data-to-pattern operator corresponding to the crack $\Gamma$,
and $S_{\Gamma}:\tilde{H}^{-1/2}(\Gamma) \to H^{1/2}(\Gamma)$ is the single layer boundary operator corresponding to the crack $\Gamma$ where we denote by
\begin{equation}
H^{1/2}(\Gamma):= \{u\bigl|_{\Gamma} : u \in H^{1/2}(\partial\Omega) \},
\end{equation}
\begin{equation}
\tilde{H}^{1/2}(\Gamma):= \{u\bigl|_{\Gamma}: u \in H^{1/2}(\partial \Omega),\ supp(u) \subset  \overline{\Gamma} \}, 
\end{equation}
and $H^{-1/2}(\Gamma)$ and $\tilde{H}^{-1/2}(\Gamma)$ the dual spaces of $\tilde{H}^{1/2}(\Gamma)$ and $H^{1/2}(\Gamma)$, respectively. We have the following inclusion relation
\begin{equation}
\tilde{H}^{1/2}(\Gamma) \subset H^{1/2}(\Gamma) \subset L^{2}(\Gamma) \subset \tilde{H}^{-1/2}(\Gamma)\subset H^{-1/2}(\Gamma).
\end{equation}
The single layer boundary operator $S_{\Gamma}$ is of the form (see Lemma 3.2 of \cite{Kirsch and Ritter})
\begin{equation}
S_{\Gamma}=C_{\Gamma}+K_{\Gamma},
\end{equation}
where $C_{\Gamma}:\tilde{H}^{-1/2}(\Gamma) \to H^{1/2}(\Gamma)$ is some positive coercive operator and $K_{\Gamma}:\tilde{H}^{-1/2}(\Gamma) \to H^{1/2}(\Gamma)$ is some compact operator. 
\par
Let $\sigma \subset \mathbb{R}^d$ be a smooth arc ($d=2$) or surface ($d=3$). Assume that $\sigma \subset \Gamma$. Then, we define $R:L^{2}(\Gamma) \to L^{2}(\sigma)$ by $Rf:=f\bigl|_{\sigma}$. Let $J:H^{1/2}(\Gamma) \to L^{2}(\Gamma)$ be the compact embedding. Since $\hat{H}_{\Gamma}^{*}=G^{Dir}_{\Gamma}S_{\Gamma}$ (see e.g., Lemma 3.4 of \cite{Kirsch and Ritter}), we have 
\begin{equation}
H_{\sigma}=RH_{\Gamma}=RJ\hat{H}_{\Gamma}=RJS_{\Gamma}^{*}G^{Dir\ *}_{\Gamma},
\end{equation}
where $H_{\Gamma}:L^{2}(\mathbb{S}^{d-1}) \to L^{2}(\Gamma)$ is the Herglotz operator corresponding to $\Gamma$ and $\hat{H}_{\Gamma}:L^{2}(\mathbb{S}^{d-1}) \to H^{1/2}(\Gamma)$ is the Herglotz operator that its range is restricted to the space $H^{1/2}(\Gamma)$. Applying {\bf (1)} of Theorem \ref{General theorem for MM} as 
\begin{equation*}
F=-F^{Dir}_{\Gamma}=G^{Dir}_{\Gamma}(C_{\Gamma}+K^{*}_{\Gamma})G^{Dir\ *}_{\Gamma},
\end{equation*}
\begin{equation*}
\tilde{F}=H^{*}_{\sigma}H_{\sigma},
\end{equation*}
we have
\begin{equation}
H^{*}_{\sigma}H_{\sigma} \leq_{\mathrm{fin}} -\mathrm{Re}F^{Dir}_{\Gamma}.
\end{equation}
\par
Assume that $\sigma \not \subset \Gamma$. Then, there exists $\sigma_0 \Subset \sigma$ such that $\sigma_{0} \cap \Gamma = \emptyset$. We set $W:=\mathrm{Ran}(H^{*}_{\sigma_0}) \subset \mathrm{Ran}(H^{*}_{\sigma})$, then, $W$ is an infinite dimensional subspace of $L^{2}(\mathbb{S}^1)$ because $H^{*}_{\sigma_0}$ is injective. From $\sigma_{0} \cap \Gamma = \emptyset$, we obtain 
\begin{equation}
W \cap \mathrm{Ran}(G^{Dir}_{\Gamma})=\{ 0 \}.
\end{equation}
(see e.g., Lemma 4.1 of \cite{Furuya3}). Applying {\bf(2)} of Theorem \ref{General theorem for MM} as 
\begin{equation*}
F=-F^{Dir}_{\Gamma}=G^{Dir}_{\Gamma}(C_{\Gamma}+K^{*}_{\Gamma})G^{Dir\ *}_{\Gamma},
\end{equation*}
\begin{equation*}
\tilde{F}=H^{*}_{\sigma}H_{\sigma},
\end{equation*}
we have
\begin{equation}
H^{*}_{\sigma}H_{\sigma} \not \leq_{\mathrm{fin}} -\mathrm{Re}F^{Dir}_{\Gamma}.
\end{equation}
From the above discussion, we conclude the following theorem, which is the same result as Theorem 1.1 of \cite{Furuya3}.
\begin{thm}[Theorem 1.1 of \cite{Furuya3}]\label{Dirichlet crack inside}
Let $\sigma \subset \mathbb{R}^d$ be a smooth arc ($d=2$) or surface ($d=3$). Then,
\begin{equation}
\sigma \subset \Gamma \ \ \ \  \Longleftrightarrow \ \ \ \  H^{*}_{\sigma}H_{\sigma} \leq_{\mathrm{fin}} -\mathrm{Re}F^{Dir}_{\Gamma}.
\end{equation}
\end{thm}
By the same argument in Theorem \ref{Dirichlet crack inside}, one can  apply {\bf(1)} and {\bf(2)} of Theorem \ref{General theorem for MM} as 
\begin{equation*}
F=H^{*}_{\partial B}H_{\partial B}=G^{Dir}_{B}\left[C_{B}J_{\partial B}C_{B} + \hat{K}_{B} \right] G^{Dir\ *}_{B},
\end{equation*}
\begin{equation*}
\tilde{F}=-F^{Dir}_{\Gamma}=G^{Dir}_{\Gamma}(C_{\Gamma}+K^{*}_{\Gamma})G^{Dir\ *}_{\Gamma}.
\end{equation*}
Then, we also conclude the following theorem, which is the same result as Theorem 1.2 of \cite{Furuya3}.
\begin{thm}[Theorem 1.2 of \cite{Furuya3}]\label{Dirichlet crack outside}
Let $B \subset \mathbb{R}^d$ be a bounded domain with the smooth boundary. Then,
\begin{equation}
\Gamma \subset B \ \ \ \  \Longleftrightarrow \ \ \ \  -\mathrm{Re}F^{Dir}_{\Gamma}  \leq_{\mathrm{fin}} H^{*}_{\partial B}H_{\partial B}.
\end{equation}
\end{thm}
\begin{rem}
We remark that the factorization reconstruction for the inverse crack scattering also does not restrict the wave number. (see Theorem 3.9 of \cite{Kirsch and Ritter}). One of the advantage of the monotonicity over the factorization is to have not only inside tests (Theorem \ref{Dirichlet crack inside}), but also outside tests (Theorem \ref{Dirichlet crack outside}). 
\end{rem}
\subsection{Mixed obstacle}
Let $F^{Mix}_{\Omega_{1}, \Omega_{2}}$ be the far-field operator for the mixed obstacle $\Omega=\Omega_{1} \cup \Omega_{2}$ with the Dirichlet part $\Omega_{1}$ and the Neumann part $\Omega_{2}$ where $\Omega_1, \Omega_2$ are bounded domains with the smooth boundary such that $\Omega_1 \cup \Omega_2 = \emptyset$. The corresponding far-field pattern is defined by solving the scattering problem (\ref{Obstacle})--(\ref{BC}) where the boundary condition (\ref{BC}) is replaced by
\begin{equation}
u=0 \ \mathrm{on} \ \partial \Omega_{1}, \ \ \ \ \ \ \ \ \frac{\partial u}{\partial \nu}=0 \ \mathrm{on} \ \partial \Omega_{2}.
\end{equation}
$F^{Mix}_{\Omega_{1}, \Omega_{2}}$ has the factorization (see Theorem 3.4 of \cite{Kirsch and Grinberg})
\begin{equation}
F^{Mix}_{\Omega_{1}, \Omega_{2}}=-G^{Mix}_{\Omega_{1}, \Omega_{2}}T^{Mix}_{\Omega_{1}, \Omega_{2}}G^{Mix\ *}_{\Omega_{1}, \Omega_{2}},
\end{equation}
where $G^{Mix}_{\Omega_{1}, \Omega_{2}}:H^{1/2}(\partial \Omega_{1}) \times H^{-1/2}(\partial \Omega_{2}) \to L^{2}(\mathbb{S}^{d-1})$ is the data-to-pattern operator for the mixed obstacle $\Omega$, i.e., defined by $G^{Mix}_{\Omega_{1}, \Omega_{2}}\left(\begin{array}{cc}
      f_{1} \\
      f_{2} 
    \end{array}\right):=v^{\infty}$ where $v^{\infty}$ is the far-field pattern of a radiating solution $v$ such that 
\begin{equation}
\Delta v+k^2v=0 \ \mathrm{in} \ \mathbb{R}^d \setminus \overline{\Omega},\ \  v=f_{1} \ \mathrm{on} \ \partial \Omega_{1}, \ \  \frac{\partial v}{\partial \nu}=f_{2} \ \mathrm{on} \ \partial \Omega_{2},
\end{equation}
and some operator $T^{Mix}_{\Omega_{1}, \Omega_{2}}:H^{-1/2}(\partial \Omega_{1}) \times H^{1/2}(\partial \Omega_{2}) \to H^{1/2}(\partial \Omega_{1}) \times H^{-1/2}(\partial \Omega_{2})$ has the form (Theorem 3.4 of \cite{Kirsch and Grinberg})
\begin{equation}
T^{Mix}_{\Omega_{1}, \Omega_{2}}=\left(\begin{array}{cc}
      C^{+}_{\Omega_{1}} &0 \\
      0 &C^{-}_{\Omega_{2}} 
    \end{array}\right) + K^{Mix}_{\Omega_{1}, \Omega_{2}},
\end{equation}
where $C^{+}_{\Omega_{1}}:H^{-1/2}(\partial \Omega_{1}) \to H^{1/2}(\partial \Omega_{1})$ is some positive coercive operator, $C^{-}_{\Omega_{2}}:H^{1/2}(\partial \Omega_{2}) \to H^{-1/2}(\partial \Omega_{2})$ is some negative coercive operator, and $K^{Mix}:H^{-1/2}(\partial \Omega_{1}) \times H^{1/2}(\partial \Omega_{2}) \to H^{1/2}(\partial \Omega_{1}) \times H^{-1/2}(\partial \Omega_{2})$ is some compact operator. 
\par
Assume that $\Omega_{1} \subset B$. Then, there exists bounded domain $\tilde{B} \Subset B$ such that $\tilde{B} \cap \Omega_{2}=\emptyset$ and $\Omega_{1} \subset \tilde{B}$. Define $R:H^{1/2}(\partial \Omega_{1}) \times H^{-1/2}(\partial \Omega_{2}) \to H^{1/2}(\partial \tilde{B}) \times H^{-1/2}(\partial \Omega_{2})$ by $R\left(\begin{array}{cc}
      f_{1} \\
      f_{2} 
    \end{array}\right):=\left(\begin{array}{cc}
      v\bigl|_{\partial \tilde{B}} \\
      f_{2} 
    \end{array}\right)$. Then, $R$ has the form $R=\left(\begin{array}{cc}
      0 & 0 \\
      0 & I
    \end{array}\right)+\mathrm{compact}$, and 
\begin{equation}    
G^{Mix}_{\Omega_{1}, \Omega_{2}} = G^{Mix}_{\tilde{B}, \Omega_{2}}R, \label{(4.4)-1}
\end{equation}
where $ G^{Mix}_{\tilde{B}, \Omega_{2}}:H^{1/2}(\partial \tilde{B}) \times H^{-1/2}(\partial \Omega_{2}) \to L^{2}(\mathbb{S}^{d-1})$ is the data-to-pattern operator corresponding to $\tilde{B} \cup \Omega_{2}$ with the Dirichlet part $\tilde{B}$ and the Neumann part $\Omega_2$. We also define $\tilde{R}:H^{1/2}(\partial \tilde{B}) \to H^{1/2}(\partial \tilde{B}) \times H^{-1/2}(\partial \Omega_{2})$ by $\tilde{R}g:=\left(\begin{array}{cc}
      g \\
      \frac{\partial w}{\partial \nu}\bigl|_{\partial \Omega_{2}}
    \end{array}\right)$ where $w$ is a radiating solution $w$ of
\begin{equation}
\Delta w+k^2w=0 \ \mathrm{in} \ \mathbb{R}^d \setminus \overline{\tilde{B}},\ \  w=g \ \mathrm{on} \ \partial \tilde{B},
\end{equation}
Then, $\tilde{R}$ has the form $\tilde{R}=\left(\begin{array}{cc}
      I \\
      0
    \end{array}\right)+\mathrm{compact}$, and 
\begin{equation}    
G^{Dir}_{\tilde{B}} = G^{Mix}_{\tilde{B}, \Omega_{2}}\tilde{R}.\label{(4.4)-2}
\end{equation}
By these and (\ref{HBHB}), we have
\begin{eqnarray}
&&F^{Mix}+ H^{*}_{\partial \tilde{B}}H_{\partial \tilde{B}}
\nonumber\\
&=& G^{Mix}_{\tilde{B}, \Omega_{2}} \left[ R(-T^{Mix}_{\Omega_{1}, \Omega_{2}})R^{*} + \tilde{R}S_{\tilde{B}}J_{\partial \tilde{B}}S^{*}_{\tilde{B}}\tilde{R}^{*}  \right]G^{Mix\ *}_{\tilde{B}, \Omega_{2}} 
\nonumber\\
&=&G^{Mix}_{\tilde{B}, \Omega_{2}}\left[ \left(\begin{array}{cc}
      C_{\tilde{B}}J_{\partial \tilde{B}}C_{\tilde{B}} &0 \\
      0 & -C_{\Omega_{2}^{-}} 
    \end{array}\right) + \tilde{K}^{Mix} \right]G^{Mix\ *}_{\tilde{B}, \Omega_{2}} ,
\end{eqnarray}
where $\left(\begin{array}{cc}
      C_{\tilde{B}}J_{\partial \tilde{B}}C_{\tilde{B}} &0 \\
      0 & -C_{\Omega_{2}^{-}} 
    \end{array}\right)$ is a positive coercive operator and $\tilde{K}^{Mix}$ is some compact operator. Applying {\bf(1)} of Theorem \ref{General theorem for MM} as 
\begin{equation*}
F=F^{Mix}_{\Omega_{1}, \Omega_{2}}+H^{*}_{\partial \tilde{B}}H_{\partial \tilde{B}}=G^{Mix}_{\tilde{B}, \Omega_{2}}\left[ \left(\begin{array}{cc}
      C_{\tilde{B}}J_{\partial \tilde{B}}C_{\tilde{B}} &0 \\
      0 & -C_{\Omega_{2}^{-}} 
    \end{array}\right) + \tilde{K}^{Mix} \right]G^{Mix\ *}_{\tilde{B}, \Omega_{2}},
\end{equation*}
\begin{equation*}
\tilde{F}=0,
\end{equation*}
we have
\begin{equation}
-\mathrm{Re}F^{Mix}_{\Omega_{1}, \Omega_{2}} \leq_{\mathrm{fin}} H^{*}_{\partial \tilde{B}}H_{\partial \tilde{B}}. \label{(4.4)-2}
\end{equation}
Since we have $\tilde{B} \subset B$, one can show  by the same argument in (\ref{(4.1)-1}) that there exists a compact operator $R_{B}:H^{1/2}(\partial \tilde{B}) \to H^{1/2}(\partial B)$ such that
\begin{equation}
G^{Dir}_{\tilde{B}}=G^{Dir}_{B}R_{B}.
\end{equation}
Then, applying {\bf (1)} of Theorem \ref{General theorem for MM} as $F=H^{*}_{\partial B}H_{\partial B}$, $\tilde{F}=H^{*}_{\partial \tilde{B}}H_{\partial \tilde{B}}$ (remark (\ref{HBHB})), it follows that
\begin{equation}
 H^{*}_{\partial \tilde{B}}H_{\partial \tilde{B}} \leq_{\mathrm{fin}} H^{*}_{\partial B}H_{\partial B},
\end{equation}
which implies that with (\ref{(4.4)-2}) we conclude that 
\begin{equation}
-\mathrm{Re}F^{Mix}_{\Omega_{1}, \Omega_{2}} \leq_{\mathrm{fin}} H^{*}_{\partial B}H_{\partial B}. \label{(4.4)-3}
\end{equation}
\par
Assume that $\Omega_{1} \not \subset B$ and $\mathbb{R}^d \setminus (\overline{B\cup \Omega})$ is connected. Then, there exists $\Gamma \Subset \partial \Omega_1$ such that $\Gamma \cap B = \emptyset$. Define $E_{\Gamma}:H^{1/2}(\Gamma) \to H^{1/2}(\partial \Omega_{1})$ by $E_{\Gamma} f=f$ on $\Gamma$, otherwise zero. Denote by $X_{\Gamma}\subset H^{1/2}(\Gamma)$ the subspace of piecewise linear continuous functions on $\Gamma$ that vanish on $\partial \Gamma$. Set $W:=\mathrm{Ran}\left(G^{Mix}R_{1}^{*}E_{\Gamma}\big|_{X_{\Gamma}}\right) \subset \mathrm{Ran}(G^{Mix}R_{1}^{*})$, it is infinite dimensional because $X_{\Gamma}$ is infinite dimensional and the operator $G^{Mix}R_{1}^{*}E_{\Gamma}$ is injective. By Lemma 4.6 of \cite{Albicker and Griesmaier}, we have 
\begin{equation}
W \cap \mathrm{Ran}(G^{Dir}_{B},G^{Mix}R_{2}^{*})=\{0 \}.
\end{equation}
Applying {\bf(3)} of Theorem \ref{General theorem for MM} as
\begin{equation*}
F^{Mix}=-F^{Mix}_{\Omega_{1}, \Omega_{2}}=G^{Mix}_{\Omega_{1}, \Omega_{2}}\left( \left(\begin{array}{cc}
      C^{+}_{\Omega_{1}} &0 \\
      0 &C^{-}_{\Omega_{2}} 
    \end{array}\right) + K^{Mix}_{\Omega_{1}, \Omega_{2}} \right)G^{Mix\ *}_{\Omega_{1}, \Omega_{2}},
\end{equation*}
\begin{equation*}
\tilde{F}=H^{*}_{\partial B}H_{\partial B}=G^{Dir}_{B}\left[C_{B}J_{\partial B}C_{B} + \hat{K}_{B} \right] G^{Dir\ *}_{B},
\end{equation*}
we have
\begin{equation}
-\mathrm{Re}F^{Mix}_{\Omega_{1}, \Omega_{2}} \not \leq_{\mathrm{fin}} H^{*}_{\partial B}H_{\partial B}.
\end{equation}
From the above discussion, we conclude the following theorem, which is the same result as Theorem 5.5 of \cite{Albicker and Griesmaier}. We remark that the monotonicity reconstruction discussed here was succeeded without assuming that the wavenumber $k^2$ is neither a Dirichlet eigenvalue of $-\Delta$ in $\Omega_1$ and $B$, nor a Neumann eigenvalue in $\Omega_2$, although Theorem 5.5 of \cite{Albicker and Griesmaier} assumed it.
\begin{thm}[Theorem 5.5 of \cite{Albicker and Griesmaier}]\label{Mixed obstacle part1}
Let $B \subset \mathbb{R}^d$ be a bounded domain with the smooth boundary such that $\mathbb{R}^d \setminus (\overline{B\cup \Omega})$ is connected. Then,
\begin{equation}
\Omega_1 \subset B \ \ \ \  \Longleftrightarrow \ \ \ \  -\mathrm{Re}F^{Mix}_{\Omega_{1}, \Omega_{2}} \leq_{\mathrm{fin}} H^{*}_{\partial B}H_{\partial B}.
\end{equation}
\end{thm}
By the same argument in Theorem \ref{Mixed obstacle part1}, one can  apply {\bf(1)} and {\bf(3)} replacing {\bf(3b)} by {\bf(3b)'} of Theorem \ref{General theorem for MM} as 
\begin{equation*}
F^{Mix}=F^{Mix}_{\Omega_{1}, \Omega_{2}}=-G^{Mix}_{\Omega_{1}, \Omega_{2}}\left( \left(\begin{array}{cc}
      C^{+}_{\Omega_{1}} &0 \\
      0 &C^{-}_{\Omega_{2}} 
    \end{array}\right) + K^{Mix}_{\Omega_{1}, \Omega_{2}} \right)G^{Dir\ *}_{\Omega_{1}, \Omega_{2}},
\end{equation*}
\begin{equation*}
\tilde{F}=H^{*}_{\partial B}H_{\partial B}=G^{Dir}_{B}\left[C_{B}J_{\partial B}C_{B} + \hat{K}_{B} \right] G^{Dir\ *}_{B},
\end{equation*}
for the reconstruction of the Neumann part $\Omega_2$. Then, we also conclude the following theorem, which is the same result as Theorem 5.5 of \cite{Albicker and Griesmaier} without the restriction of the wave number $k>0$ as well as Theorem \ref{Mixed obstacle part1}.
\begin{thm}[Theorem 5.5 of \cite{Albicker and Griesmaier}]\label{Mixed obstacle part2}
Let $B \subset \mathbb{R}^d$ be a bounded domain with the smooth boundary such that $\mathbb{R}^d \setminus (\overline{B\cup \Omega})$ is connected. Then,
\begin{equation}
\Omega_2 \subset B \ \ \ \  \Longleftrightarrow \ \ \ \  \mathrm{Re}F^{Mix}_{\Omega_{1}, \Omega_{2}} \leq_{\mathrm{fin}} H^{*}_{\partial B}H_{\partial B}.
\end{equation}
\end{thm}
\begin{rem}
We remark that the factorization reconstruction for the mixed obstacle has to assume that one component is covered by some artificial domain which is disjoint with the other one we want to reconstruct, and furthermore assume that $k^2$ is neither a Dirichlet eigenvalue of $-\Delta$ in $\Omega_1$ and artificial covering domain, nor a Neumann eigenvalue in $\Omega_2$ (see Lemma 3.5 of \cite{Kirsch and Grinberg}). However, the monotonicity reconstruction (Theorems \ref{Mixed obstacle part1} and \ref{Mixed obstacle part2}) does not require both of them.
\end{rem}
\subsection{Mixed crack}
Let $F^{Mix}_{\Gamma}$ be the far-field operator for the mixed crack $\Gamma$ (The assumption for $\Gamma$ is the same as Section 5.3). The corresponding far-field pattern is defined by solving the scattering problem (\ref{Obstacle})--(\ref{BC}) where $\Omega$ in (\ref{Obstacle}) is replaced by $\Gamma$ and the boundary condition (\ref{BC}) is replaced by
\begin{equation}
u_{-}=0 \ \mathrm{on} \ \Gamma, \ \ \ \ \ \ \ \ \frac{\partial u_{+}}{\partial \nu}=0 \ \mathrm{on} \ \Gamma.
\end{equation}
$F^{Mix}_{\Gamma}$ has the factorization (see (3.6) and (3.13) of \cite{Wu and Yan})
\begin{equation}
F^{Mix}_{\Gamma}=-G^{Mix}_{\Gamma}M^{Mix}_{\Gamma}G^{Mix\ *}_{\Gamma},
\end{equation}
where $G^{Mix}_{\Gamma}:H^{1/2}(\Gamma) \times H^{-1/2}(\Gamma) \to L^{2}(\mathbb{S}^{d-1})$ is the data-to-pattern operator for the mixed crack $\Gamma$, i.e., defined by $G^{Mix}_{\Gamma}\left(\begin{array}{cc}
      f_{1} \\
      f_{2} 
    \end{array}\right):=v^{\infty}$ where $v^{\infty}$ is the far-field pattern of a radiating solution $v$ such that 
\begin{equation}
\Delta v+k^2v=0 \ \mathrm{in} \ \mathbb{R}^d \setminus \Gamma, \ \  v_{-}=f_{1} \ \mathrm{on} \ \Gamma, \ \  \frac{\partial v_{+}}{\partial \nu}=f_{2} \ \mathrm{on} \ \Gamma,
\end{equation}
and $M^{Mix}_{\Gamma}:\tilde{H}^{-1/2}(\Gamma) \times \tilde{H}^{1/2}(\Gamma) \to H^{1/2}(\Gamma) \times H^{-1/2}(\Gamma)$ has the form (see (3.12) of \cite{Wu and Yan})
\begin{equation}
M^{Mix}_{\Gamma}=\left(\begin{array}{cc}
      C^{+}_{\Gamma} & -I \\
      -I & C^{-}_{\Gamma} 
    \end{array}\right) + K^{Mix}_{\Gamma},
\end{equation}
where $C^{+}_{\Gamma}:\tilde{H}^{-1/2}(\Gamma) \to H^{1/2}(\Gamma)$ is some positive coercive operator, $C^{-}_{\Gamma}: \tilde{H}^{1/2}(\Gamma) \to H^{-1/2}(\Gamma)$ is some negative coercive operator, and $K^{Mix}_{\Gamma}:\tilde{H}^{-1/2}(\Gamma) \times \tilde{H}^{1/2}(\Gamma) \to H^{1/2}(\Gamma) \times H^{-1/2}(\Gamma)$ is some compact operator. 
\par
Assume that $\Gamma \subset B$. Define $R:H^{1/2}(\Gamma) \times H^{-1/2}(\Gamma) \to H^{1/2}(\partial B)$ by $R\left(\begin{array}{cc}
      f_{1} \\
      f_{2} 
    \end{array}\right):= v\bigl|_{\partial B}$, then, $R$ is a compact operator and 
\begin{equation}    
G^{Mix}_{\Gamma} = G^{Dir}_{B}R.
\end{equation}
Applying {\bf(1)} of Theorem \ref{General theorem for MM} as 
\begin{equation*}
F=H^{*}_{\partial B}H_{\partial B}=G^{Dir}_{B}\left[C_{B}J_{\partial B}C_{B} + \hat{K}_{B} \right] G^{Dir\ *}_{B},
\end{equation*}
\begin{equation*}
\tilde{F}=-F^{Mix}_{\Gamma}=G^{Mix}_{\Gamma}M^{Mix}_{\Gamma}G^{Mix\ *}_{\Gamma},
\end{equation*}
we have
\begin{equation}
-\mathrm{Re}F^{Mix}_{\Gamma} \leq_{\mathrm{fin}} H^{*}_{\partial B}H_{\partial B}.
\end{equation}
\par
Assume that $\Gamma \not \subset B$. Then, there exists $\tilde{\Gamma} \Subset \Gamma$ such that $\tilde{\Gamma} \cap B = \emptyset$. Setting $W:=\mathrm{Ran}\left(G^{Mix}_{\Gamma}R_{1}^{*}E_{\tilde{\Gamma}}\big|_{X_{\tilde{\Gamma}}}\right) \subset \mathrm{Ran}(G^{Mix}_{\Gamma}R_{1}^{*})$, it is infinite dimensional because $X_{\tilde{\Gamma}}$ is infinite dimensional and the operator $G^{Mix}_{\Gamma}R_{1}^{*}E_{\tilde{\Gamma}}$ is injective. By the same argument in Lemma 4.6 of \cite{Albicker and Griesmaier}, we have 
\begin{equation}
W \cap \mathrm{Ran}(G^{Dir}_{B},G^{Mix}_{\Gamma}R_{2}^{*})=\{0 \}.
\end{equation}
Applying {\bf(3)} of Theorem \ref{General theorem for MM} as
\begin{equation*}
F^{Mix}=-F^{Mix}_{\Gamma}=G^{Mix}_{\Gamma}\left( \left(\begin{array}{cc}
      C^{+}_{\Gamma} & -I \\
      -I & C^{-}_{\Gamma} 
    \end{array}\right) + K^{Mix}_{\Gamma} \right)G^{Mix\ *}_{\Gamma},
\end{equation*}
\begin{equation*}
\tilde{F}=H^{*}_{\partial B}H_{\partial B}=G^{Dir}_{B}\left[C_{B}J_{\partial B}C_{B} + \hat{K}_{B} \right] G^{Dir\ *}_{B},
\end{equation*}
we have
\begin{equation}
-\mathrm{Re}F^{Mix}_{\Gamma} \not \leq_{\mathrm{fin}} H^{*}_{\partial B}H_{\partial B}.
\end{equation}
From the above discussion, we conclude the following theorem.
\begin{thm}\label{Mixed crack part1}
Let $B \subset \mathbb{R}^d$ be a bounded domain. Then,
\begin{equation}
\Gamma \subset B \ \ \ \  \Longleftrightarrow \ \ \ \  -\mathrm{Re}F^{Mix}_{\Gamma} \leq_{\mathrm{fin}} H^{*}_{\partial B}H_{\partial B}.
\end{equation}
\end{thm}
By the same argument in Theorem \ref{Mixed crack part1}, one can apply {\bf(1)} and {\bf(3)} replacing {\bf(3b)} by {\bf(3b)'} of Theorem \ref{General theorem for MM} as
\begin{equation*}
F^{Mix}=F^{Mix}_{\Gamma}=-G^{Mix}_{\Gamma}\left( \left(\begin{array}{cc}
      C^{+}_{\Gamma} & -I \\
      -I & C^{-}_{\Gamma} 
    \end{array}\right) + K^{Mix}_{\Gamma} \right)G^{Mix\ *}_{\Gamma},
\end{equation*}
\begin{equation*}
\tilde{F}=H^{*}_{\partial B}H_{\partial B}=G^{Dir}_{B}\left[C_{B}J_{\partial B}C_{B} + \hat{K}_{B} \right] G^{Dir\ *}_{B},
\end{equation*}
Then, we also conclude the following theorem.
\begin{thm}\label{Mixed crack part2}
Let $B \subset \mathbb{R}^d$ be a bounded domain. Then,
\begin{equation}
\Gamma \subset B \ \ \ \  \Longleftrightarrow \ \ \ \  \mathrm{Re}F^{Mix}_{\Gamma} \leq_{\mathrm{fin}} H^{*}_{\partial B}H_{\partial B}.
\end{equation}
\end{thm}
\begin{rem}
The factorization reconstruction for the mixed crack has been studied in \cite{Wu and Yan}, but an extensive closed curve of the unknown crack should be known, which is a very restrictive assumption (see Theorem 3.3 of \cite{Wu and Yan}). However, the monotonicity reconstruction (Theorems \ref{Mixed crack part1} and \ref{Mixed crack part2}) does not require it. We also remark that our works in this section would be a new extension of the monotonicity to the inverse acoustic mixed crack scattering.
\end{rem}

\section{Numerical examples}
In this section, we study numerical examples of our theoretical results in 2 dimensions. The far-field operator $F$ is approximated by the matrix
\begin{equation}
F \approx \frac{2\pi}{N} \bigl(u^{\infty}(\hat{x}_l, \theta_m) \bigr)_{1 \leq l,m \leq N} \in \mathbb{C}^{N \times N},
\end{equation}
where $\hat{x}_l=\bigl(\mathrm{cos}(\frac{2\pi l}{N}), \mathrm{sin}(\frac{2\pi l}{N}) \bigr)$ and $\theta_m=\bigl(\mathrm{cos}(\frac{2\pi m}{N}), \mathrm{sin}(\frac{2\pi m}{N}) \bigr)$. For the far-field pattern $u^{\infty}(\hat{x}, \theta)$ of the Dirichlet obstacle and the Dirichlet crack, we numerically compute the integral
\begin{equation}
u^{\infty}(\hat{x},\theta)=\frac{e^{i\frac{\pi}{4}}}{\sqrt{8\pi k}}\int_{ \Gamma}e^{-ik\hat{x}y}\varphi_{\theta}(y)ds(y), \end{equation}
where $\varphi_{\theta}$ is given by solving
\begin{equation}
-e^{ikx\cdot \theta}=\int_{\Gamma}\Phi(x,y) \varphi_{\theta}(y)ds(y), \ x \in \Gamma.
\end{equation}
For the inhomogeneous medium, we numerically compute the integral
\begin{equation}
u^{\infty}(\hat{x},\theta)=\frac{e^{i\frac{\pi}{4}}}{\sqrt{8\pi k}}\int_{ \Omega}e^{-ik\hat{x}y}q(y)(u_{\theta}(y)+e^{iky\cdot \theta})dy, 
\end{equation}
where $u_{\theta}$ is given by solving 
\begin{equation}
u(x)=k^{2}\int_{\Omega}\Phi(x,y)q(y)( u_{\theta}(y)+e^{iky\cdot \theta})dy, \ x\in \Omega.
\end{equation}
\par
\begin{figure}[h]
  \centering
  \includegraphics[scale=0.35]{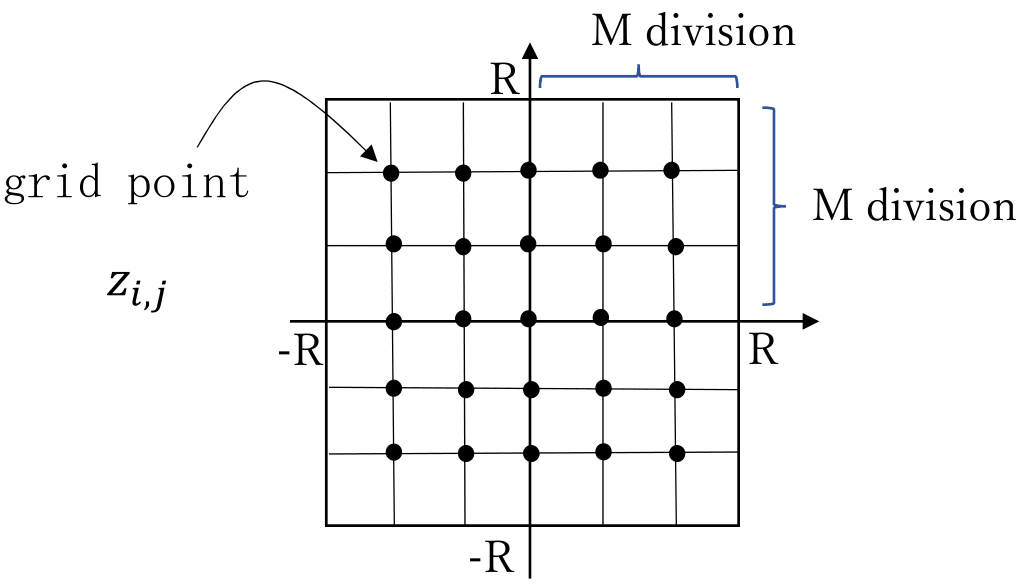}
  \caption{The grid points in the sampling square region.}\label{grid}
\end{figure}

For a bounded domain (or a smooth curve) $B$, the operator $H^{*}_{B}H_{B}$ is approximated by
\begin{equation}
H^{*}_{B}H_{B} \approx \frac{2\pi}{N} \biggl(\int_{B}e^{iky\cdot(\theta_m-\hat{x}_l)}dy \biggr)_{1 \leq l,m \leq N} \in \mathbb{C}^{N \times N}.
\end{equation}
We denote by the sampling square region $[-R, R]^{2}$, which includes the unknown target. We also denote by a grid point $z_{i,j}:=(\frac{Ri}{M}, \frac{Rj}{M})$ ($i,j = -M, -M+1, ..., M$) in the sampling square region $[-R, R]^{2}$ (see Figure \ref{grid}). Throughout our examples, we fix parameters $R=1.5$, $M=100$, and $N=20$.
\subsection{Dirichlet obstacle and inhomogeneous medium}
\begin{figure}[h]
  \centering
  \hspace{-2cm}
  \includegraphics[scale=0.7]{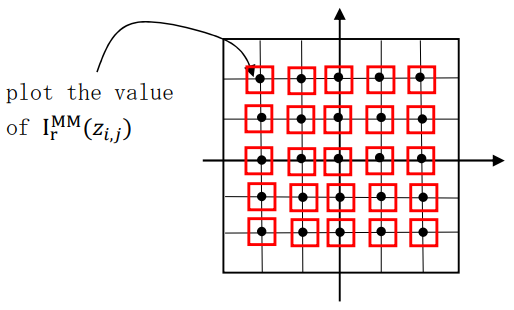}
  \caption{The square test.}\label{box}
\end{figure}

\begin{figure}[htbp]
\begin{tabular}{c}
\begin{minipage}{0.5\hsize}
  \begin{center}
   \includegraphics[scale=0.5]{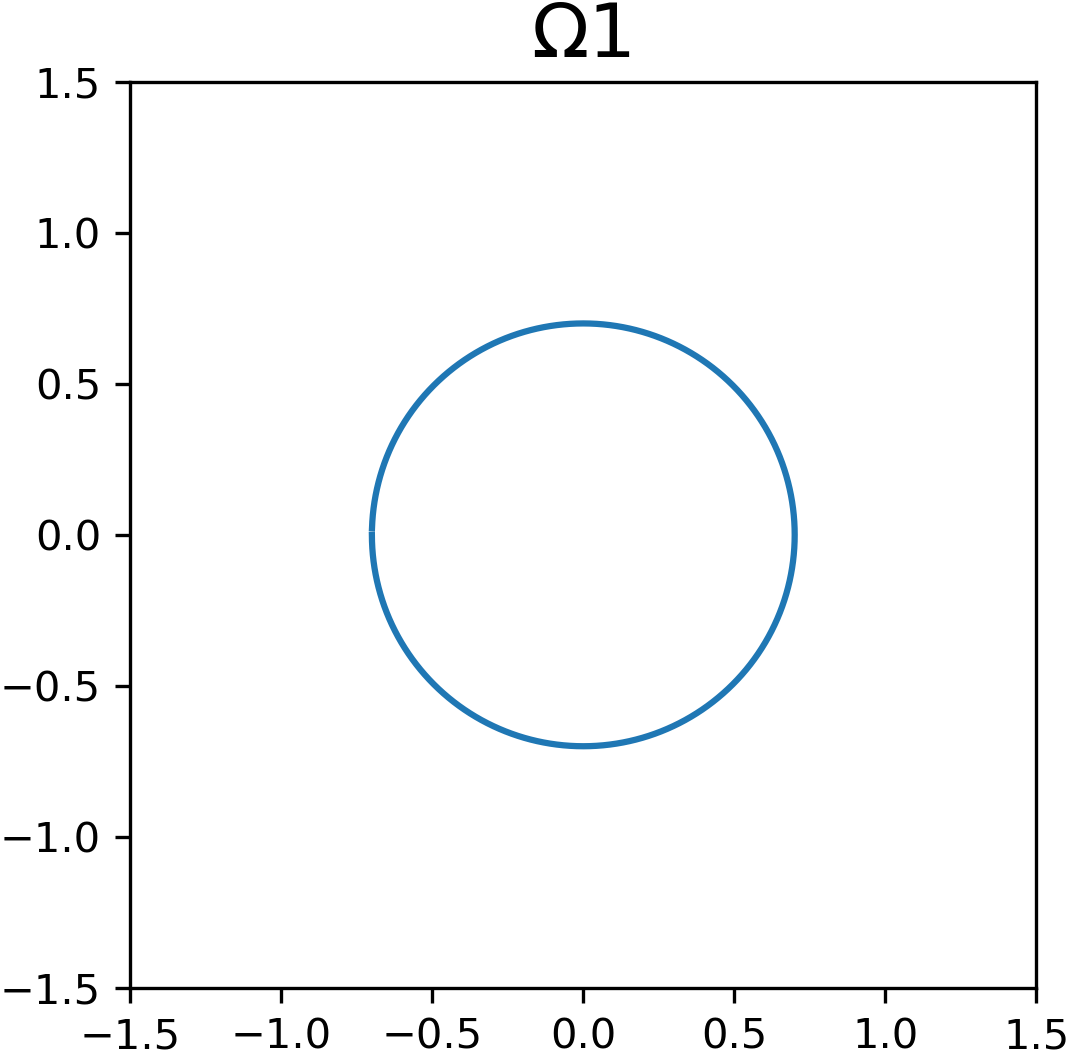}
  \end{center}
 \end{minipage}
 \begin{minipage}{0.5\hsize}
 \begin{center}
  \includegraphics[scale=0.5]{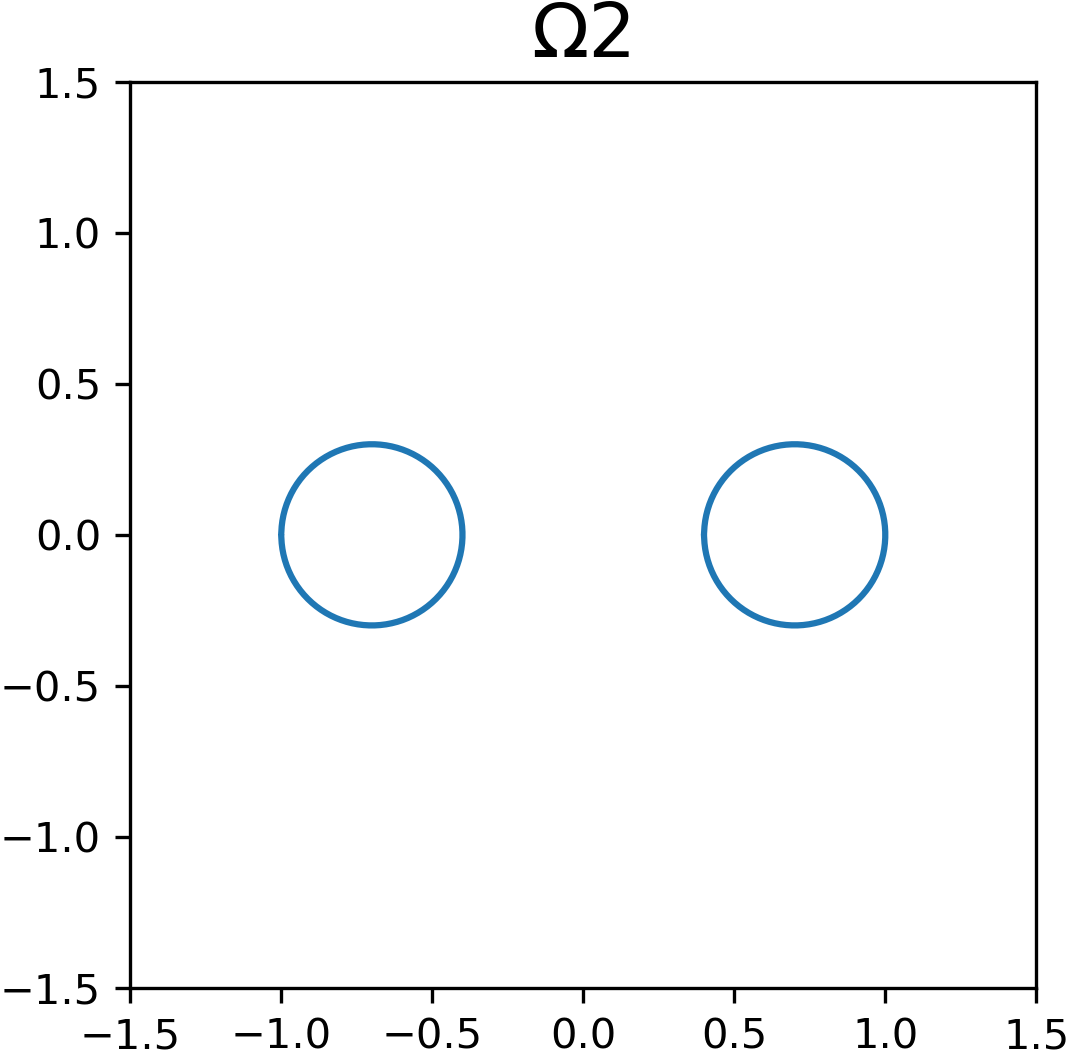}
 \end{center}
 \end{minipage}
\end{tabular}
\caption{The original domains $\Omega_1$ (left) and $\Omega_2$ (right).}\label{The original domains-1}
\vspace{3cm}
\end{figure}

We consider the Dirichlet obstacle and the inhomogeneous medium detections discussed in Sections 5.1 and 5.2. The following shapes $\partial \Omega_j$ ($j=1,2$) are considered (see Figure \ref{The original domains-1}):
\begin{equation}
\partial \Omega_1=\left\{\left(0.7\mathrm{cos}(\pi s), 0.7\mathrm{sin}(\pi s)\right) | -1\leq s \leq 1 \right\}.
\end{equation}
\begin{equation*}
\partial \Omega_2=\left\{\left(0.3\mathrm{cos}(\pi s)-0.7, \ 0.3\mathrm{sin}(\pi s) \right) | -1\leq s \leq 1 \right\}
\end{equation*}
\begin{equation}
\cup\left\{\left(0.3\mathrm{cos}(\pi s)+0.7, \ 0.3\mathrm{sin}(\pi s) \right) | -1\leq s \leq 1 \right\}.
\end{equation}
Based on Theorems \ref{Dirichlet obstacle inside} and \ref{Medium inside}, the indicator functions for the Dirichlet obstacle and inhomogeneous medium are given by
\begin{equation}
I^{MM}_{dir}(B):= \# \left\{\mathrm{negative} \ \mathrm{eigenvalues}  \ \mathrm{of} \  -\mathrm{Re}F^{Dir}_{\Omega}-H^{*}_{\partial B}H_{\partial B} \right\},
\end{equation}
\begin{equation}
I^{MM}_{med}(B):= \# \left\{\mathrm{negative} \ \mathrm{eigenvalues}  \ \mathrm{of} \ \mathrm{Re}F^{Med}_{\Omega}-\alpha H^{*}_{\partial B}H_{\partial B} \right\},
\end{equation}
for a bounded domain $B$, respectively. For the medium, we always consider $q=1$ in $\Omega_j$, and $\alpha=1$. Here, $B$ is chosen as a square, i.e., $B=B_{i,j}(r):=z_{i,j}+[-\frac{r}{2},\frac{r}{2}]\times[-\frac{r}{2},\frac{r}{2}]$ where a grid point $z_{i,j}$ is the center of a square (see Figure \ref{box}). 

Then, we can compute the integral
\begin{equation}
\int_{B_{i,j}(r)}e^{iky\cdot(\theta_m-\hat{x}_l)}dy=r^{2} e^{ik(\theta_m-\hat{x}_l)\cdot z_{i,j}}\mathrm{sinc}\left(\frac{kr}{2}(\theta_m-\hat{x}_l)_{1} \right)\mathrm{sinc}\left(\frac{kr}{2}(\theta_m-\hat{x}_l)_{2} \right).
\end{equation}
Figures \ref{MM obstacle} and \ref{MM medium} are given by plotting the values of the indicator functions 
\begin{equation}
I^{MM}_{dir, r}(z_{i,j}):=I^{MM}_{dir}(B_{i,j}(r)), \ \ \mathrm{for \ each} \ i, j = -100, -99, ..., 100,
\end{equation}
\begin{equation}
I^{MM}_{med, r}(z_{i,j}):=I^{MM}_{med}(B_{i,j}(r)), \ \ \mathrm{for \ each} \ i, j = -100, -99, ..., 100,
\end{equation}
respectively, for different lengths $r=0.1, 0.5$, wavenumbers $k=1,5$, and shapes $\Omega_1, \Omega_2$.
\par
We also plot the values of the indicator function for the factorization method
\begin{equation}
I^{FM}(z_{i,j}):=\left( \sum_{n=1}^{\infty}\frac{|\langle \phi_{z_{i,j}}, \varphi_{n} \rangle_{L^{2}(\mathbb{S}^{1})}|^{2}}{\mu_{n}} \right)^{-1}, \ \ \mathrm{for \ each} \ i, j = -100, -99, ..., 100,
\end{equation}
with $\left\{ \mu_n, \phi_n \right\}$ an eigensystem of the self-adjoint compact operator $|\mathrm{Re}F|+|\mathrm{Im}F|$ where $F$ is some far-field operator, and $\phi_{z_{i,j}}$ is defined by
\begin{equation}
\phi_{z_{i,j}}(\hat{x}):=e^{-ik\hat{x}\cdot z_{i,j}}, \hat{x} \in \mathbb{S}^{1}, \label{FFP of Green}
\end{equation}
Figures \ref{FM obstacle} and \ref{FM medium} are corresponding to $F=F^{Dir}_{\Omega}$ for the Dirichlet obstacle and $F=F^{Med}_{\Omega}$ for the inhomogeneous medium, respectively. For details of introductions of the indicator function $I^{FM}$, we refer to Corollary 2.16 and Theorem 4.9 of \cite{Kirsch and Grinberg}.
\subsection{Dirichlet crack}
\begin{figure}[h]
  \centering
  \hspace{-2cm}
  \includegraphics[scale=0.7]{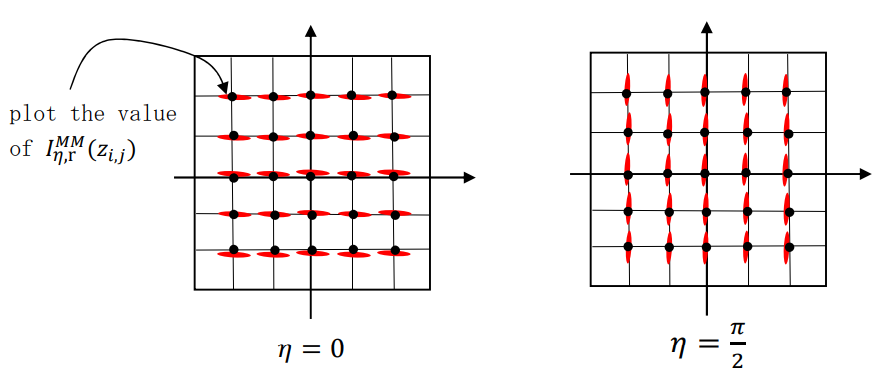}
  \caption{The line segment test.}\label{line segment test}
\end{figure}

\begin{figure}[htbp]
\begin{tabular}{c}
\begin{minipage}{0.5\hsize}
  \begin{center}
   \includegraphics[scale=0.5]{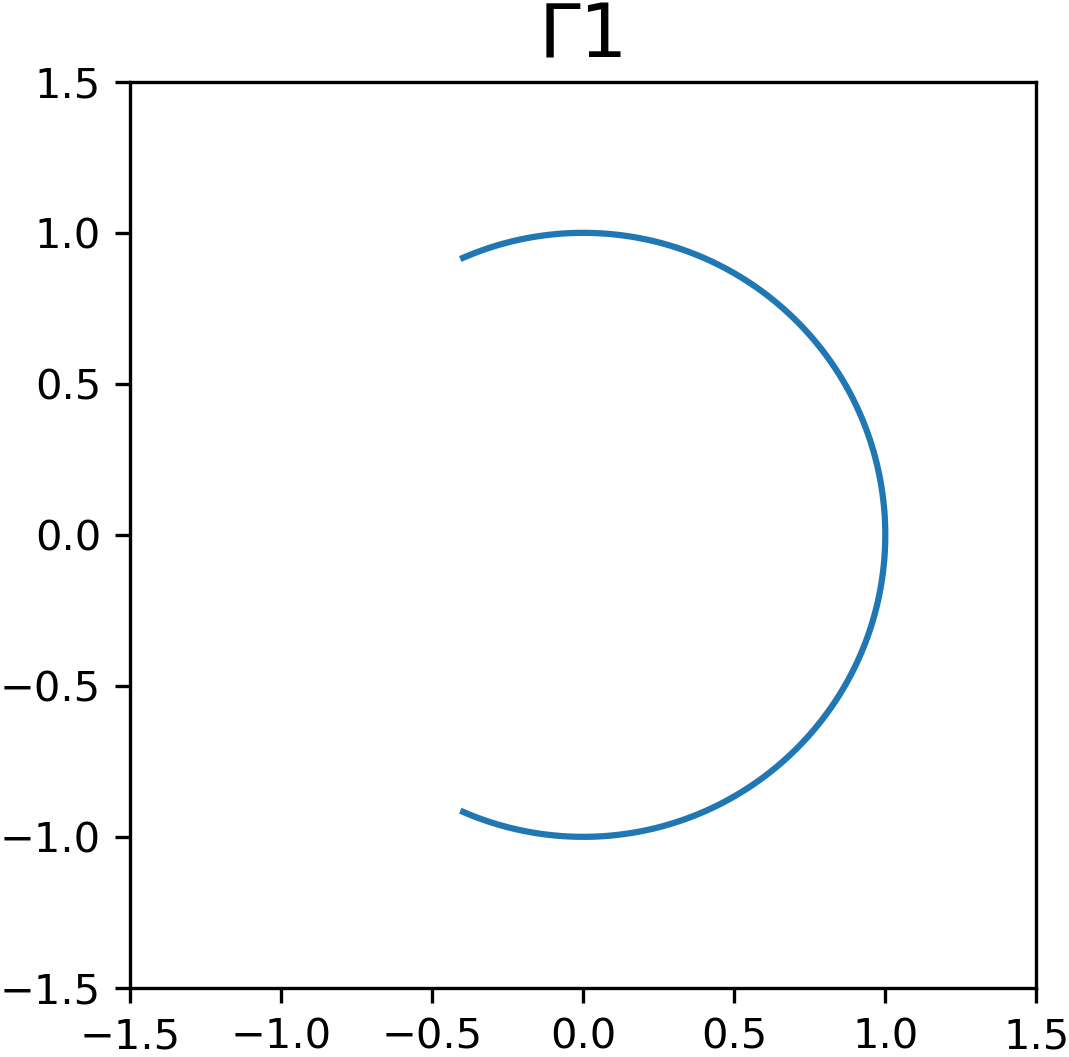}
  \end{center}
 \end{minipage}
 \begin{minipage}{0.5\hsize}
 \begin{center}
  \includegraphics[scale=0.5]{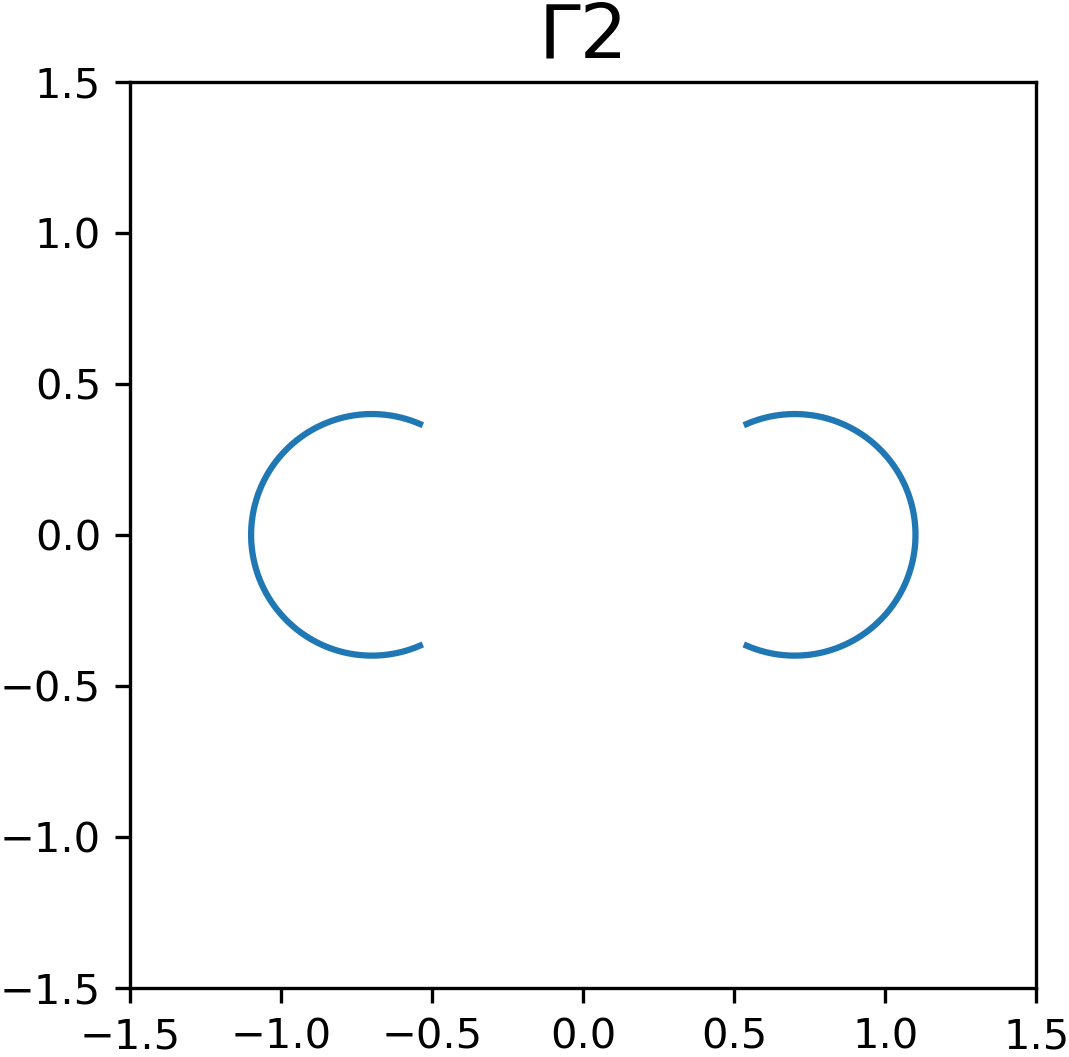}
 \end{center}
 \end{minipage}
\end{tabular}
\caption{The original open arcs $\Gamma_1$ (left) and $\Gamma_2$ (right).}\label{The original open arcs-1}
\vspace{3cm}
\end{figure}

We consider the Dirichlet crack detection discussed in Section 5.3. The following shapes $\Gamma_j$ ($j=1,2$) are considered (see Figure \ref{The original open arcs-1}):
\begin{equation}
    \Gamma_1=\left\{\left(\mathrm{cos}(2s), \mathrm{sin}(2s)\right) | -1\leq s \leq 1 \right\}.
\end{equation}
\begin{equation*}
\Gamma_2=\left\{\left(-0.4\mathrm{cos}(2s)-0.7, \ 0.4\mathrm{sin}(2s) \right) | -1\leq s \leq 1 \right\}
\end{equation*}
\begin{equation}
\cup\left\{\left(0.4\mathrm{cos}(2s)+0.7, \ 0.4\mathrm{sin}(2s) \right) | -1\leq s \leq 1 \right\}.
\end{equation}
Based on Theorem \ref{Dirichlet crack inside}, the indicator function is given by
\begin{equation}
I^{MM}_{dir}(\sigma):= \# \left\{\mathrm{negative} \ \mathrm{eigenvalues}  \ \mathrm{of} -\mathrm{Re}F^{Dir}_{\Gamma}-H^{*}_{\sigma}H_{\sigma} \right\}.
\end{equation}
for a smooth arc $\sigma$. Here, $\sigma$ is chosen as a line segment, i.e., $\sigma=\sigma_{i,j}(\eta,r):=z_{i,j}+L(\eta,r)$ where a grid point $z_{i,j}$ is the center of line segments, and $L(\eta,r)$ is defined by
\begin{equation}
L(\eta,r):=\left\{ (s, s\mathrm{tan}(\eta))\bigl| -\frac{r}{2}\mathrm{cos}(\eta) \leq s \leq \frac{r}{2}\mathrm{cos}(\eta) \right\},
\end{equation}
where $\eta \in [0, \pi]$ is the angle and $r>0$ is the length of the line segment (see Figure \ref{line segment test}).

Then, we can compute the integral
\begin{equation}
\int_{\sigma_{i,j}(\eta,r)}e^{iky\cdot(\theta_m-\hat{x}_l)}ds(y)=re^{ik(\theta_m-\hat{x}_l)\cdot z_{i,j}}\mathrm{sinc}\biggl(\frac{rk}{2} \Bigl( \mathrm{cos}(\eta)(\theta_m-\hat{x}_l)_{1}+\mathrm{sin}(\eta)(\theta_m-\hat{x}_l)_{2} \Bigr) \biggr).
\end{equation}
Figures \ref{MM crack-1} and \ref{MM crack-2} are given by plotting the values of the indicator function for $\Gamma_1$ and $\Gamma_2$
\begin{equation}
I^{MM}_{dir,\eta,r}(z_{i,j}):=I^{MM}_{dir}(\sigma_{i,j}(\eta,r)), \ \ \mathrm{for \ each} \ i, j = -100, -99, ..., 100,
\end{equation}
respectively, for different angles $\eta=0, \pi/2$, lengths $r=0.01, 0.1$, and wavenumbers $k=1,5$. 
\par
In Figure \ref{FM crack} we also plot the values of the indicator function for the factorization method
\begin{equation}
I^{FM}_{dir}(z_{i,j}):=\left( \sum_{n=1}^{\infty}\frac{|\langle \phi_{z_{i,j}}, \varphi_{n} \rangle_{L^{2}(\mathbb{S}^{1})}|^{2}}{\mu_{n}} \right)^{-1}, \ \ \mathrm{for \ each} \ i, j = -100, -99, ..., 100,
\end{equation}
with $\left\{ \mu_n, \phi_n \right\}$ an eigensystem of the self-adjoint compact operator $|\mathrm{Re}F^{Dir}_{\Gamma}|+|\mathrm{Im}F^{Dir}_{\Gamma}|$ and $\phi_{z_{i,j}}$ is defined in (\ref{FFP of Green}). For details of introductions of the indicator function $I^{FM}_{dir}$, we refer to Theorem 3.9 of \cite{Kirsch and Ritter}.
\subsection{Mixed crack}
\begin{figure}[htbp]
\hspace{-10mm}
\begin{tabular}{c}
\begin{minipage}{0.55\hsize}
  \begin{center}
   \includegraphics[scale=0.5]{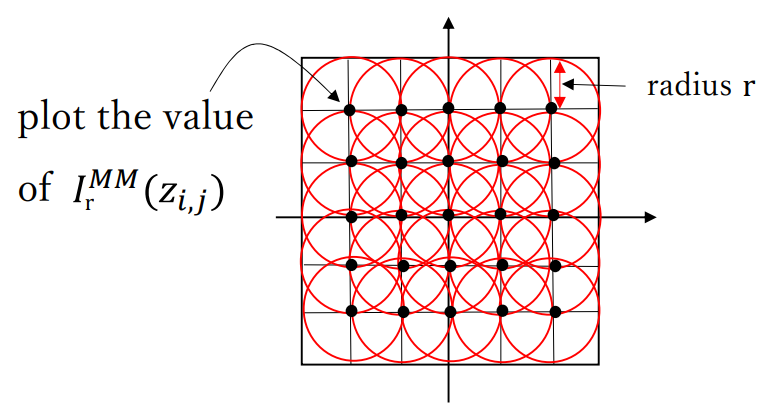}
  \end{center}
 \end{minipage}
 \hspace{-1.5cm}
 \begin{minipage}{0.55\hsize}
 \begin{center}
  \includegraphics[scale=0.5]{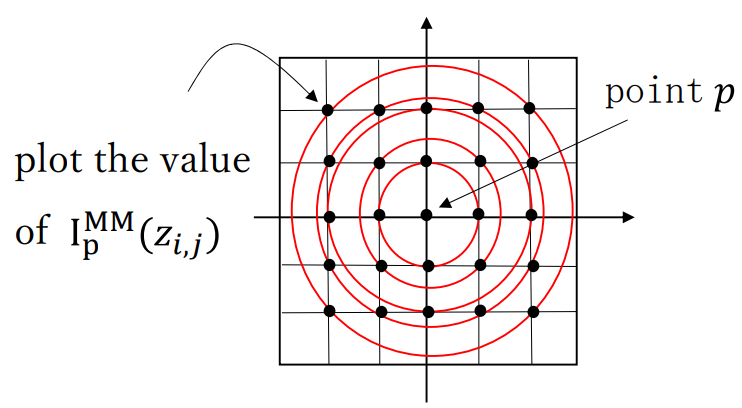}
 \end{center}
 \end{minipage}
\end{tabular}
\caption{The shifting circle test(left) and the shrinking circle test(right).}\label{circle test}
\end{figure}

\begin{figure}[h] 
  \begin{center}
   \includegraphics[scale=0.5]{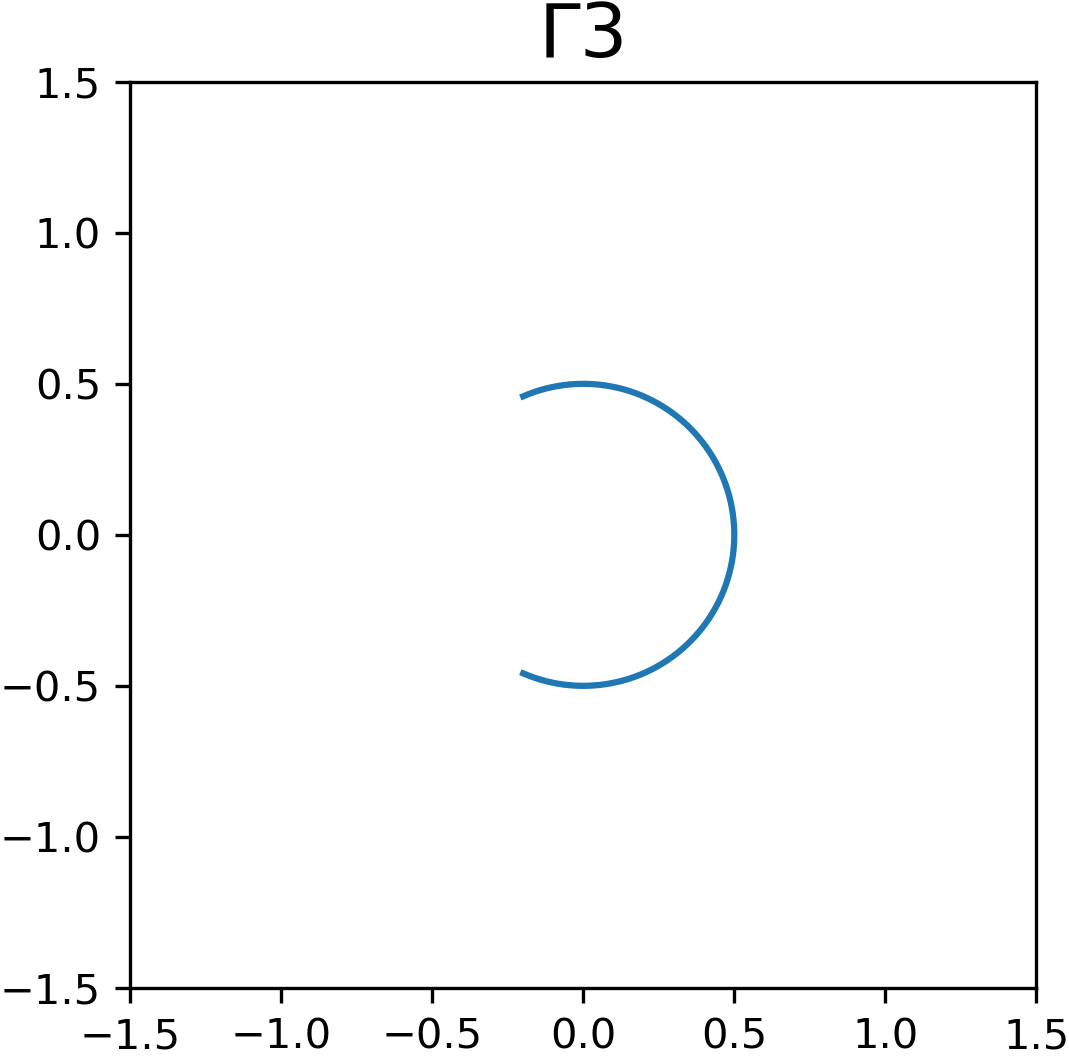}
  \end{center}
   \caption{The original open arc $\Gamma_3$.}\label{The original open arc-2}
\vspace{3cm}
\end{figure}

We consider the mixed crack detection discussed in Section 5.5. The following shape $\Gamma_3$ is considered (see Figure \ref{The original open arc-2}):
\begin{equation}
\Gamma_3=\left\{\left(0.5\mathrm{cos}(2s), 0.5\mathrm{sin}(2s)\right) | -1\leq s \leq 1 \right\}.
\end{equation}
Based on Theorem \ref{Mixed crack part1}, the indicator function is given by
\begin{equation}
I^{MM}_{mix}(B):= \# \left\{\mathrm{negative} \ \mathrm{eigenvalues}  \ \mathrm{of}\ \mathrm{Re}F^{Mix}_{\Gamma}+H^{*}_{\partial B}H_{\partial B} \right\}.
\end{equation}
for a bounded domain $B$. Here, $B$ is chosen as a circle, i.e., $B=B_{r}(z)$ is an open circle with center $z \in \mathbb{R}^{2}$ and radius $r>0$. Then, we can compute the integral
\begin{equation}
\int_{\partial B_{r}(z)}e^{iky\cdot(\theta_m-\hat{x}_l)}ds(y)=2\pi r e^{ik(\theta_m-\hat{x}_l)\cdot z}J_{0}(kr|\theta_m-\hat{x}_l|),
\end{equation}
where $J_0$ is the Bessel function of the first kind for the zero order. Figures \ref{MM mixed crack-shifting} and \ref{MM mixed crack-shrinking} are given by plotting the values of two types (see Figure \ref{circle test}) of indicator functions
\begin{equation}
I^{MM}_{mix,r}(z_{i,j}):= I^{MM}_{mix}(B_{r}(z_{i,j})), \ \ \mathrm{for \ each} \ i, j = -100, -99, ..., 100,
\end{equation}
\begin{equation}
I^{MM}_{mix,p}(z_{i,j}):=I^{MM}_{mix}(B_{|z_{i,j}-p|}(p)), \ \ \mathrm{for \ each} \ i, j = -100, -99, ..., 100,
\end{equation}
respectively, for different radiuses $r=0.25, 1$, points $p=(0,0)$, $(1,0)$, wavenumbers $k=1, 5$.

In Figures \ref{FM mixed crack} we also plot the values of the indicator function for the factorization method
\begin{equation}
I^{FM}_{mix}(z_{i,j}):=\left( \sum_{n=1}^{\infty}\frac{|\langle \phi_{z_{i,j}}, \varphi_{n} \rangle_{L^{2}(\mathbb{S}^{1})}|^{2}}{\mu_{n}} \right)^{-1}, \ \ \mathrm{for \ each} \ i, j = -100, -99, ..., 100,
\end{equation}
with $\left\{ \mu_n, \phi_n \right\}$ an eigensystem of the self-adjoint compact operator $|\mathrm{Re}\tilde{F}|+|\mathrm{Im}\tilde{F}|$. The operator $\tilde{F}$ is defined by adding the Herglotz operator correspoding to the auxiliary $\partial \Omega_{3}$ to the original far-field operator
\begin{equation}
\tilde{F}:=F^{Mix}_{\Gamma} - p H^{*}_{\partial \Omega_3}H_{\partial \Omega_3}, \label{auxiliary}
\end{equation}
for $p \in \mathbb{C}\setminus\{0\}$. We choose $p=0.01+0.01i$ or $0.1+0.1i$ in our numerical examples. Here, the auxiliary $\partial \Omega_{3}$ is defined by
\begin{equation}
\partial \Omega_{3}=\left\{\left(0.5\mathrm{cos}(\pi s), 0.5\mathrm{sin}(\pi s)\right) | -1\leq s \leq 1 \right\},
\end{equation}
which is an extension of $\Gamma_3$, that is, $\Gamma_3 \subset \partial \Omega_{3}$. For details of introductions of the indicator function $I^{FM}_{mix}$, we refer to Theorem 3.4 of \cite{Wu and Yan}.
\section*{Conclusions}
In this paper, we studied the factorization and monotonicity method for inverse acoustic scattering problems. The main contribution was to give a new general functional analysis theorem (Theorem \ref{General theorem for MM}) for the monotonicity method, which can provide reconstruction schemes under weaker {\it a priori} assumptions rather than the factorization method (see the assumptions in Theorems \ref{General theorem for FM} and \ref{General theorem for MM}). Furthermore, we observed that the factorization method needs the real and imaginary parts of the far-field operator (see (\ref{range identity})), while the monotonicity needs only the real part (see (\ref{conclusion(1)}), (\ref{conclusion(2)}), and (\ref{conclusion(3)})), which is also the advantage over the factorization in terms of less data. After proving the general theorem, we also showed how the general theorem is applied to three typical inverse scattering problems (obstacle in Sections 5.1 and 5.4, medium in Section 5.2, and crack in Sections 5.3 and 5.5). However, it can be applied to other inverse problems as well (especially for inverse problems that the factorization method already studied, e.g., inverse scatterings by a layered medium \cite{Bondarenko and Kirsch and Liu}, a mixed-type scatterer of a obstacle and a medium \cite{Kirsch and Liu}, and an obstacle in an homogeneous half-space \cite{Grinberg}).
\par
\begin{figure}[htbp]
\vspace{-3cm}
\begin{tabular}{c}
\begin{minipage}{0.5\hsize}
  \begin{center}
   \includegraphics[scale=0.5]{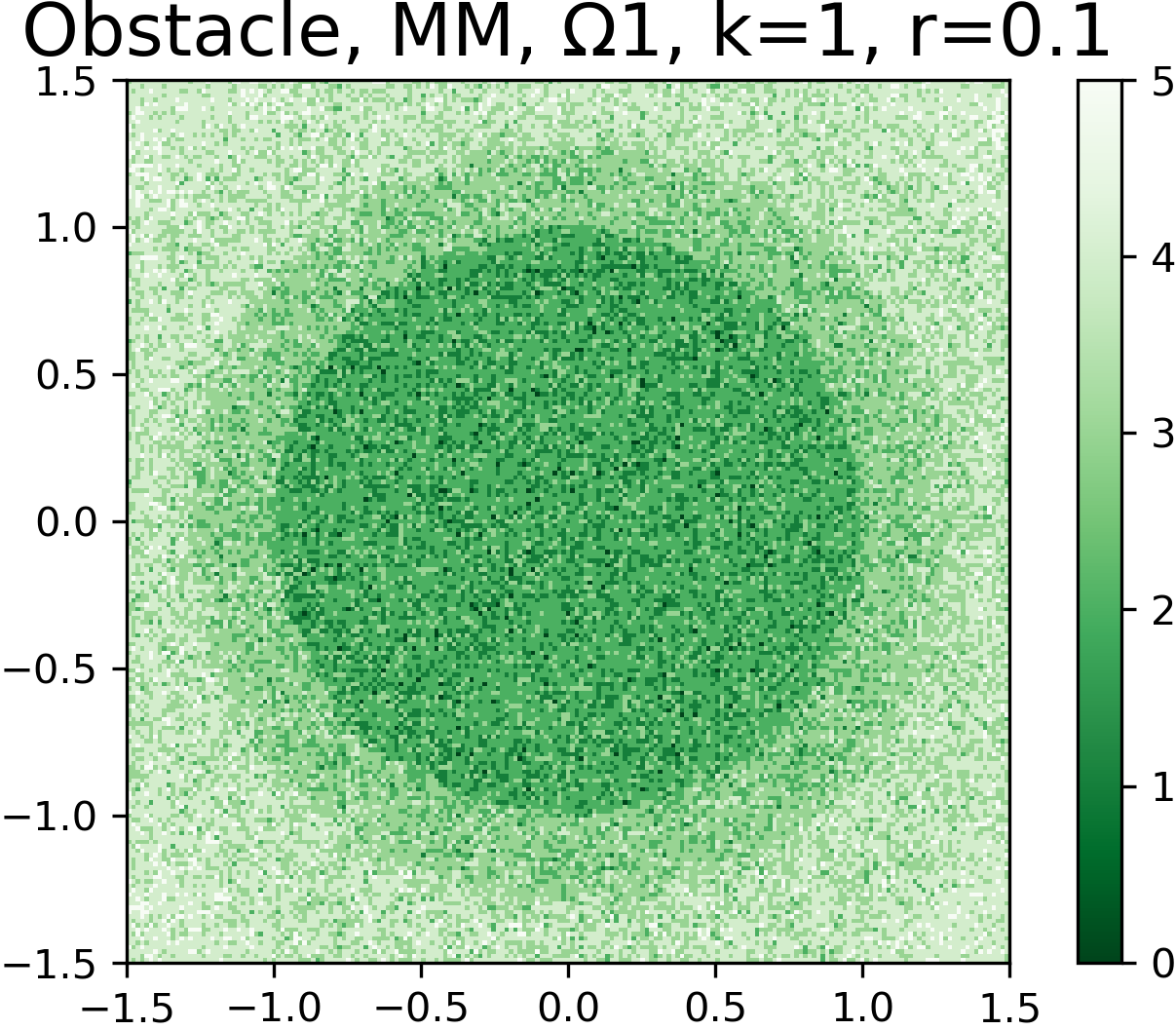}
  \end{center}
 \end{minipage}
 \begin{minipage}{0.5\hsize}
 \begin{center}
  \includegraphics[scale=0.5]{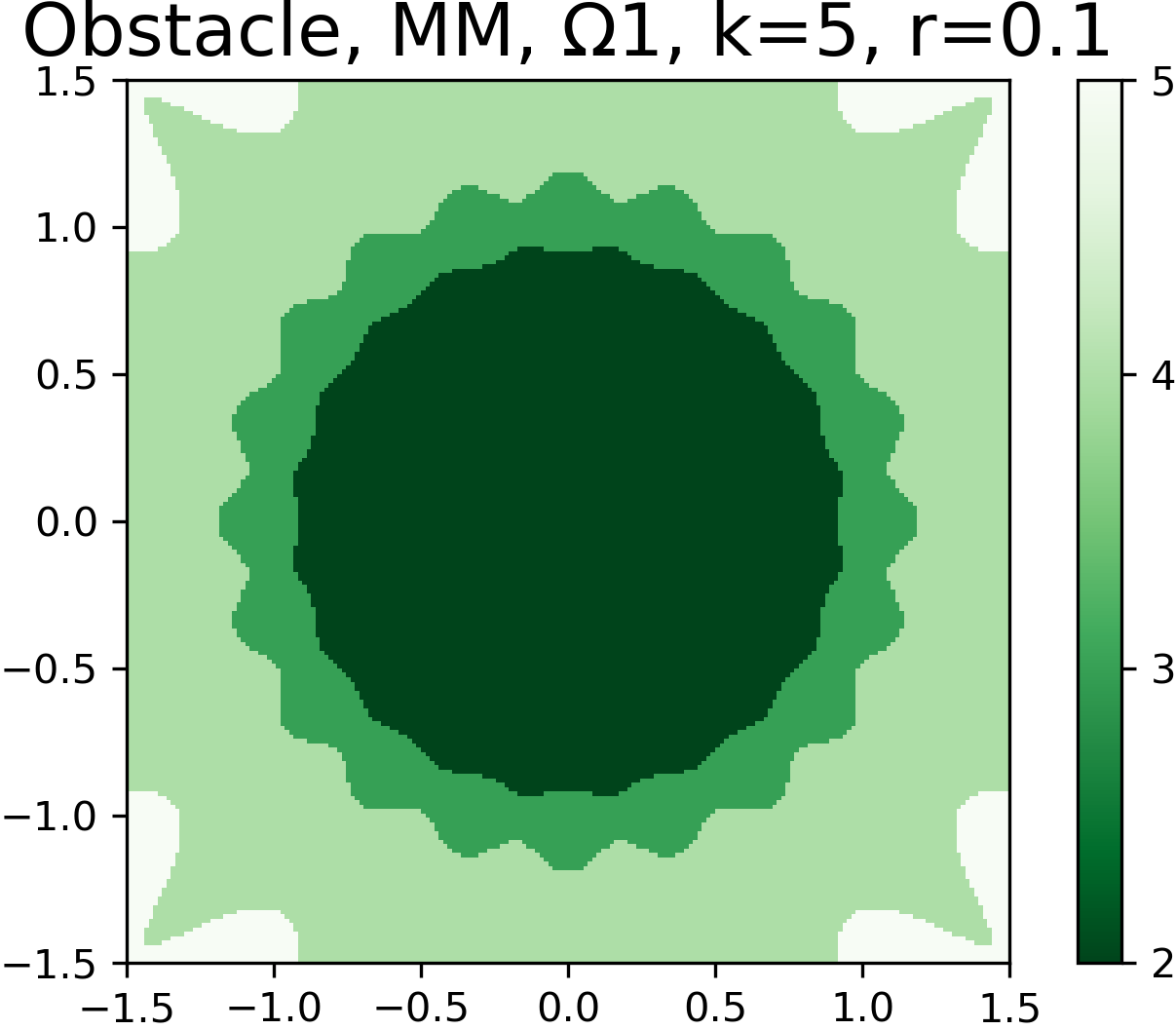}
 \end{center}
 \end{minipage}
 \vspace{1cm} \\ 
\begin{minipage}{0.5\hsize}
  \begin{center}
   \includegraphics[scale=0.5]{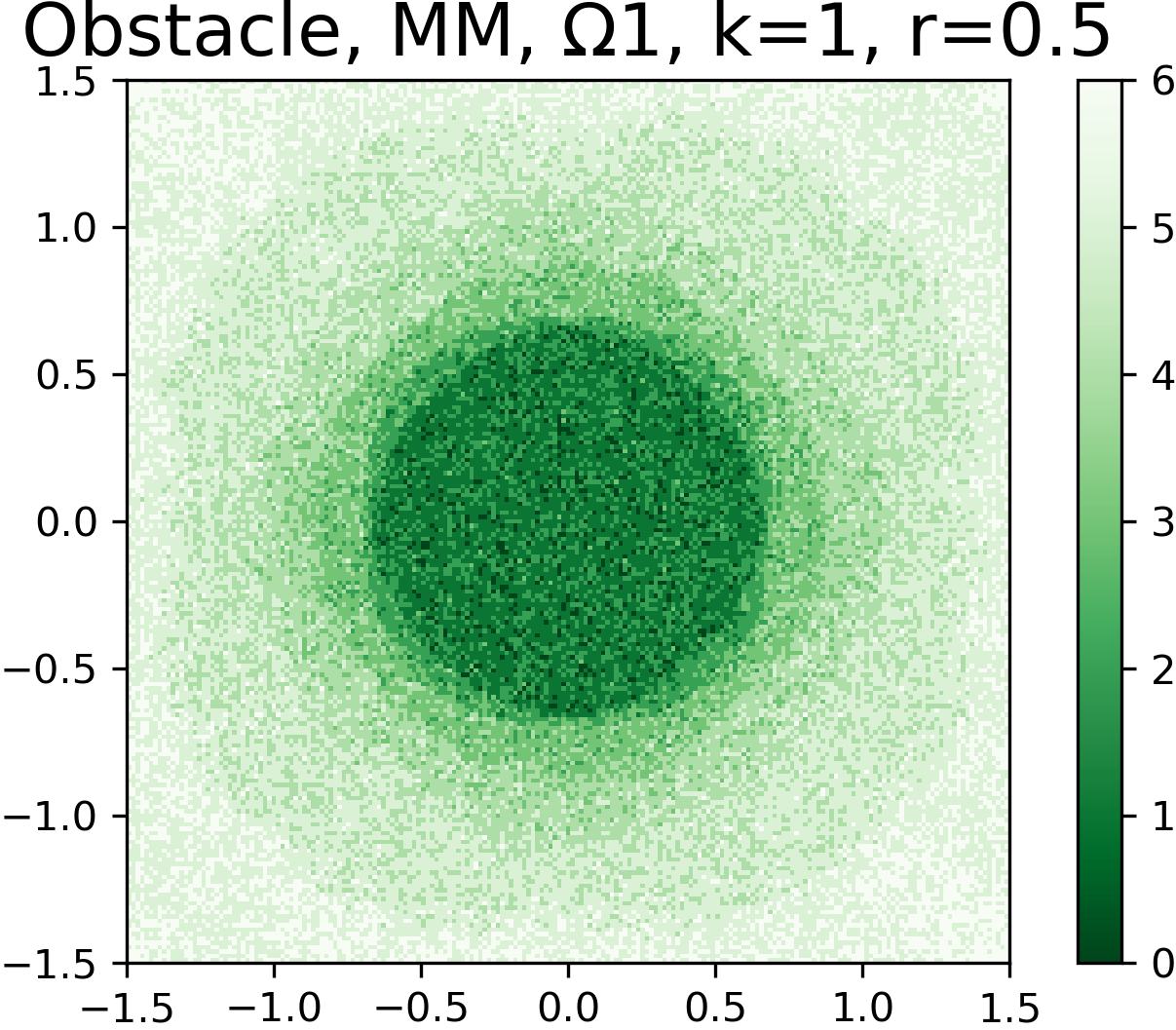}
  \end{center}
 \end{minipage}
 \begin{minipage}{0.5\hsize}
 \begin{center}
  \includegraphics[scale=0.5]{figs/MMobst1k1r05}
 \end{center}
 \end{minipage}
 \vspace{1cm} \\ 
\begin{minipage}{0.5\hsize}
  \begin{center}
   \includegraphics[scale=0.5]{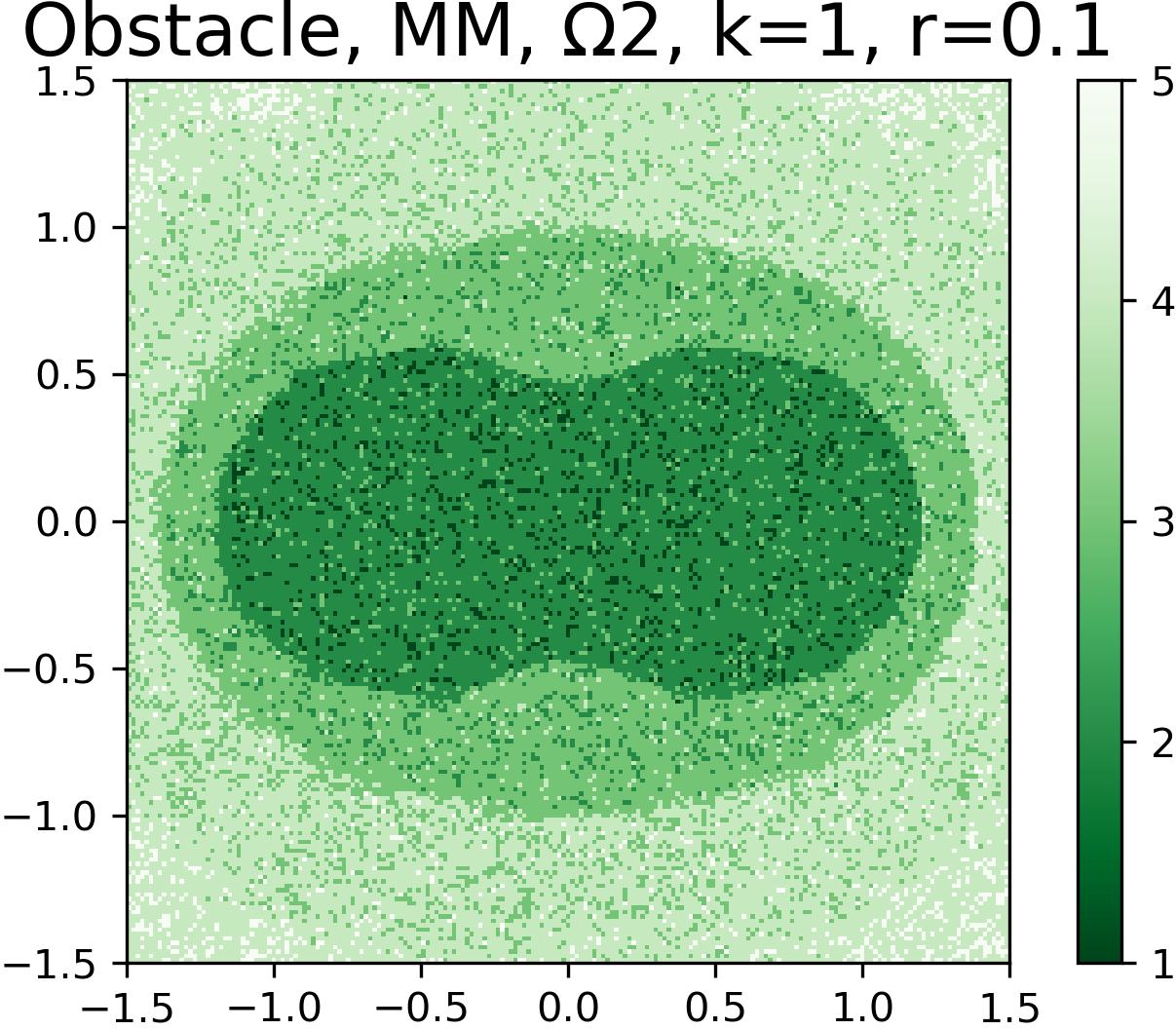}
  \end{center}
 \end{minipage}
 \begin{minipage}{0.5\hsize}
 \begin{center}
  \includegraphics[scale=0.5]{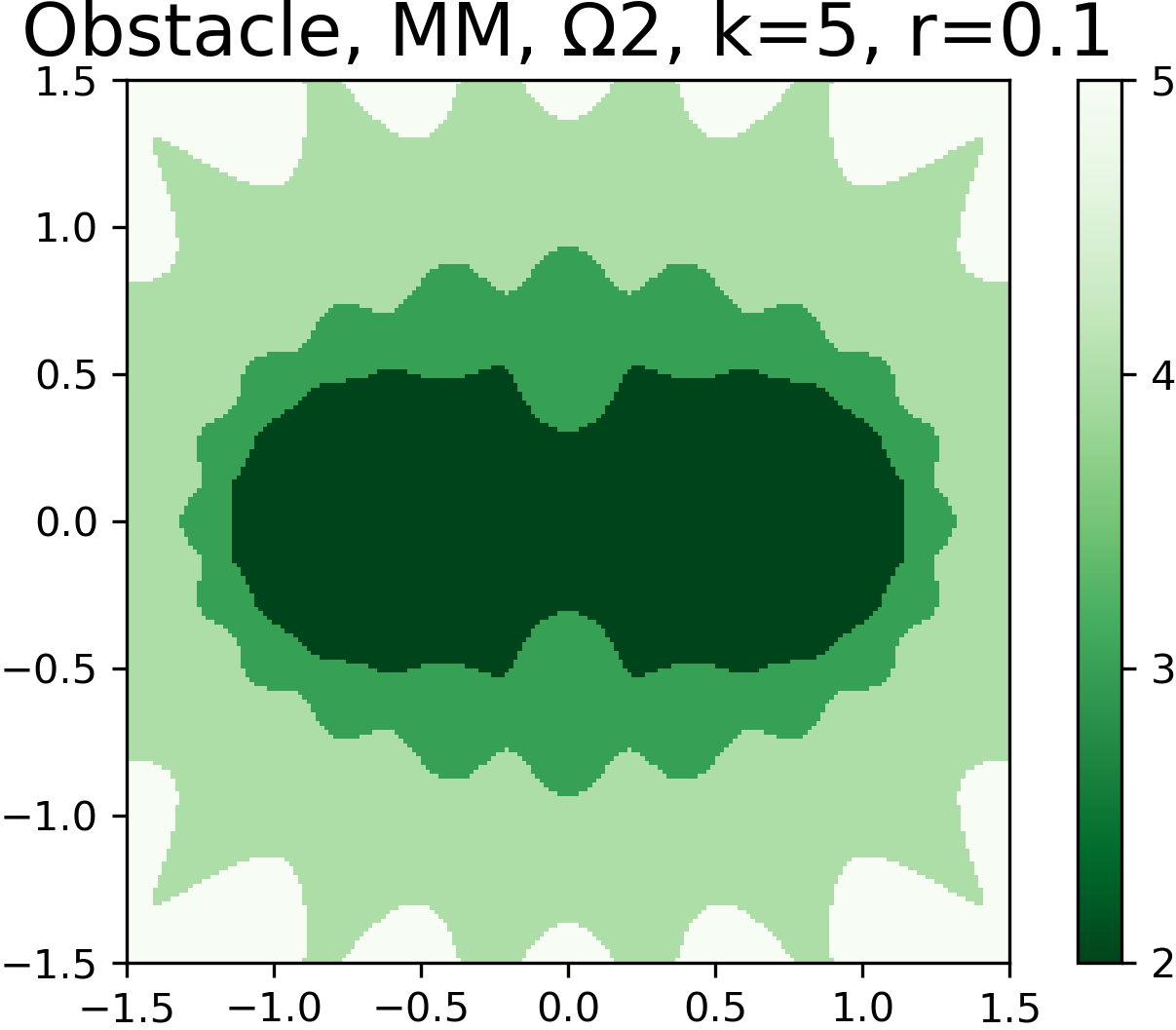}
 \end{center}
 \end{minipage}
 \vspace{1cm} \\ 
\begin{minipage}{0.5\hsize}
  \begin{center}
   \includegraphics[scale=0.5]{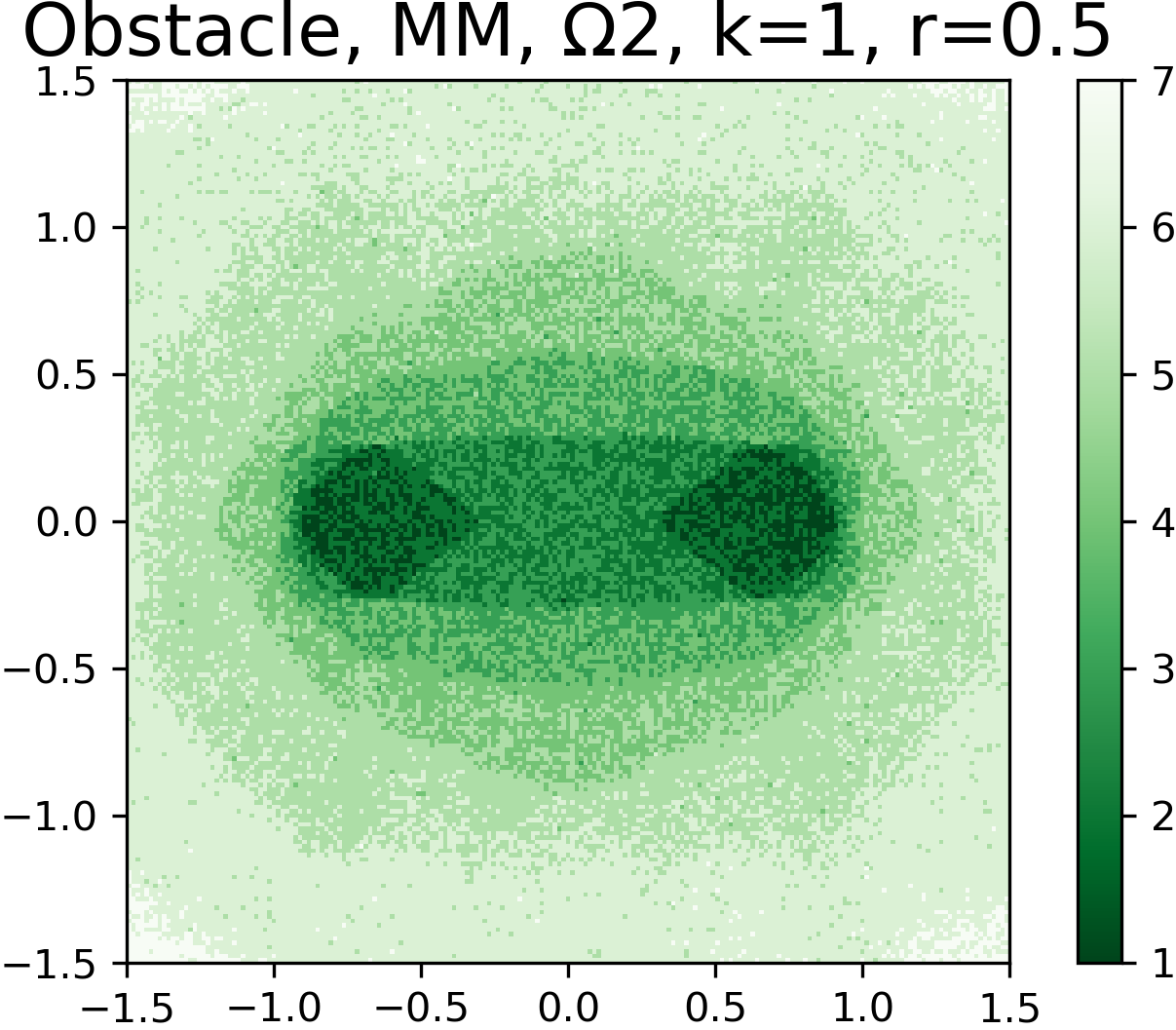}
  \end{center}
 \end{minipage}
 \begin{minipage}{0.5\hsize}
 \begin{center}
  \includegraphics[scale=0.5]{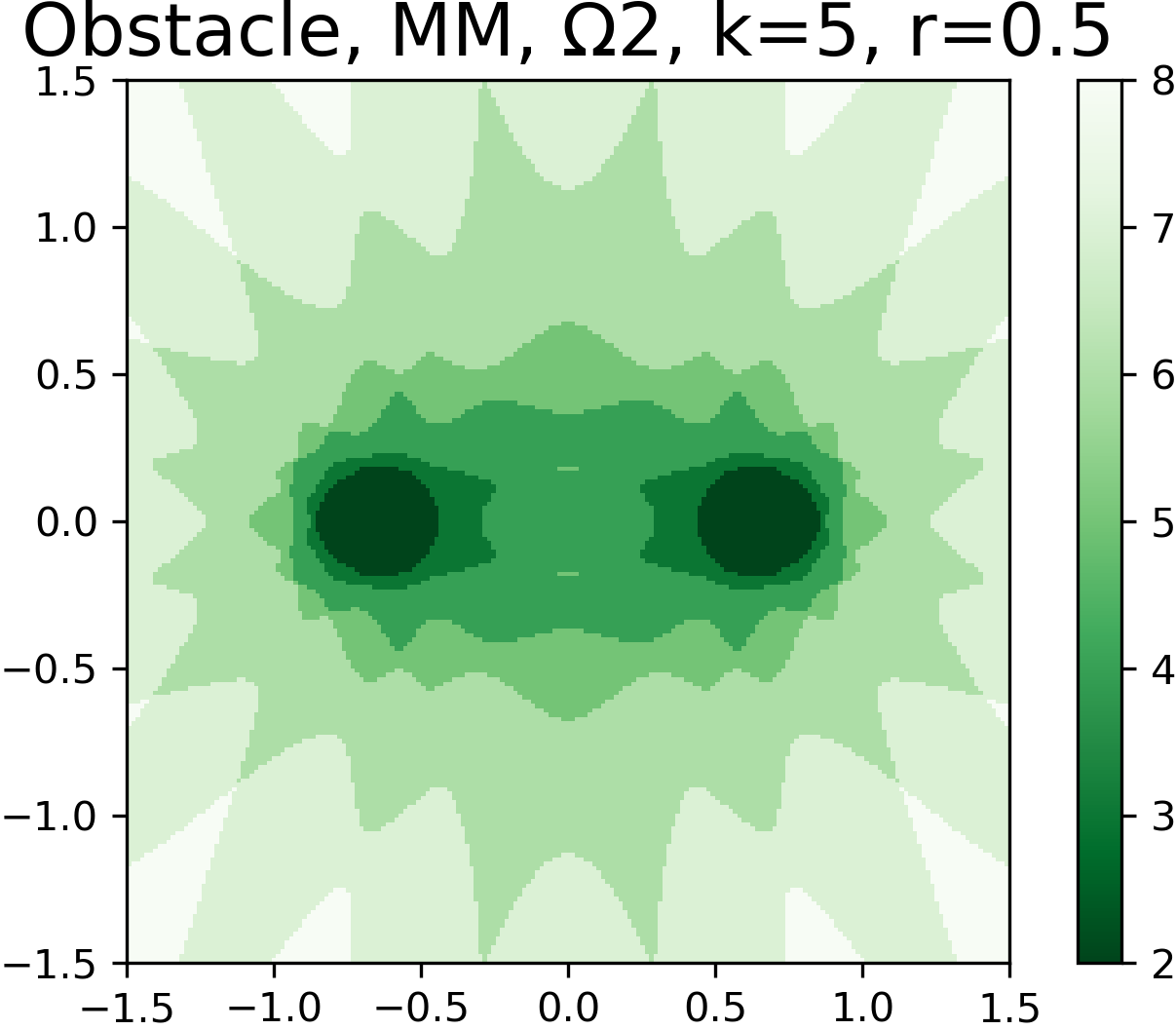}
 \end{center}
 \end{minipage}
 \vspace{1cm} \\ 
\end{tabular}
\caption{Reconstruction for the Dirichlet obstacle by the monotonicity method for different lengths $r=0.1, 0.5$, wavenumbers $k=1,5$, and shapes $\Omega_1, \Omega_2$.}\label{MM obstacle}
\end{figure}
\begin{figure}[htbp]
\vspace{-3cm}
\begin{tabular}{c}
\begin{minipage}{0.5\hsize}
  \begin{center}
   \includegraphics[scale=0.5]{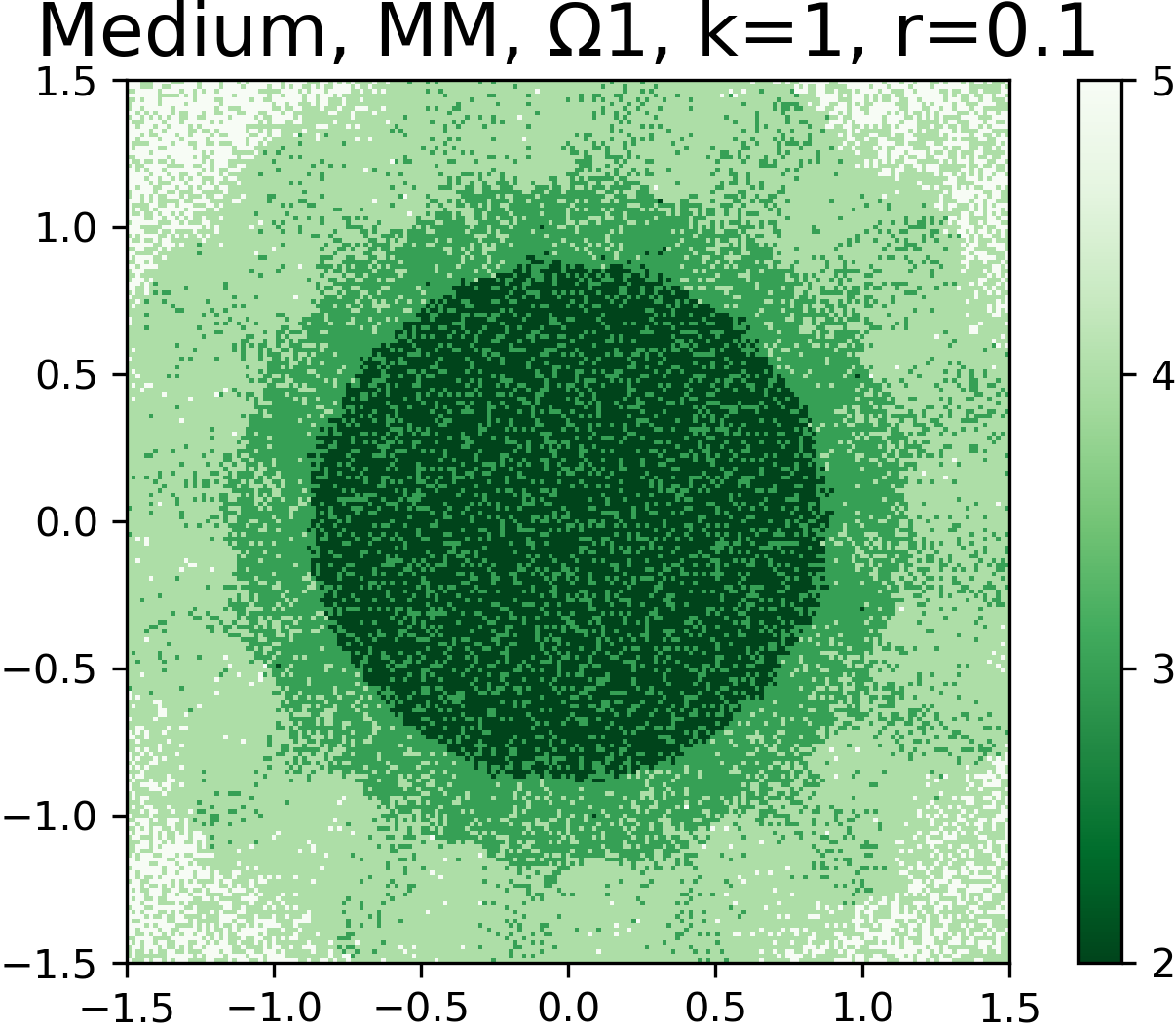}
  \end{center}
 \end{minipage}
 \begin{minipage}{0.5\hsize}
 \begin{center}
  \includegraphics[scale=0.5]{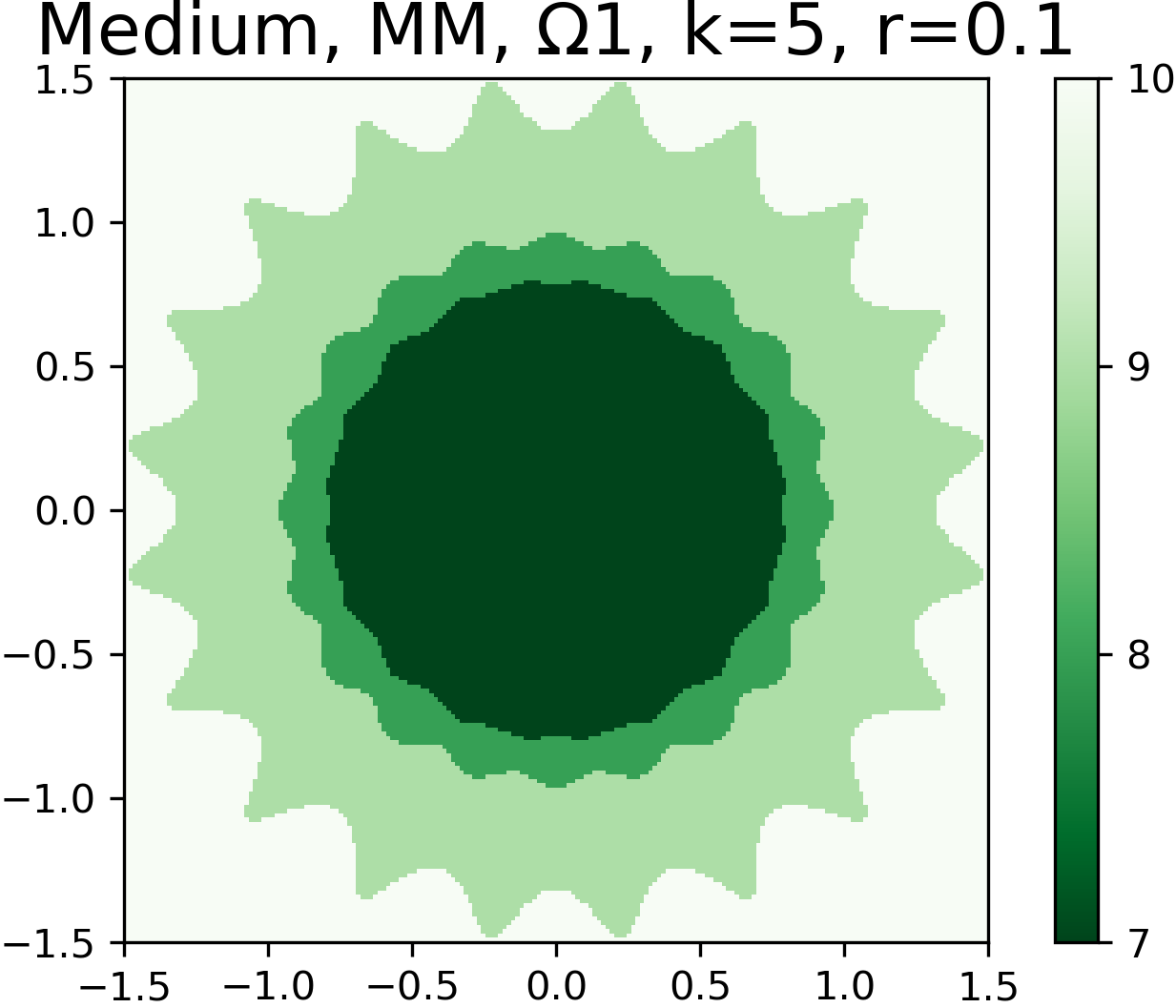}
 \end{center}
 \end{minipage}
 \vspace{1cm} \\ 
\begin{minipage}{0.5\hsize}
  \begin{center}
   \includegraphics[scale=0.5]{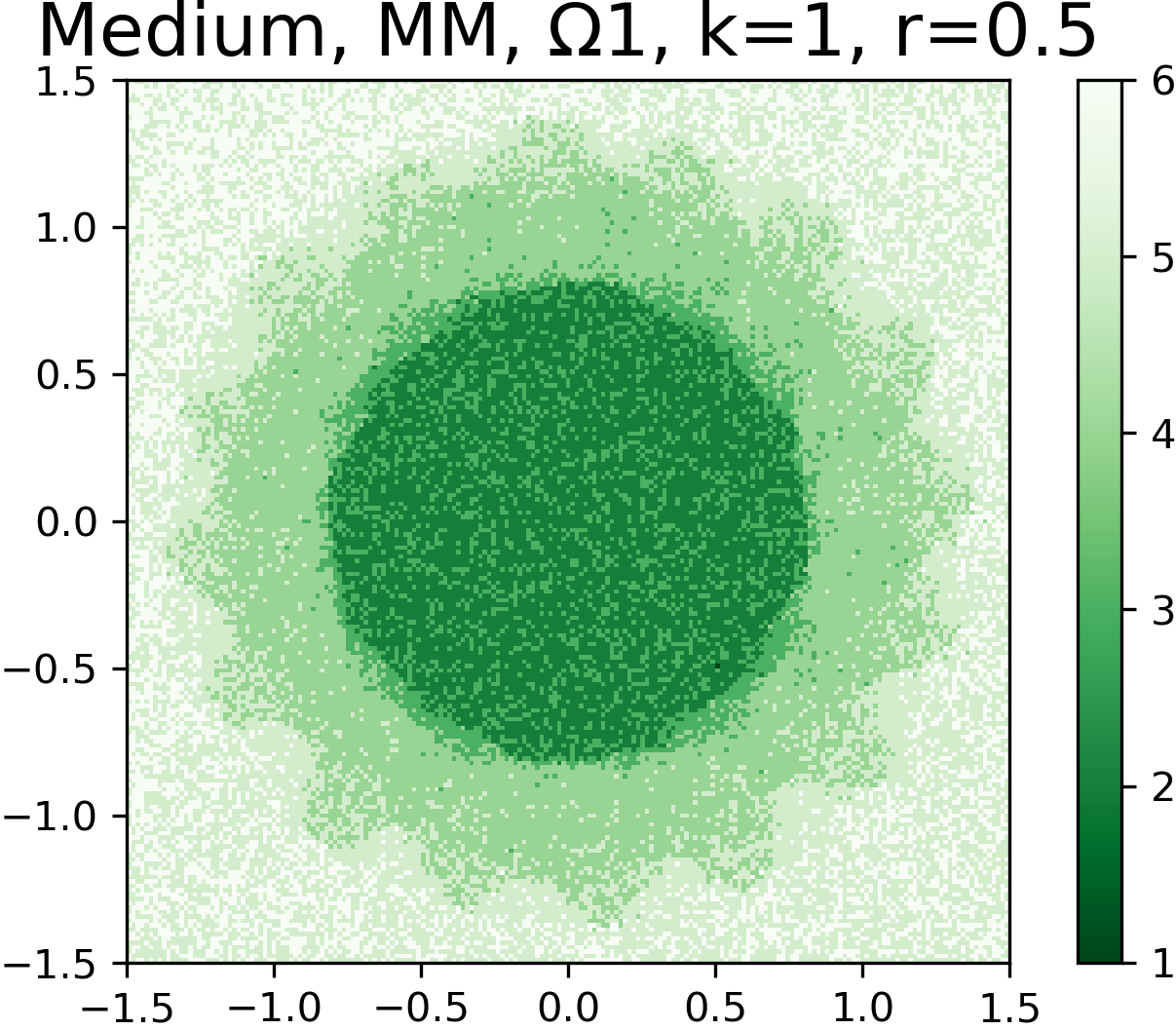}
  \end{center}
 \end{minipage}
 \begin{minipage}{0.5\hsize}
 \begin{center}
  \includegraphics[scale=0.5]{figs/MMmedt1k1r05}
 \end{center}
 \end{minipage}
 \vspace{1cm} \\ 
\begin{minipage}{0.5\hsize}
  \begin{center}
   \includegraphics[scale=0.5]{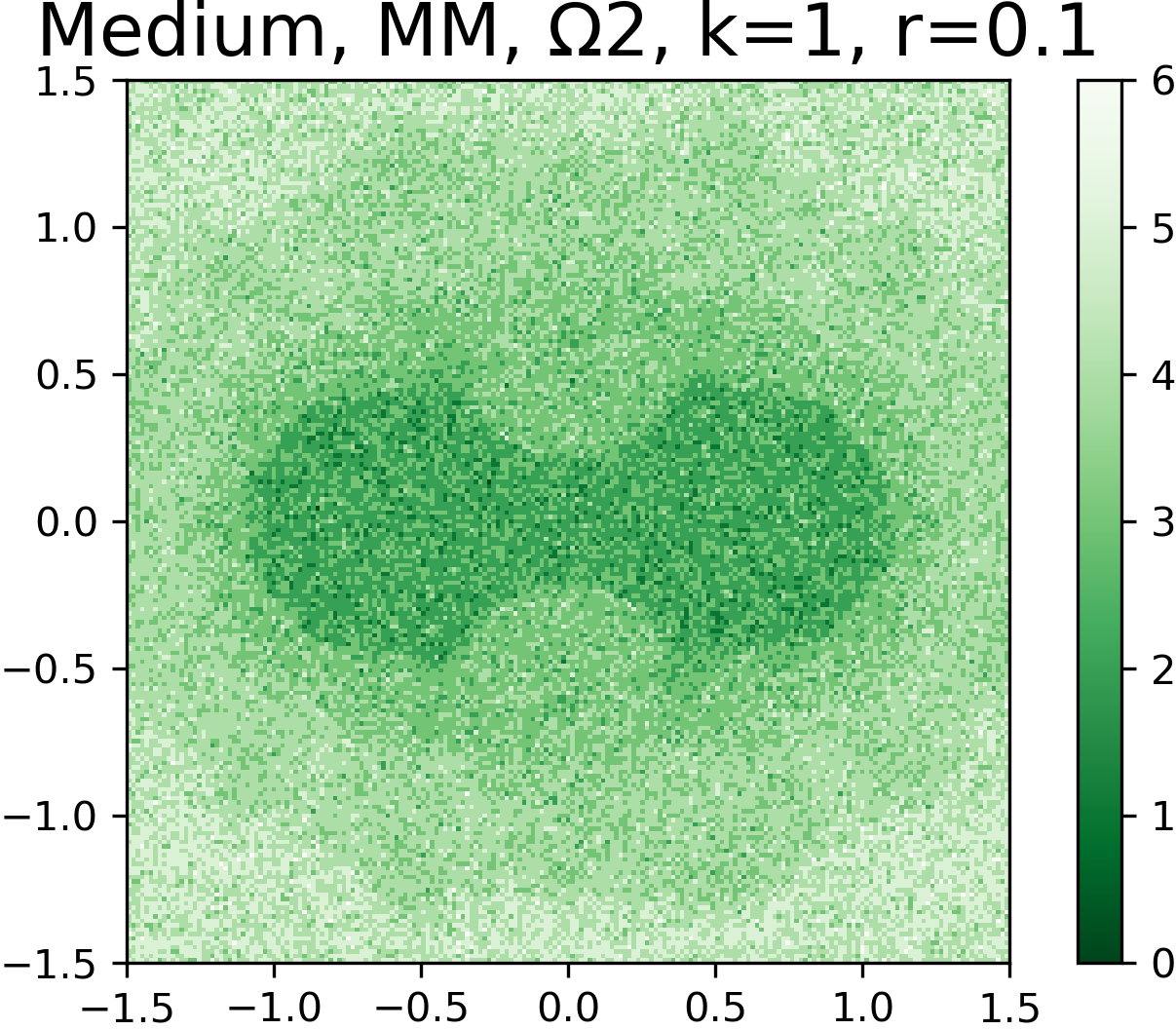}
  \end{center}
 \end{minipage}
 \begin{minipage}{0.5\hsize}
 \begin{center}
  \includegraphics[scale=0.5]{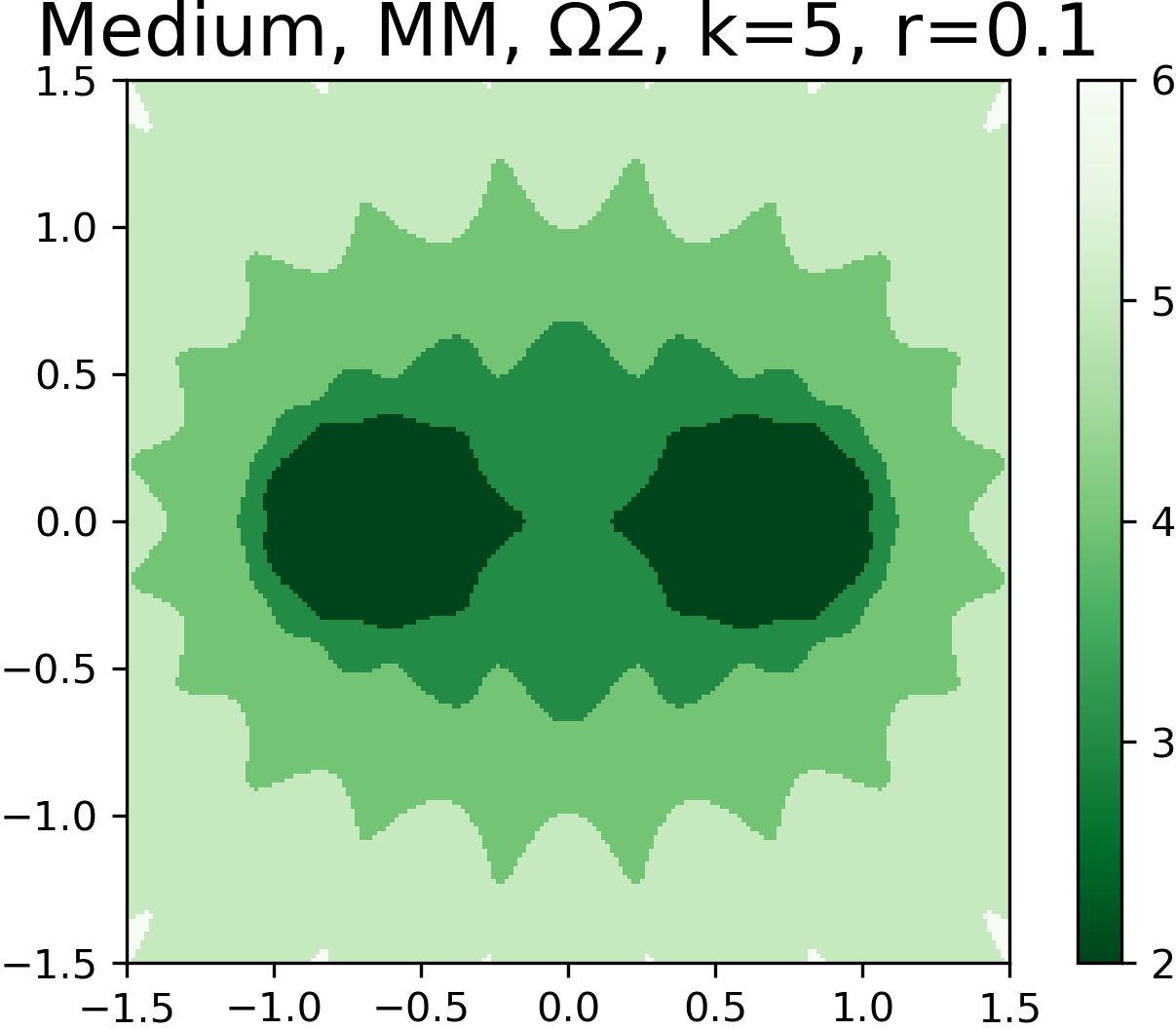}
 \end{center}
 \end{minipage}
 \vspace{1cm} \\ 
\begin{minipage}{0.5\hsize}
  \begin{center}
   \includegraphics[scale=0.5]{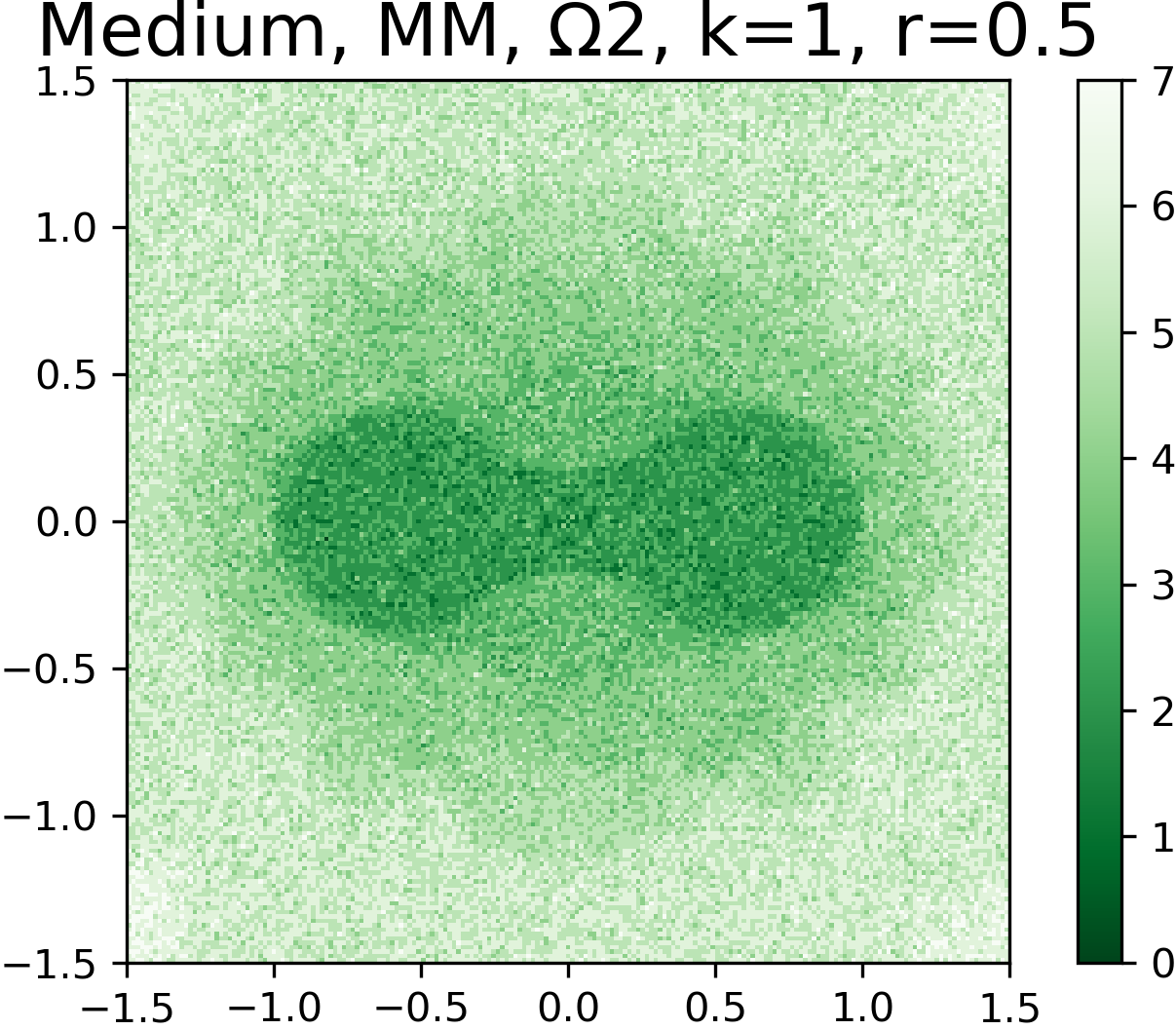}
  \end{center}
 \end{minipage}
 \begin{minipage}{0.5\hsize}
 \begin{center}
  \includegraphics[scale=0.5]{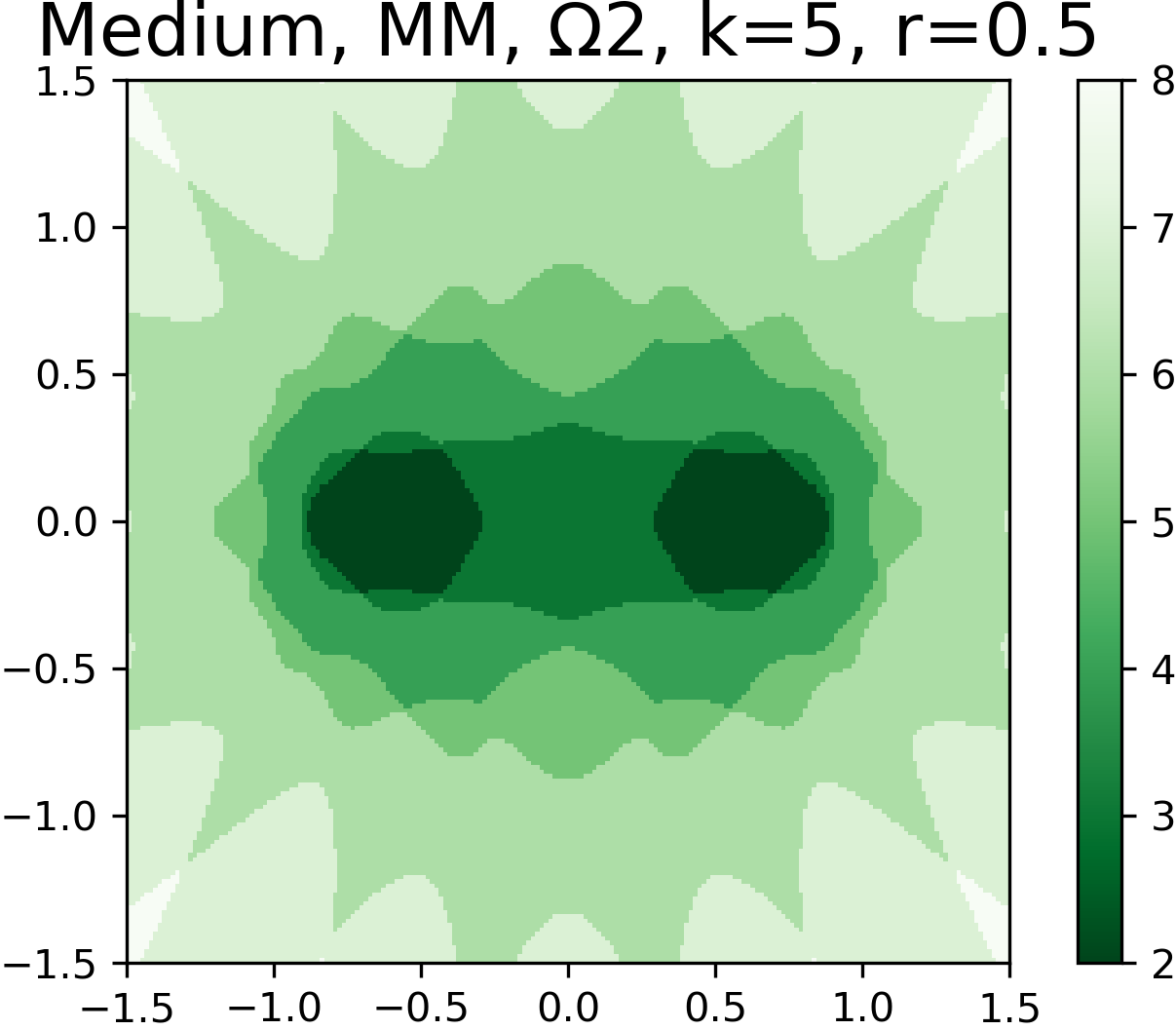}
 \end{center}
 \end{minipage}
 \vspace{1cm} \\ 
\end{tabular}
\caption{Reconstruction for the inhomogeneous medium by the monotonicity method for different lengths $r=0.1, 0.5$, wavenumbers $k=1,5$, and shapes $\Omega_1, \Omega_2$}\label{MM medium}
\end{figure}
\begin{figure}[htbp]
\begin{tabular}{c}
\begin{minipage}{0.5\hsize}
  \begin{center}
   \includegraphics[scale=0.5]{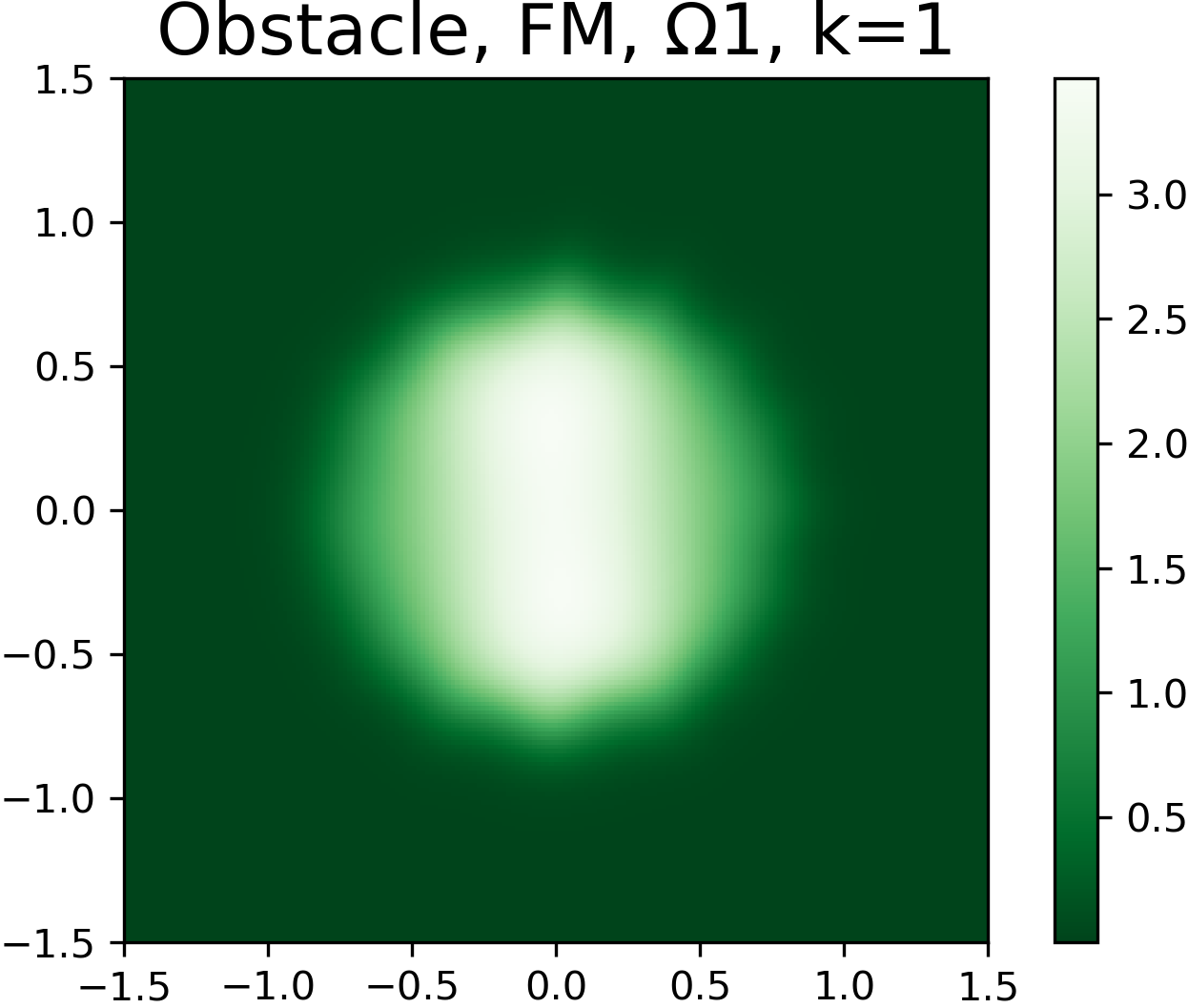}
  \end{center}
 \end{minipage}
 \begin{minipage}{0.5\hsize}
 \begin{center}
  \includegraphics[scale=0.5]{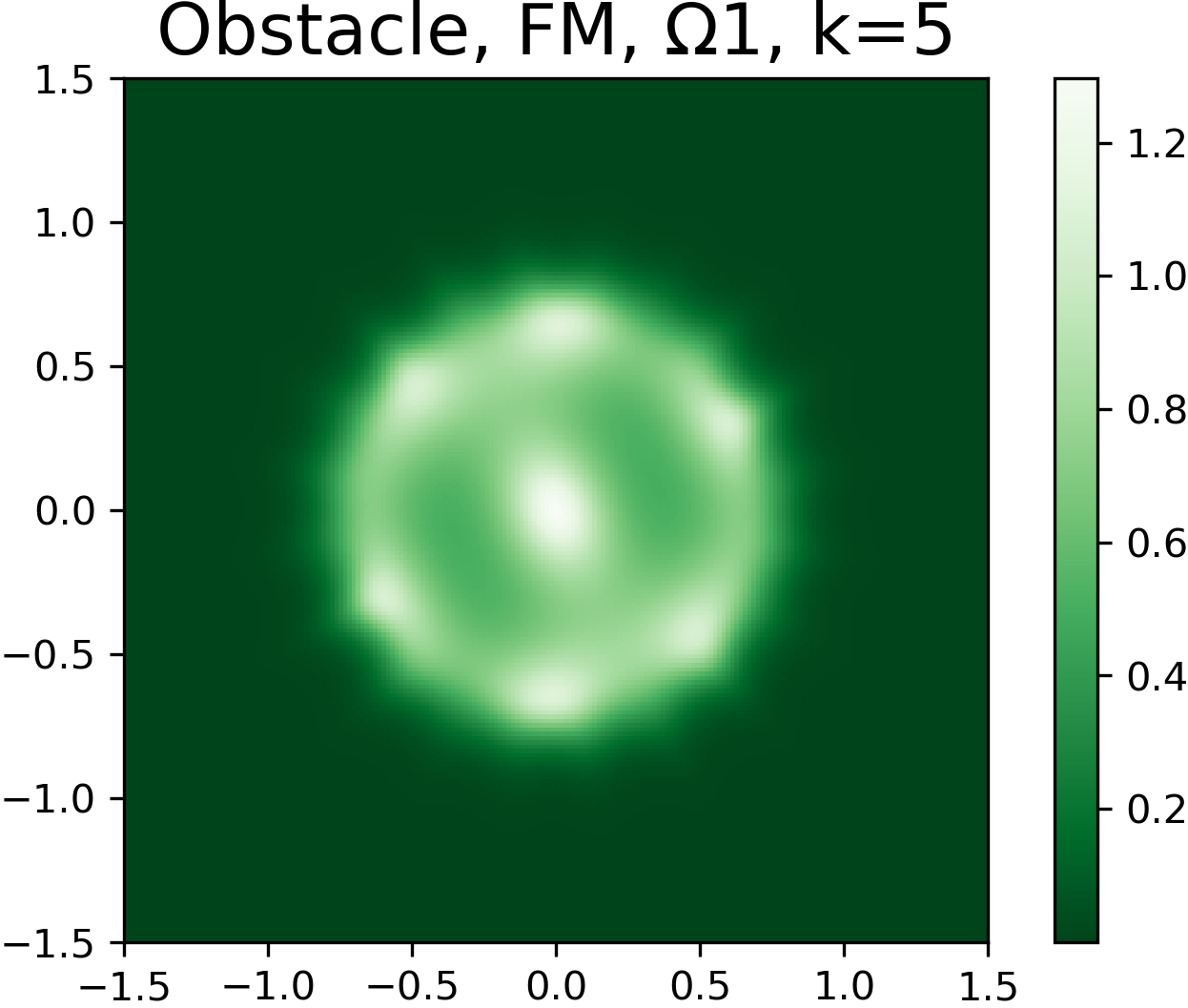}
 \end{center}
 \end{minipage}
 \vspace{1cm} \\ 
\begin{minipage}{0.5\hsize}
  \begin{center}
   \includegraphics[scale=0.5]{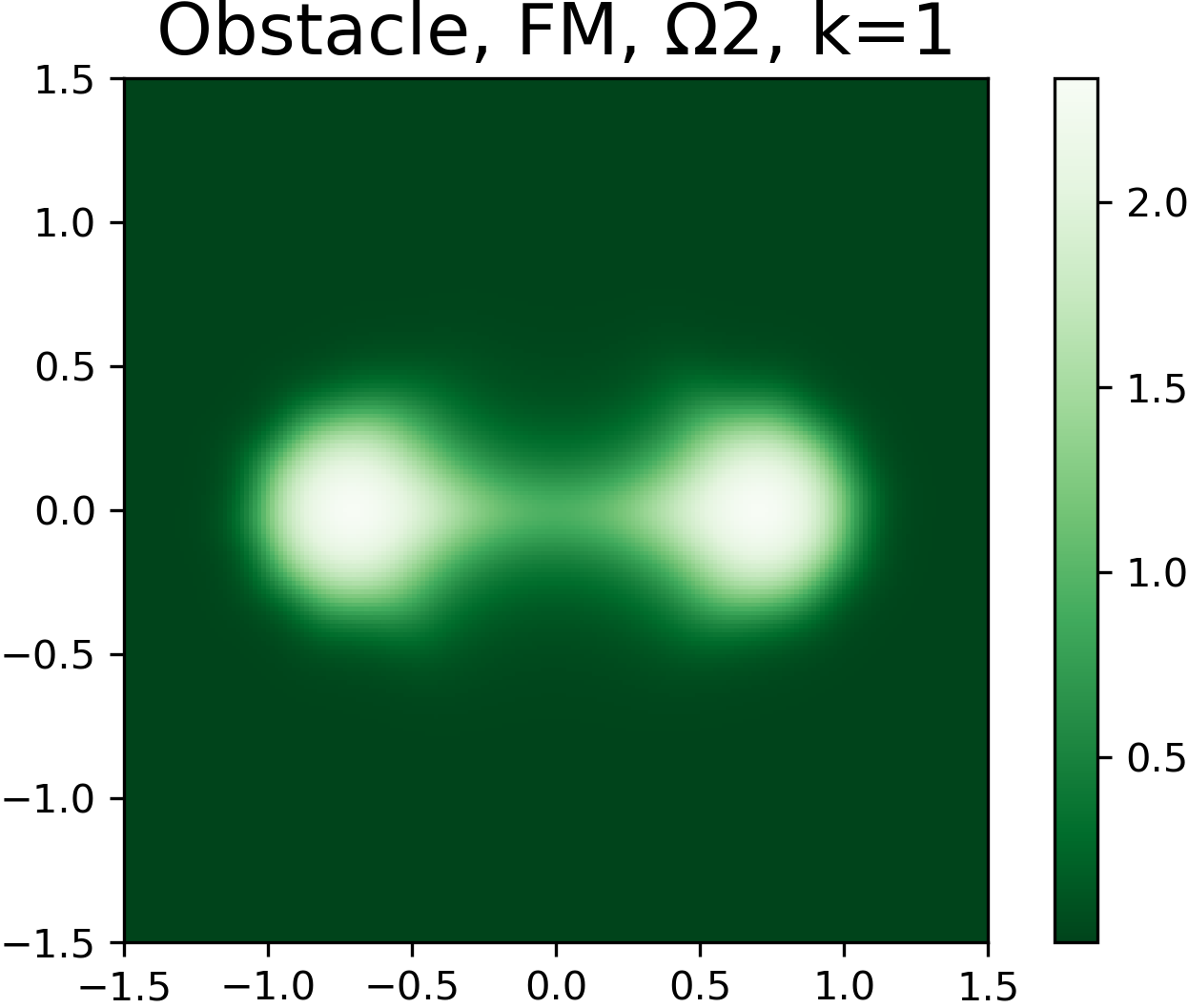}
  \end{center}
 \end{minipage}
 \begin{minipage}{0.5\hsize}
 \begin{center}
  \includegraphics[scale=0.5]{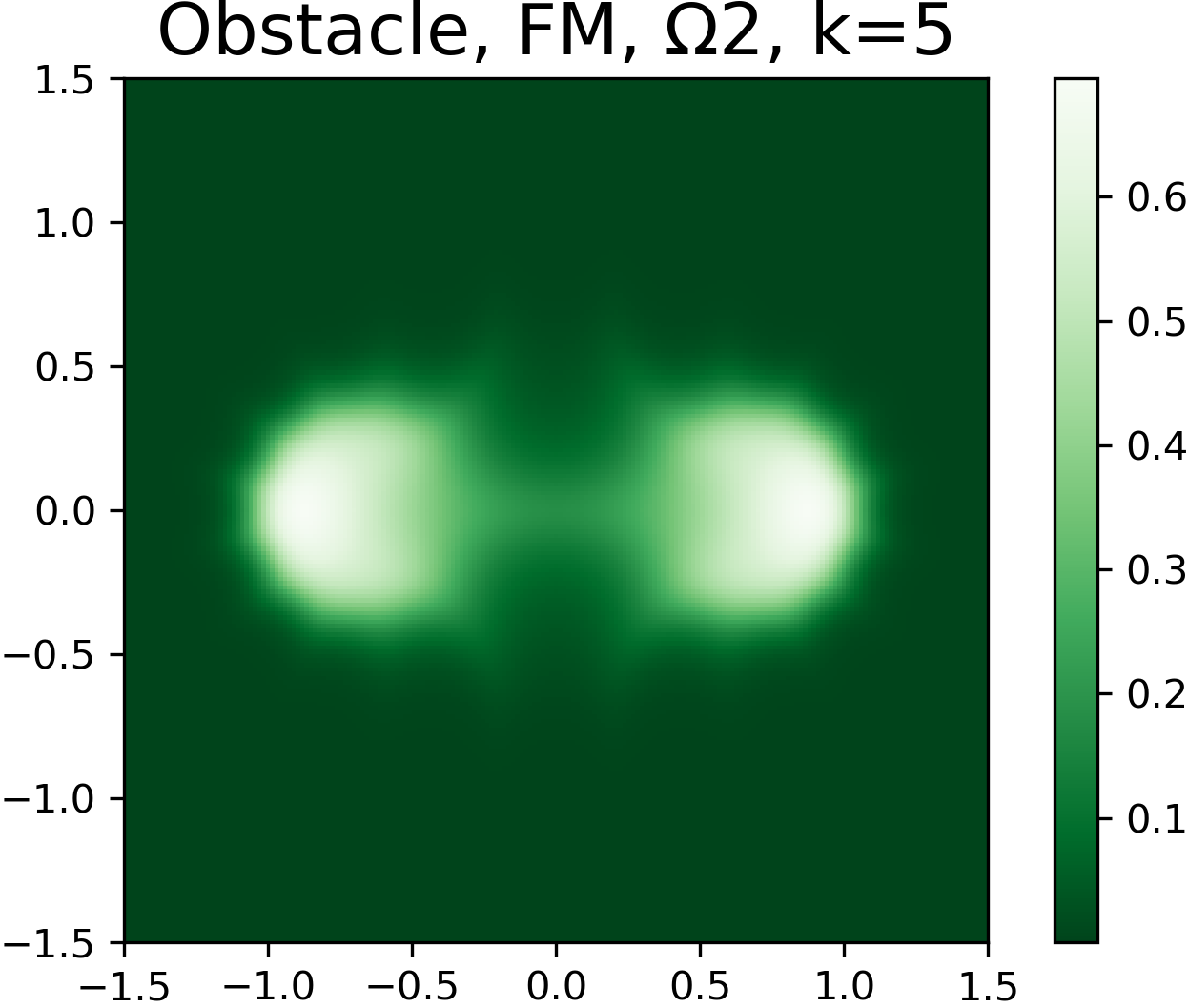}
 \end{center}
 \end{minipage}
\end{tabular}
\caption{Reconstruction for the Dirichlet obstacle by the factorization method for different wavenumbers $k=1,5$ and shapes $\Omega_1, \Omega_2$.}\label{FM obstacle}
\end{figure}
\begin{figure}[htbp]
\begin{tabular}{c}
\begin{minipage}{0.5\hsize}
  \begin{center}
   \includegraphics[scale=0.5]{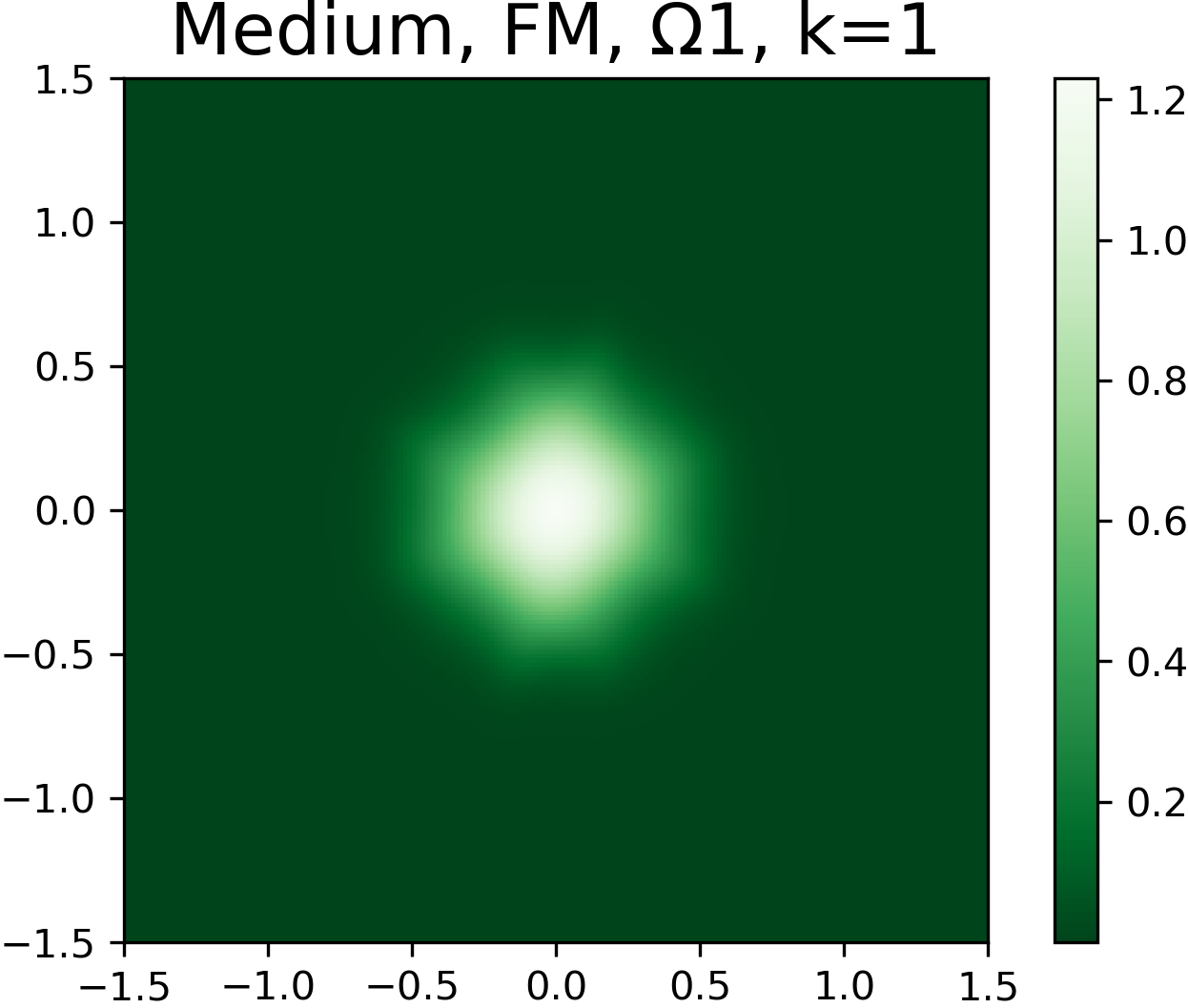}
  \end{center}
 \end{minipage}
 \begin{minipage}{0.5\hsize}
 \begin{center}
  \includegraphics[scale=0.5]{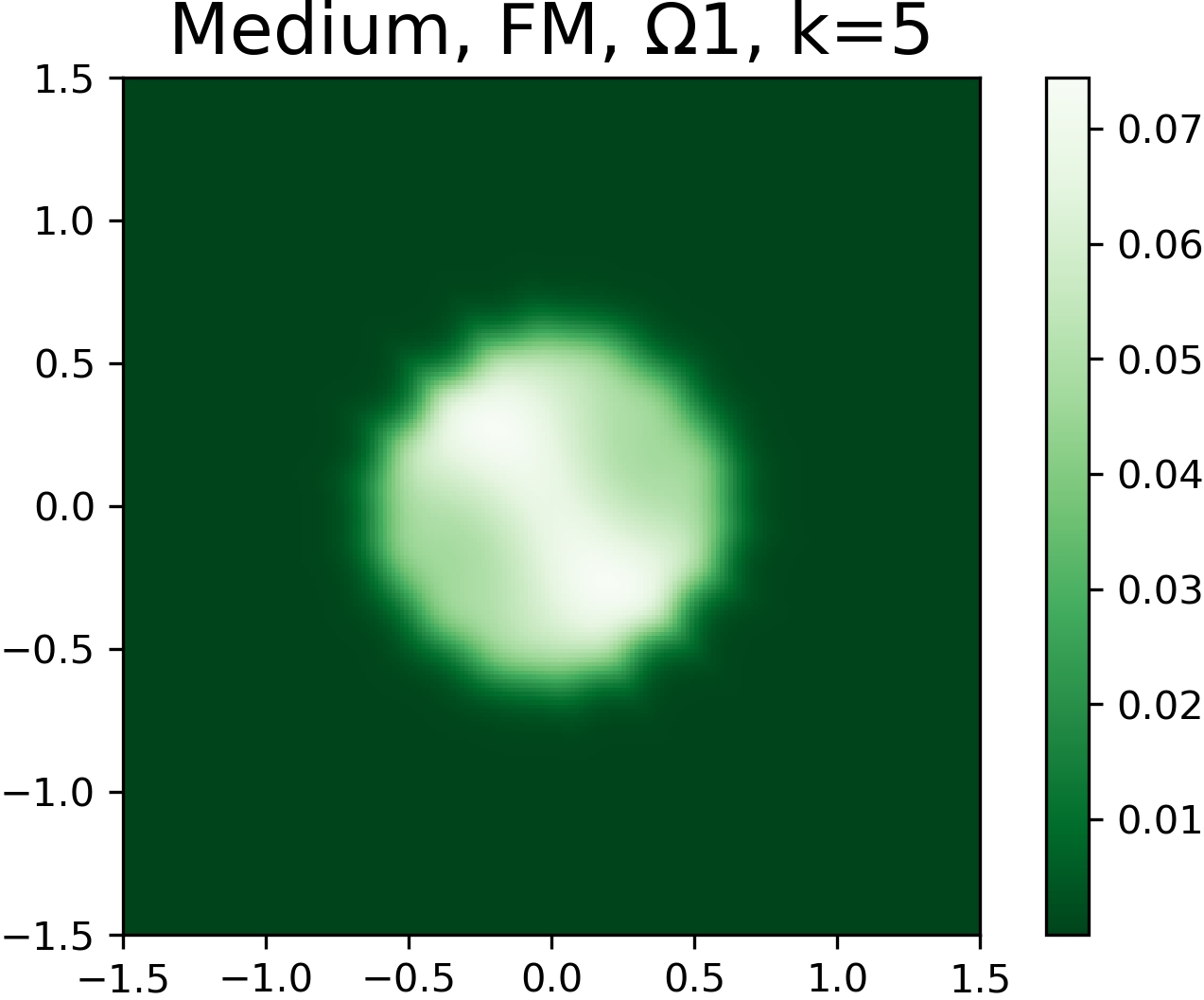}
 \end{center}
 \end{minipage}
 \vspace{1cm} \\ 
\begin{minipage}{0.5\hsize}
  \begin{center}
   \includegraphics[scale=0.5]{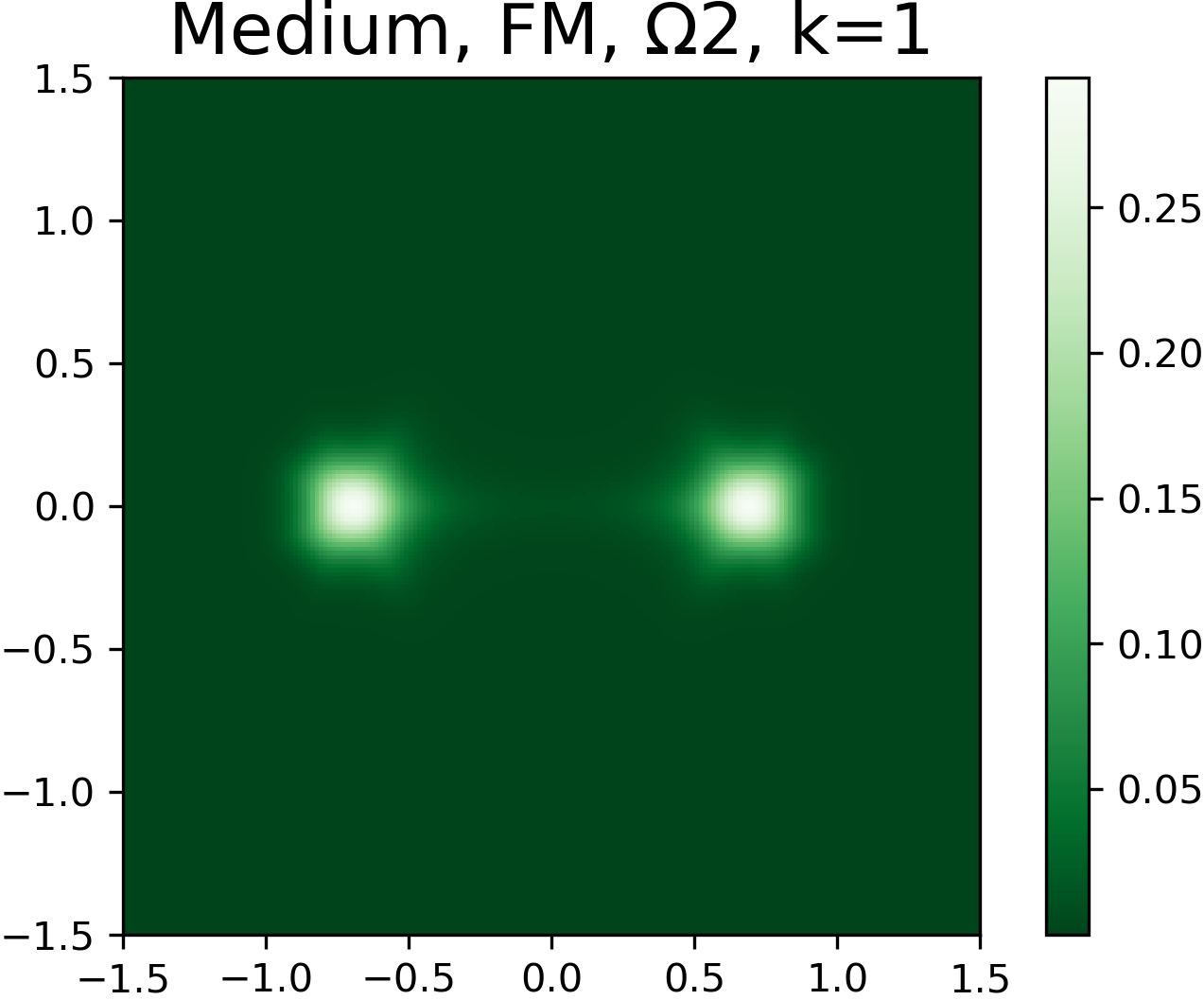}
  \end{center}
 \end{minipage}
 \begin{minipage}{0.5\hsize}
 \begin{center}
  \includegraphics[scale=0.5]{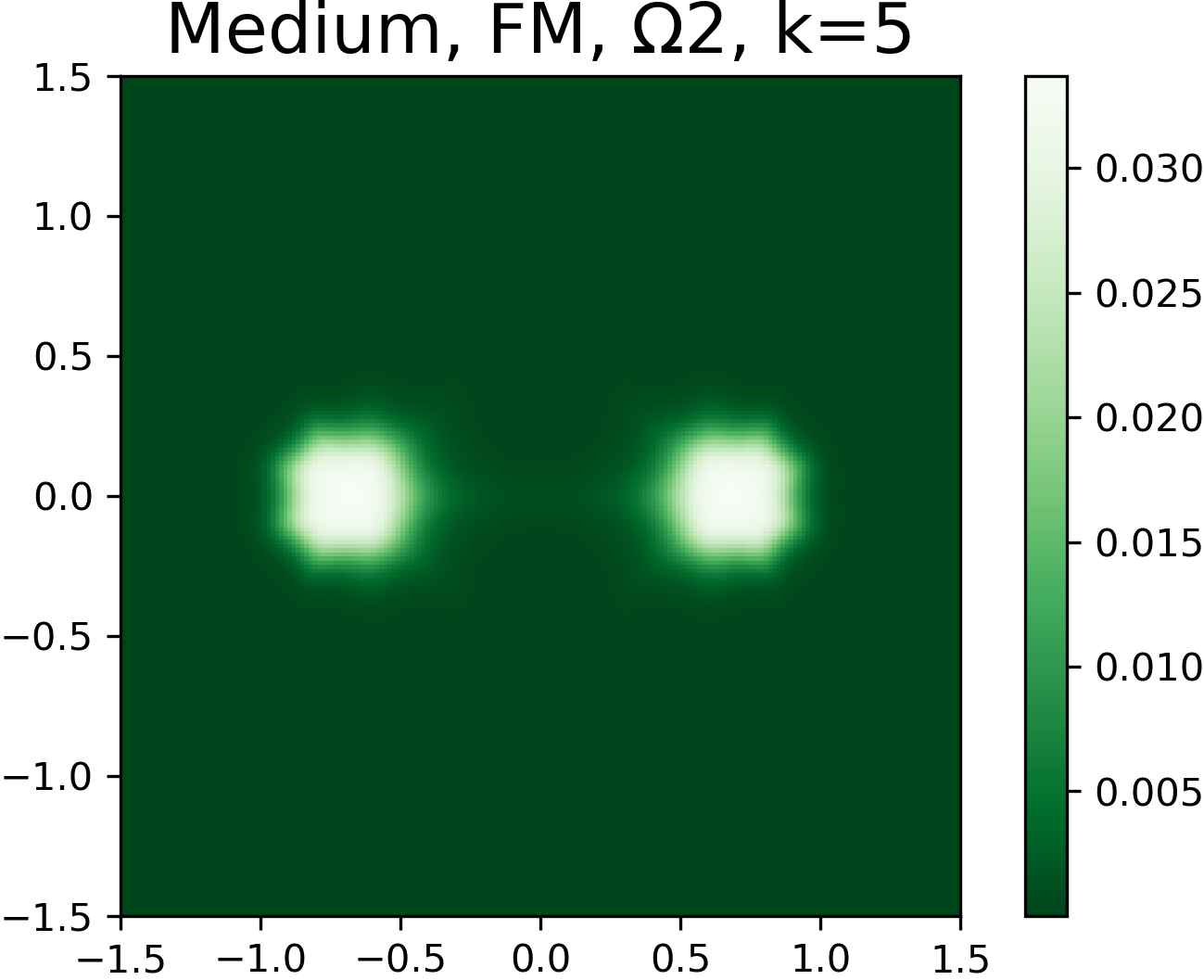}
 \end{center}
 \end{minipage}
\end{tabular}
\caption{Reconstruction for the inhomogeneous medium by the factorization method for different wavenumbers $k=1,5$ and shapes $\Omega_1, \Omega_2$.}\label{FM medium}
\end{figure}
\begin{figure}[htbp]
\vspace{-3cm}
\begin{tabular}{c}
\begin{minipage}{0.5\hsize}
  \begin{center}
   \includegraphics[scale=0.5]{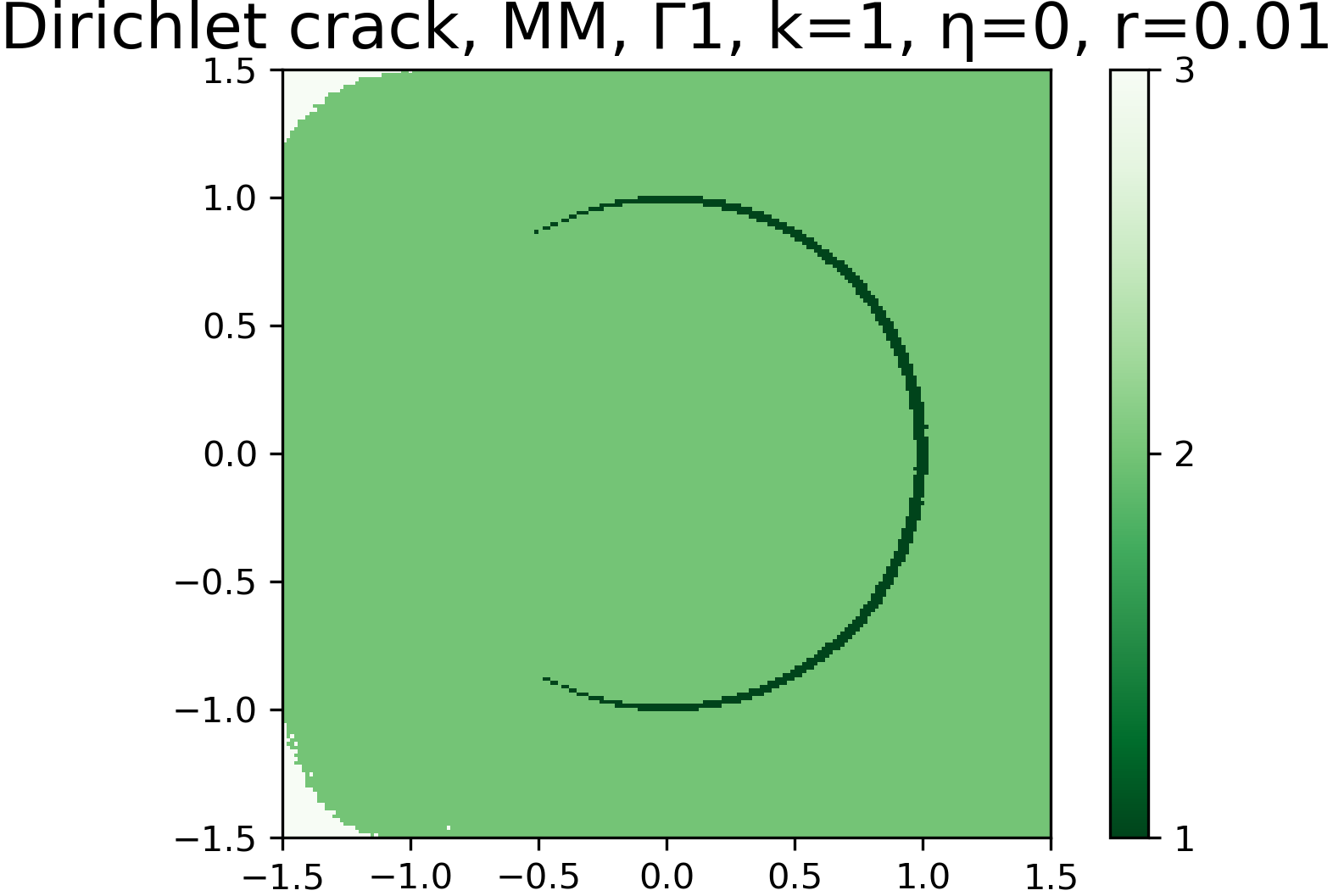}
  \end{center}
 \end{minipage}
 \begin{minipage}{0.5\hsize}
 \begin{center}
  \includegraphics[scale=0.5]{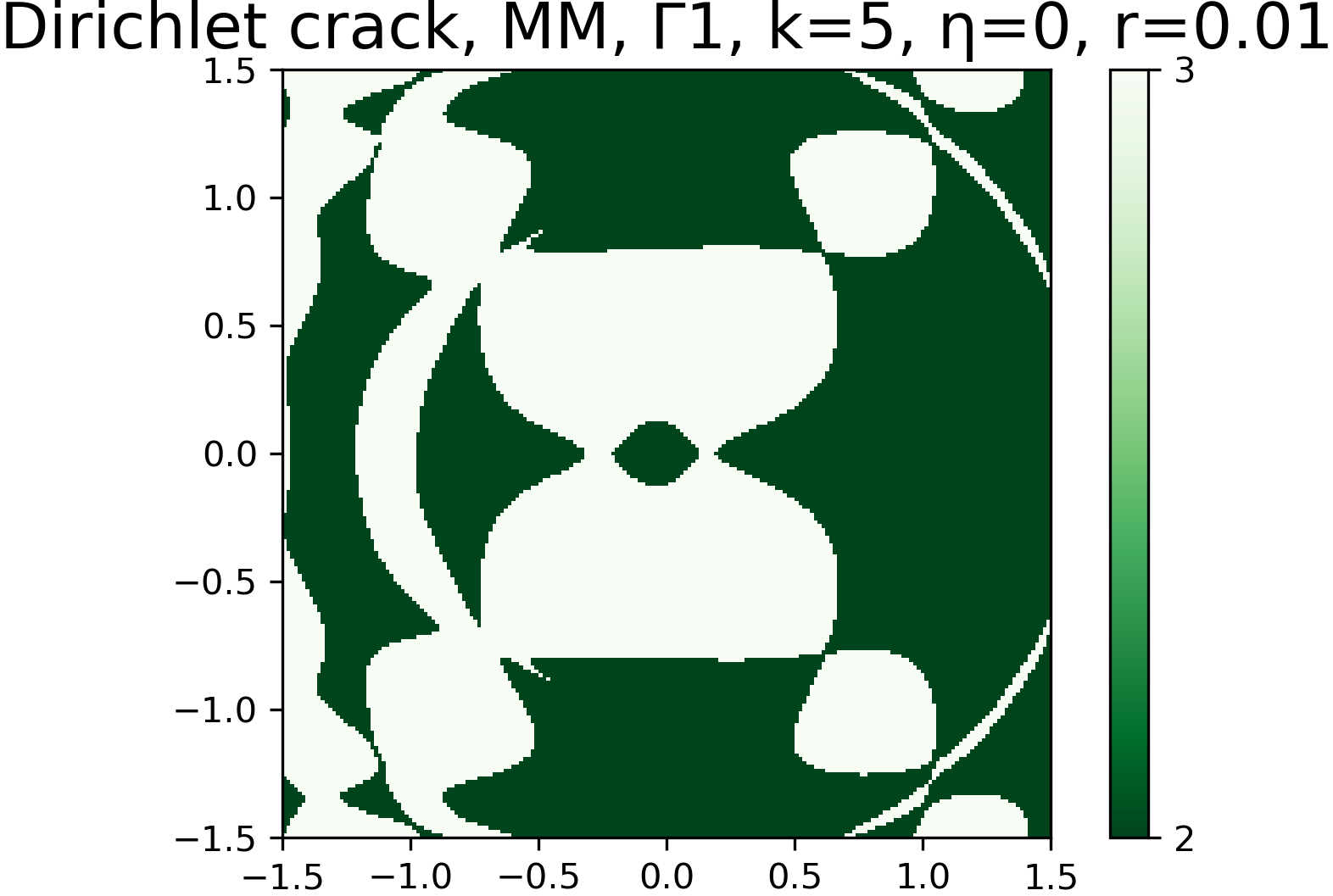}
 \end{center}
 \end{minipage}
 \vspace{1cm} \\ 
\begin{minipage}{0.5\hsize}
  \begin{center}
   \includegraphics[scale=0.5]{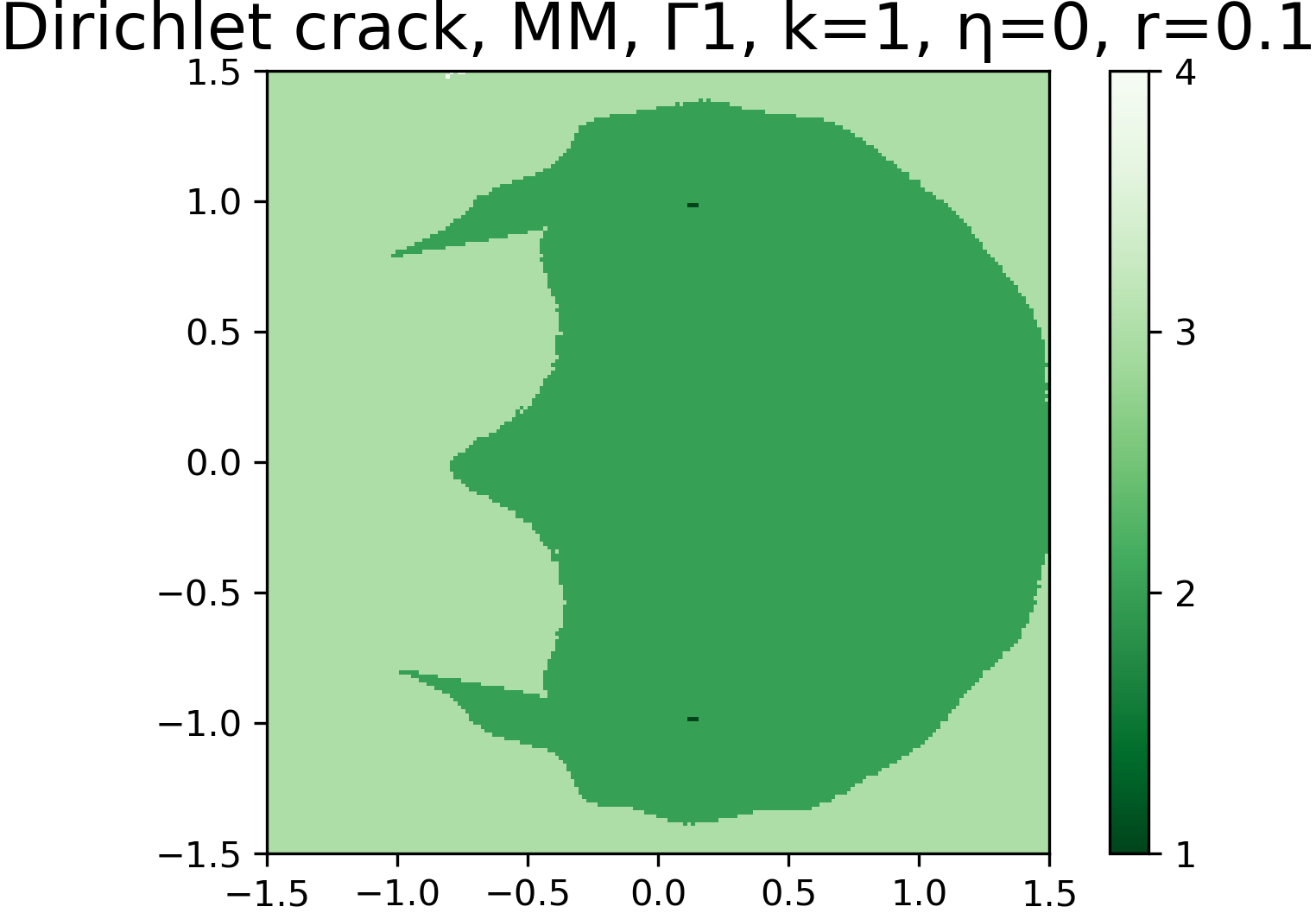}
  \end{center}
 \end{minipage}
 \begin{minipage}{0.5\hsize}
 \begin{center}
  \includegraphics[scale=0.5]{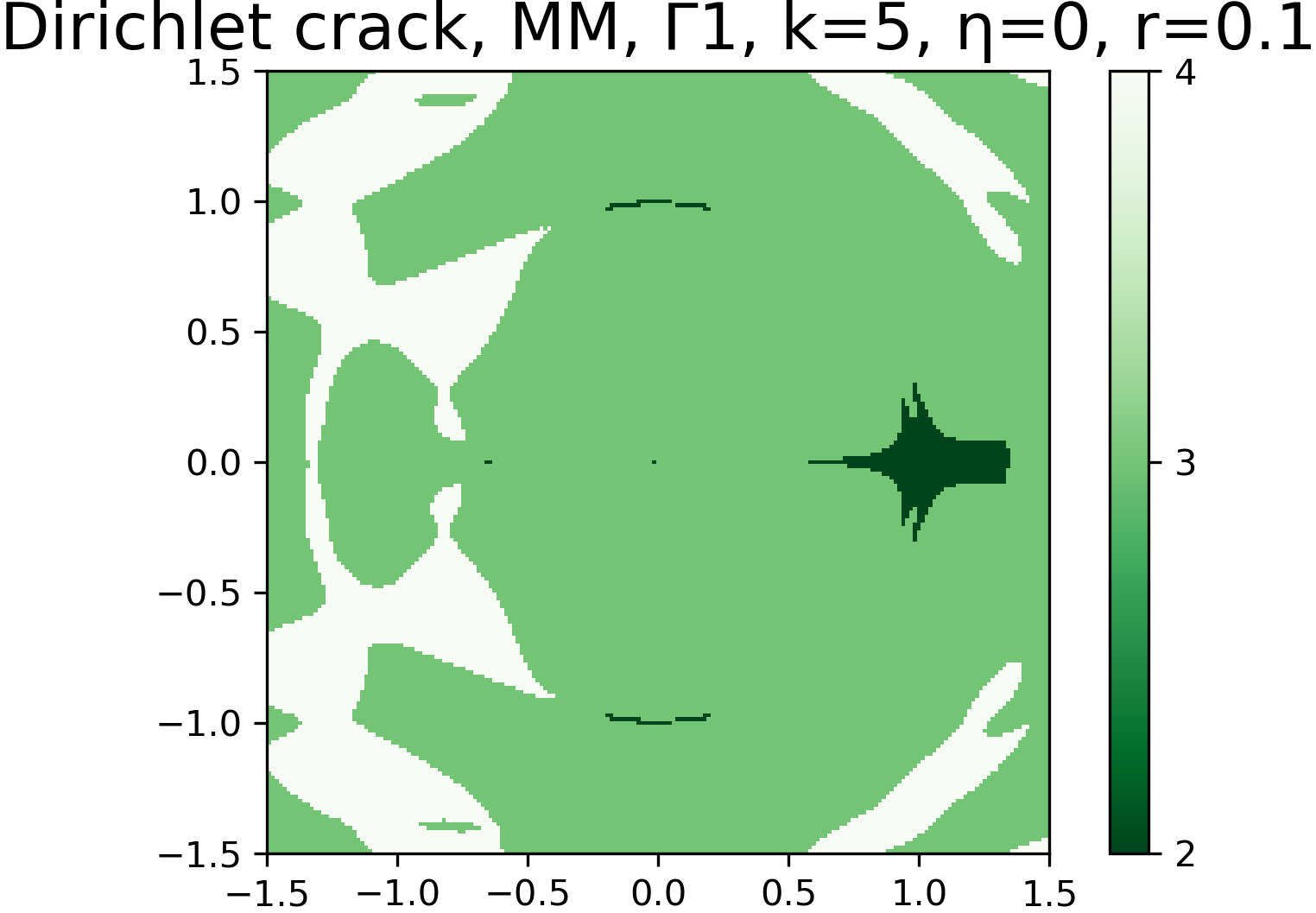}
 \end{center}
 \end{minipage}
 \vspace{1cm} \\ 
\begin{minipage}{0.5\hsize}
  \begin{center}
   \includegraphics[scale=0.5]{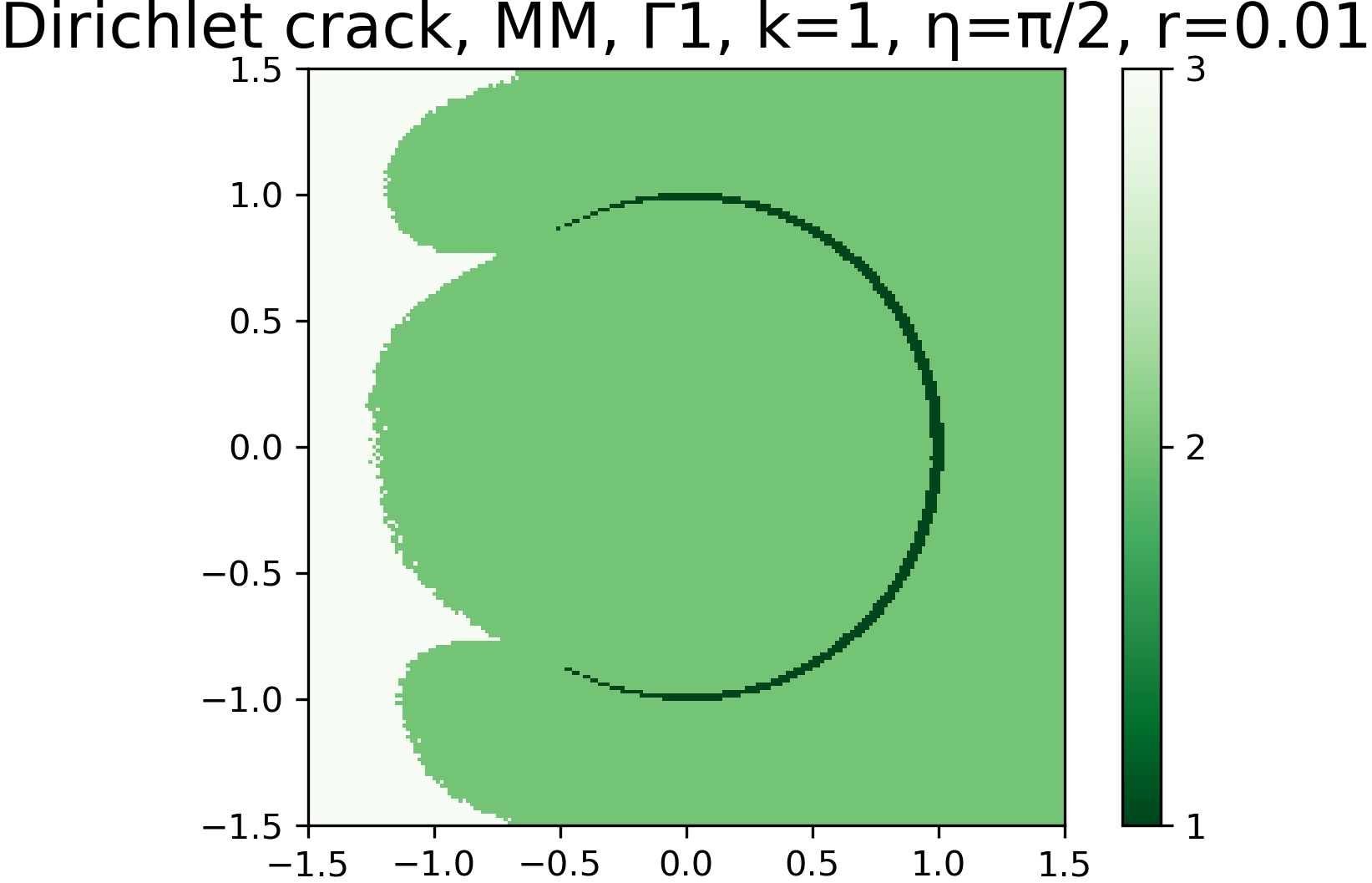}
  \end{center}
 \end{minipage}
 \begin{minipage}{0.5\hsize}
 \begin{center}
  \includegraphics[scale=0.5]{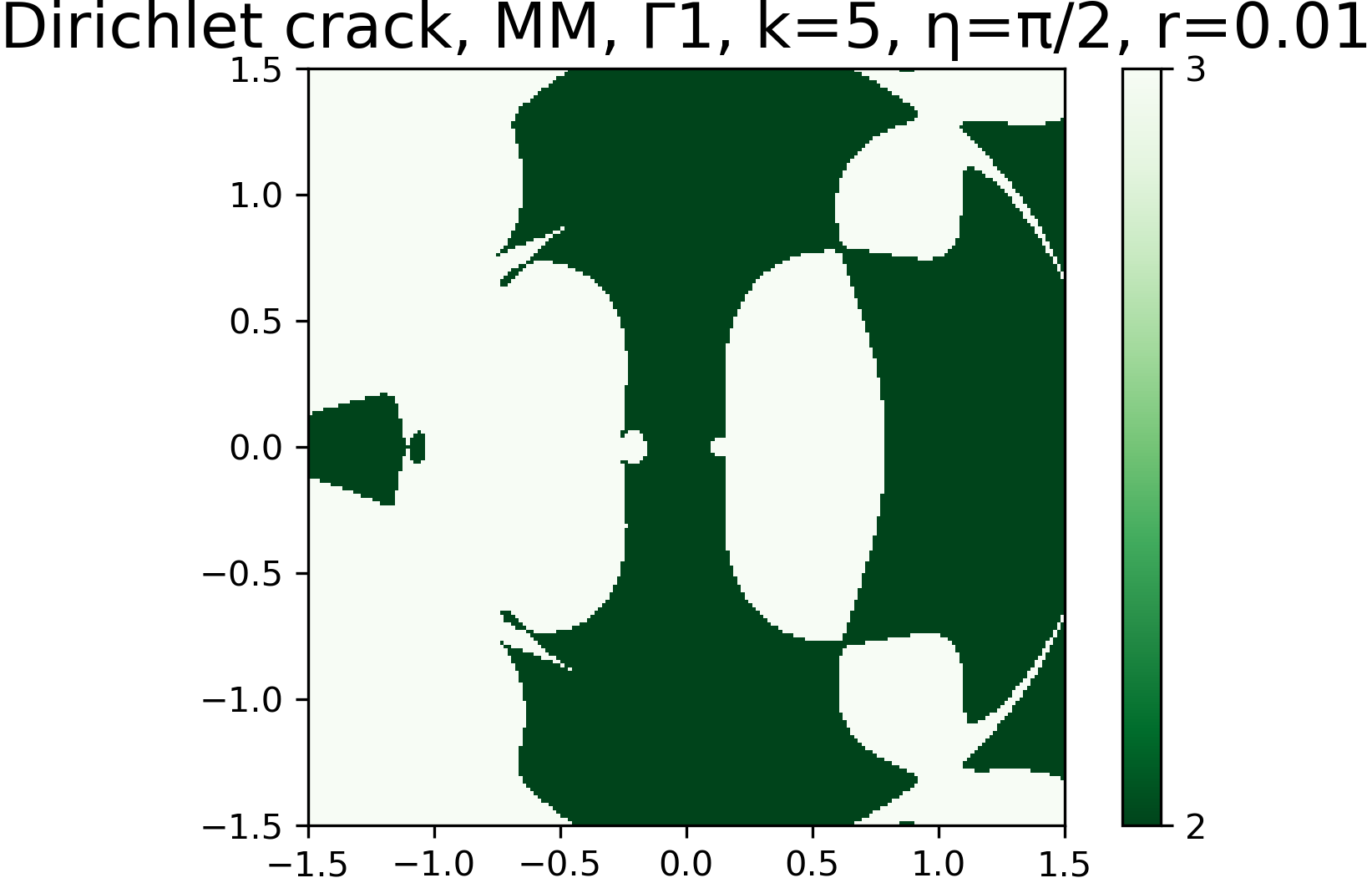}
 \end{center}
 \end{minipage}
 \vspace{1cm} \\ 
\begin{minipage}{0.5\hsize}
  \begin{center}
   \includegraphics[scale=0.5]{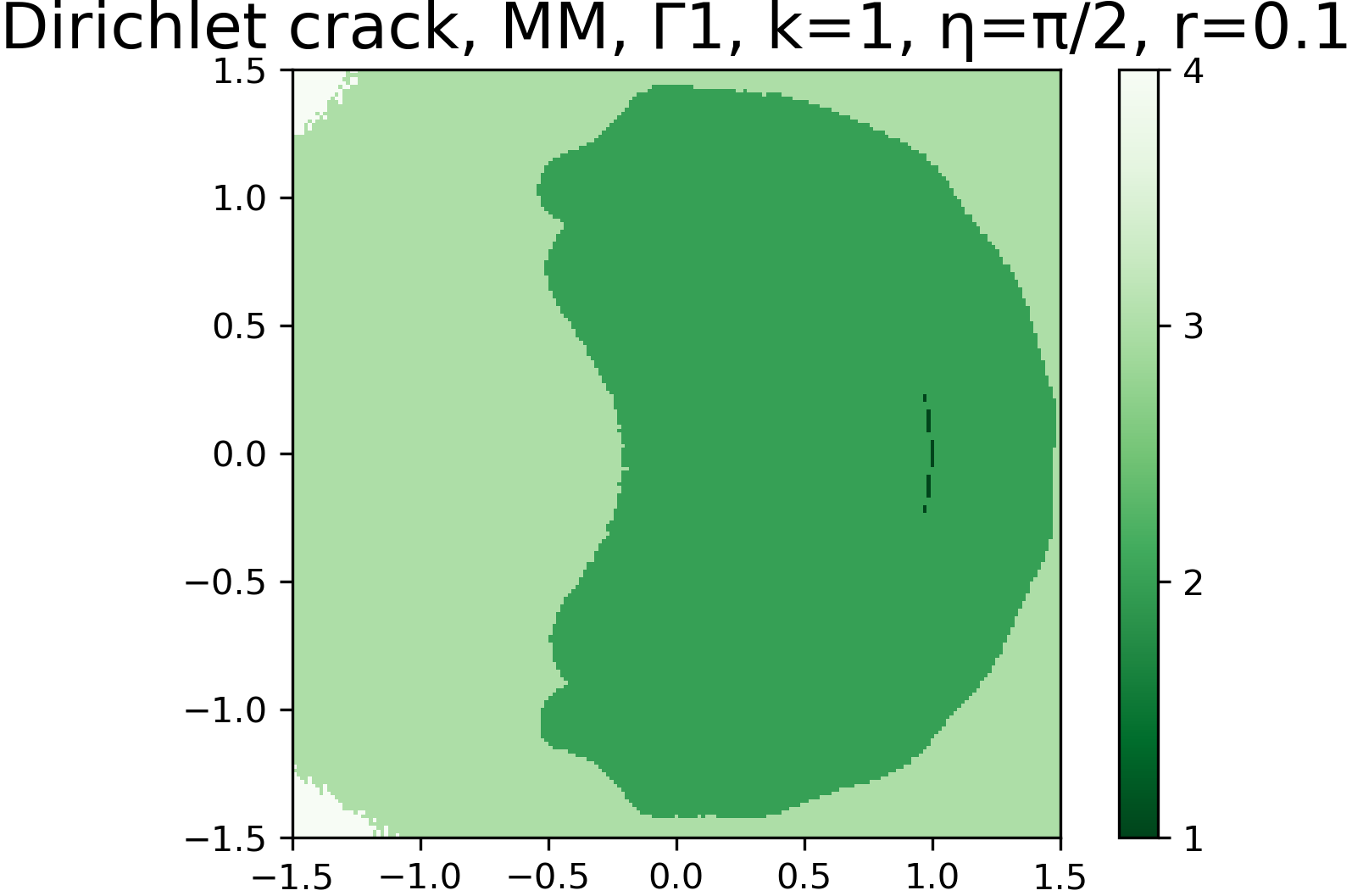}
  \end{center}
 \end{minipage}
 \begin{minipage}{0.5\hsize}
 \begin{center}
  \includegraphics[scale=0.5]{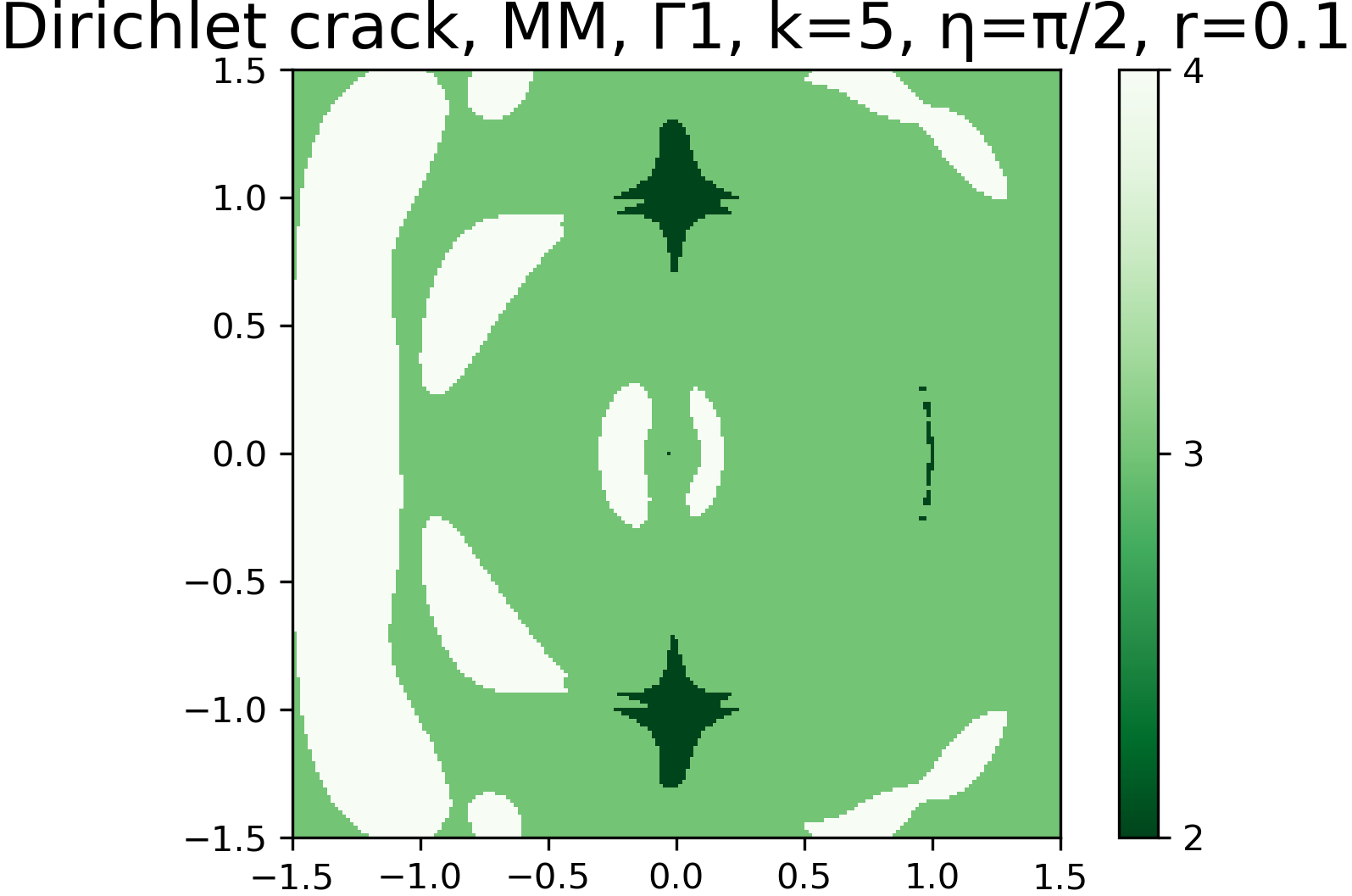}
 \end{center}
 \end{minipage}
 \vspace{1cm} \\ 
\end{tabular}
\caption{Reconstruction for the Dirichlet crack $\Gamma_1$ by the monotonicity method for different angles $\eta=0, \pi/2$, lengths $r=0.01, 0.1$, and wavenumbers $k=1,5$.}\label{MM crack-1}
\end{figure}
\begin{figure}[htbp]
\vspace{-3cm}
\begin{tabular}{c}
\begin{minipage}{0.5\hsize}
  \begin{center}
   \includegraphics[scale=0.5]{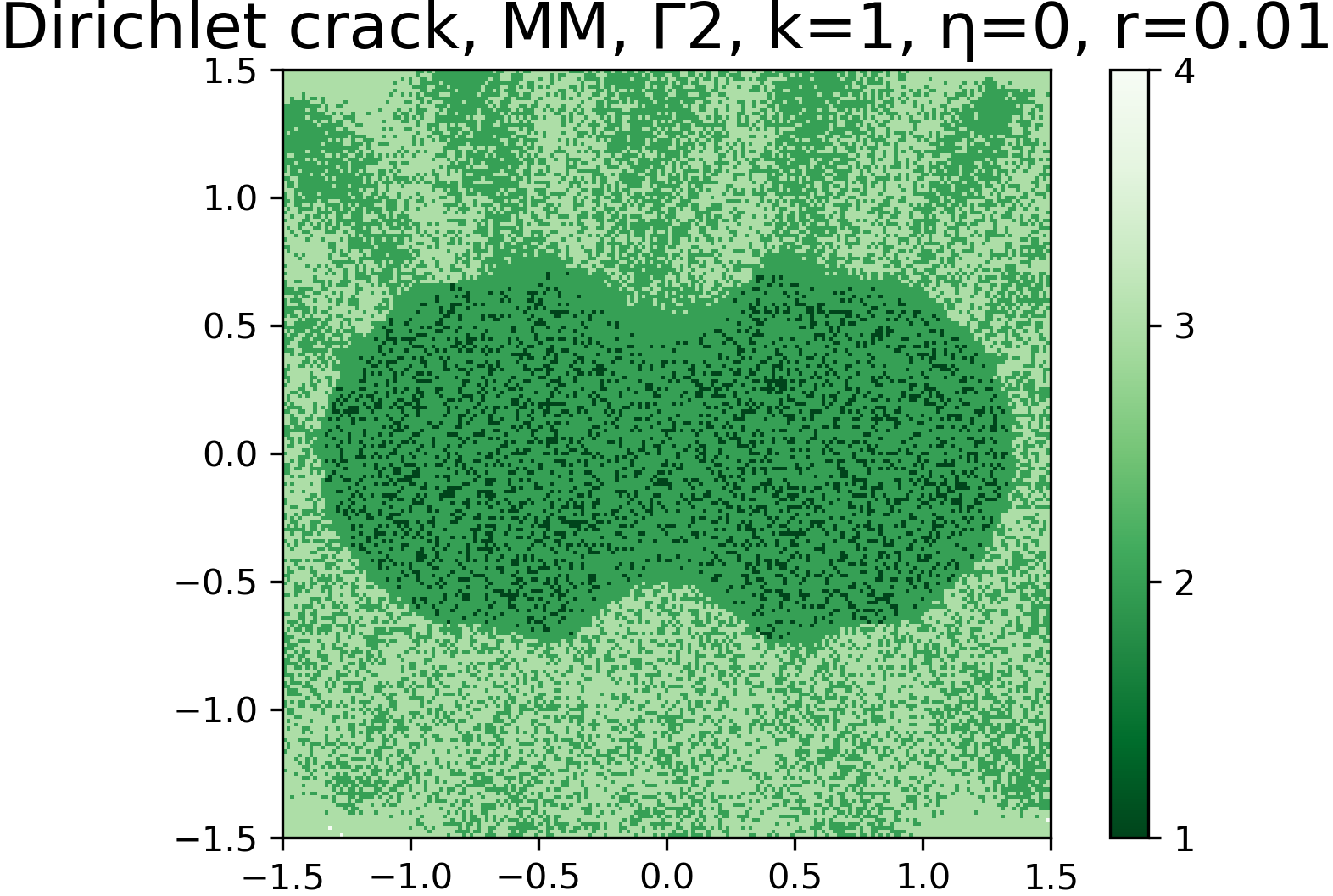}
  \end{center}
 \end{minipage}
 \begin{minipage}{0.5\hsize}
 \begin{center}
  \includegraphics[scale=0.5]{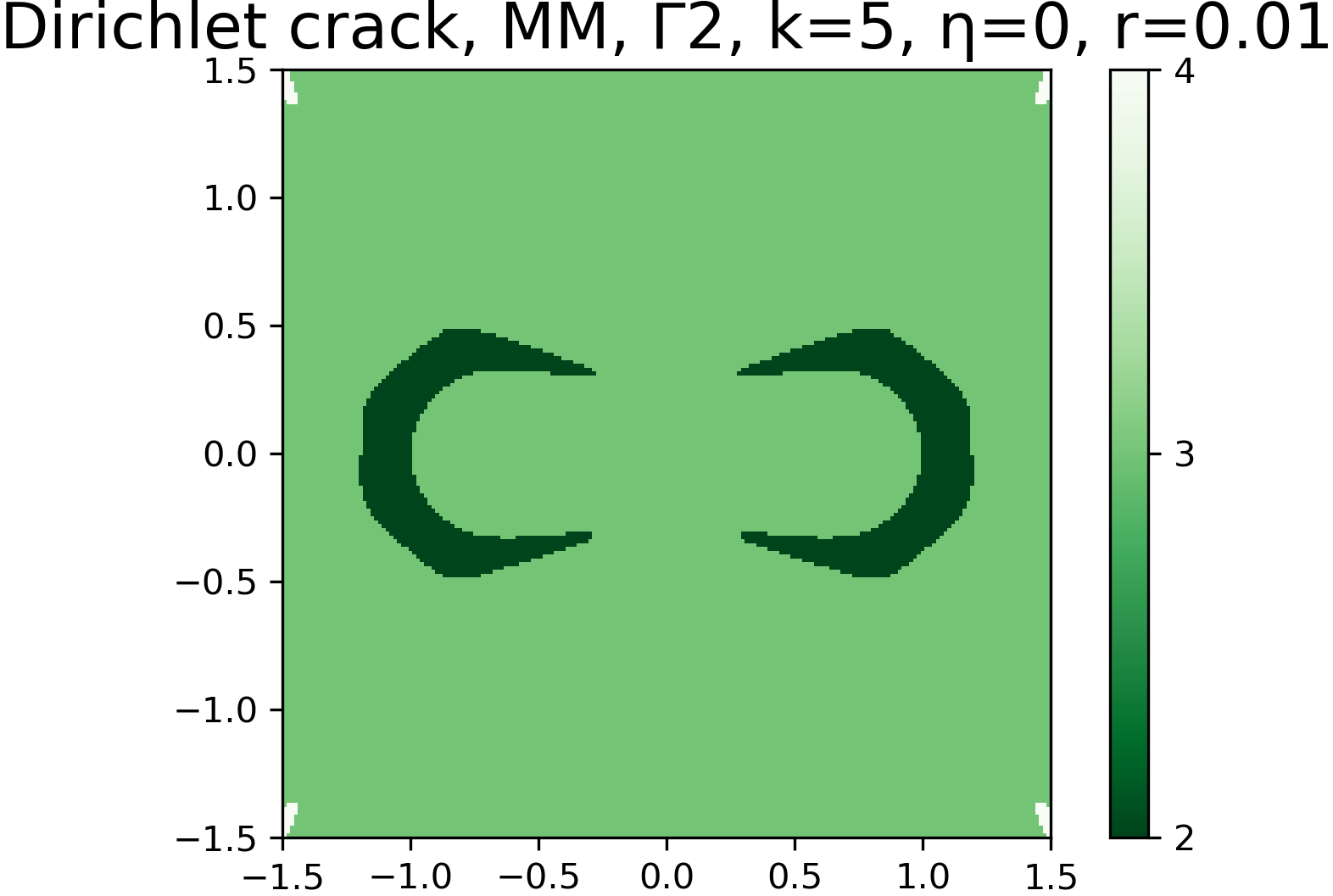}
 \end{center}
 \end{minipage}
 \vspace{1cm} \\ 
\begin{minipage}{0.5\hsize}
  \begin{center}
   \includegraphics[scale=0.5]{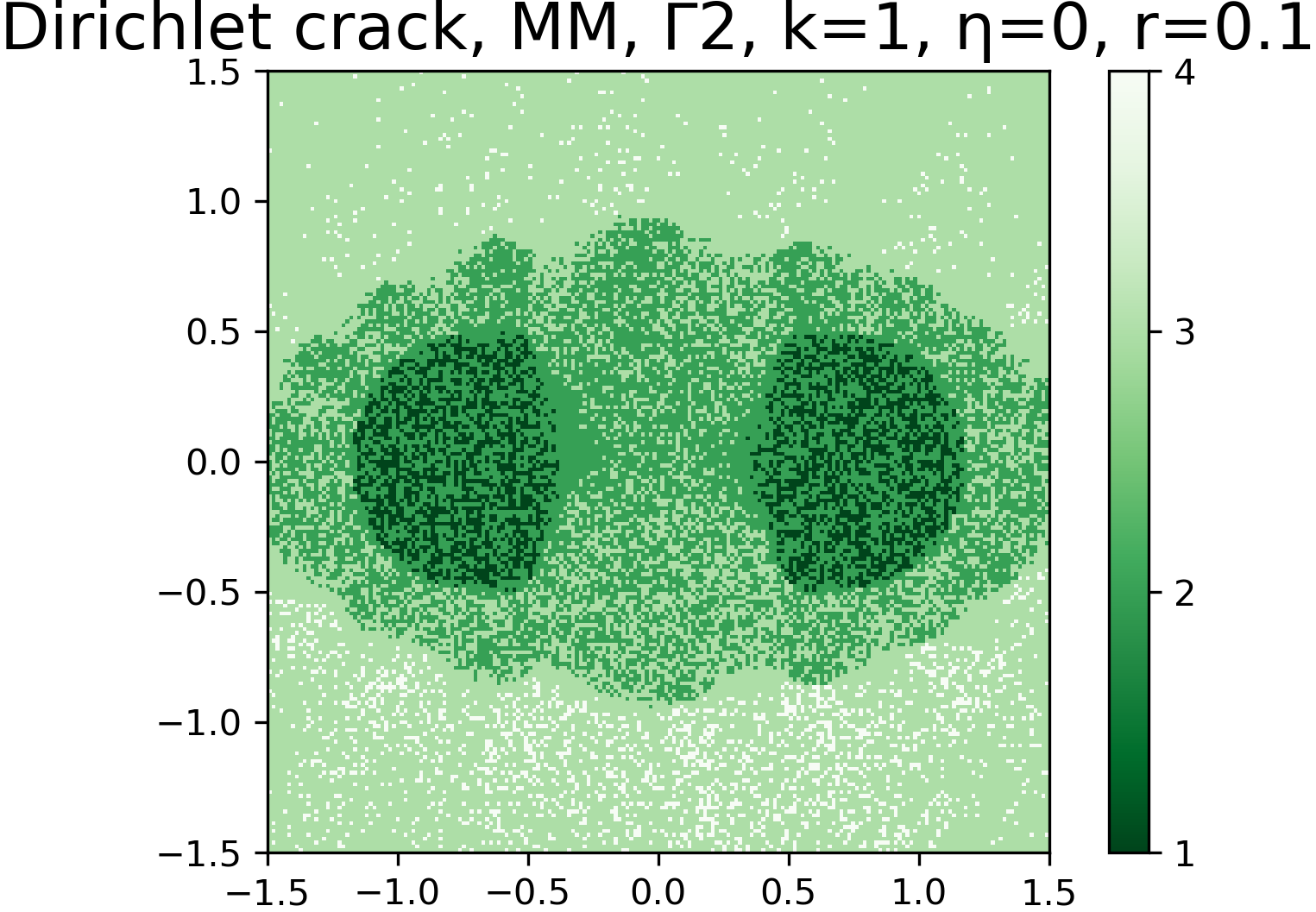}
  \end{center}
 \end{minipage}
 \begin{minipage}{0.5\hsize}
 \begin{center}
  \includegraphics[scale=0.5]{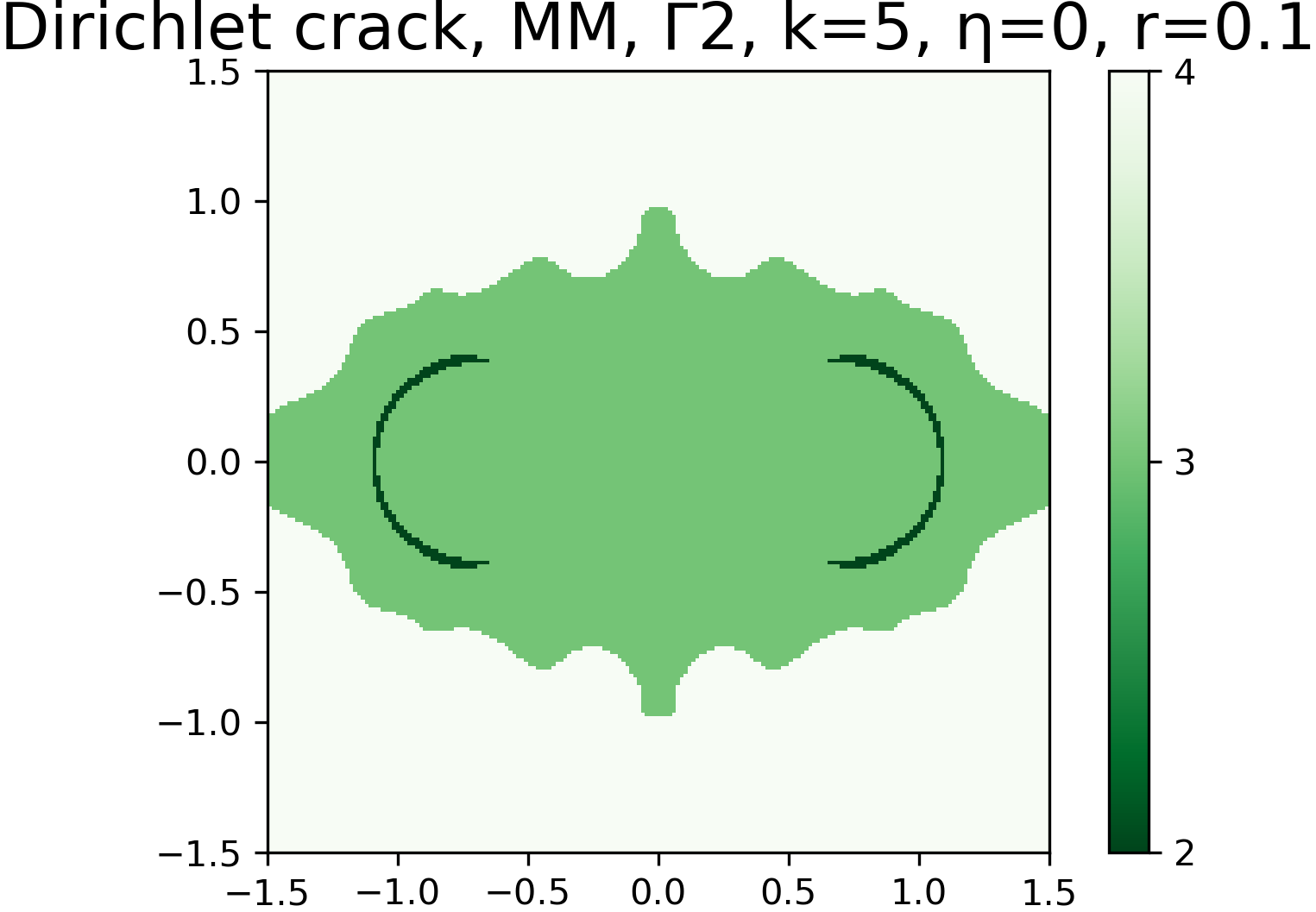}
 \end{center}
 \end{minipage}
 \vspace{1cm} \\ 
\begin{minipage}{0.5\hsize}
  \begin{center}
   \includegraphics[scale=0.5]{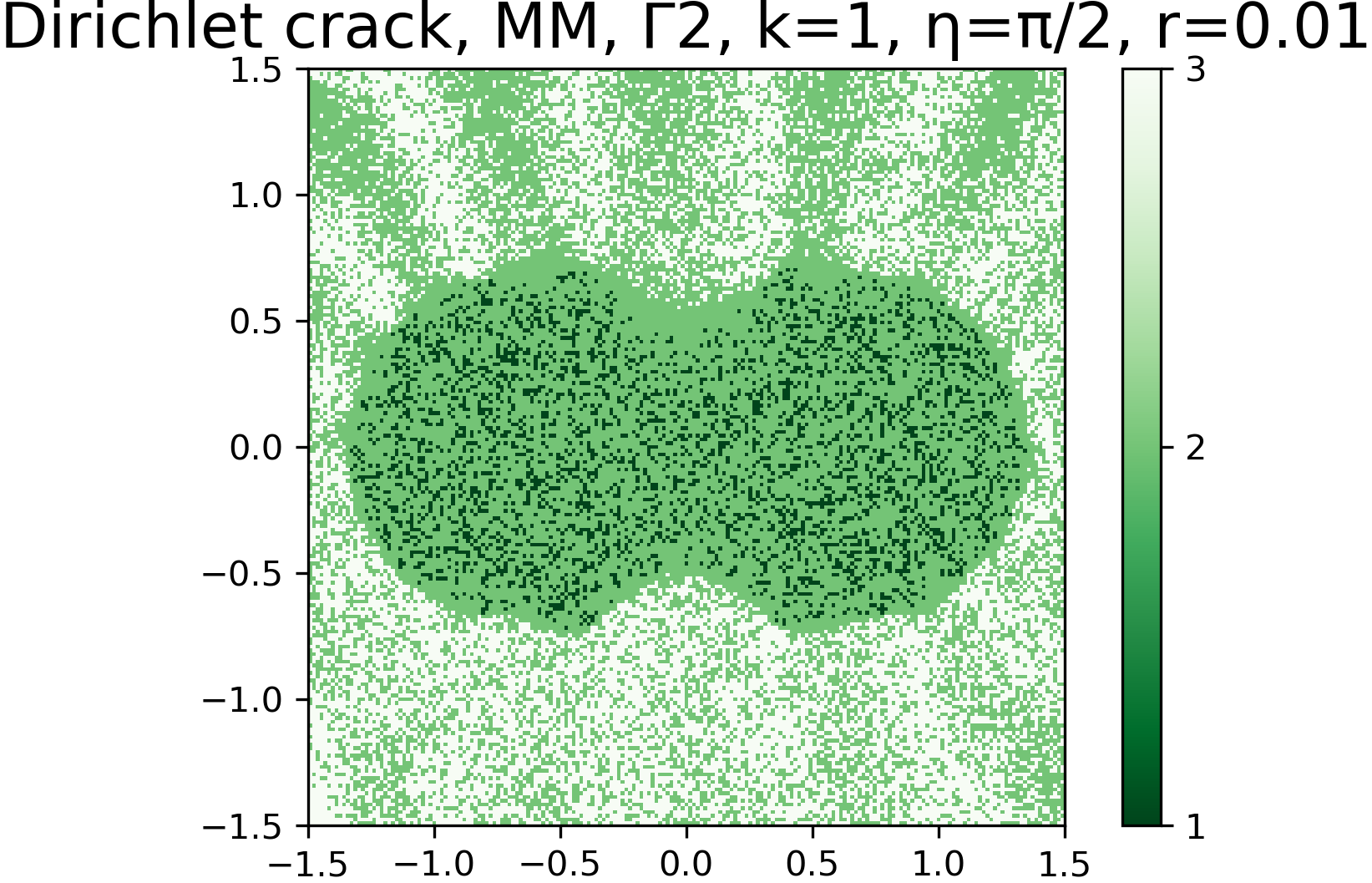}
  \end{center}
 \end{minipage}
 \begin{minipage}{0.5\hsize}
 \begin{center}
  \includegraphics[scale=0.5]{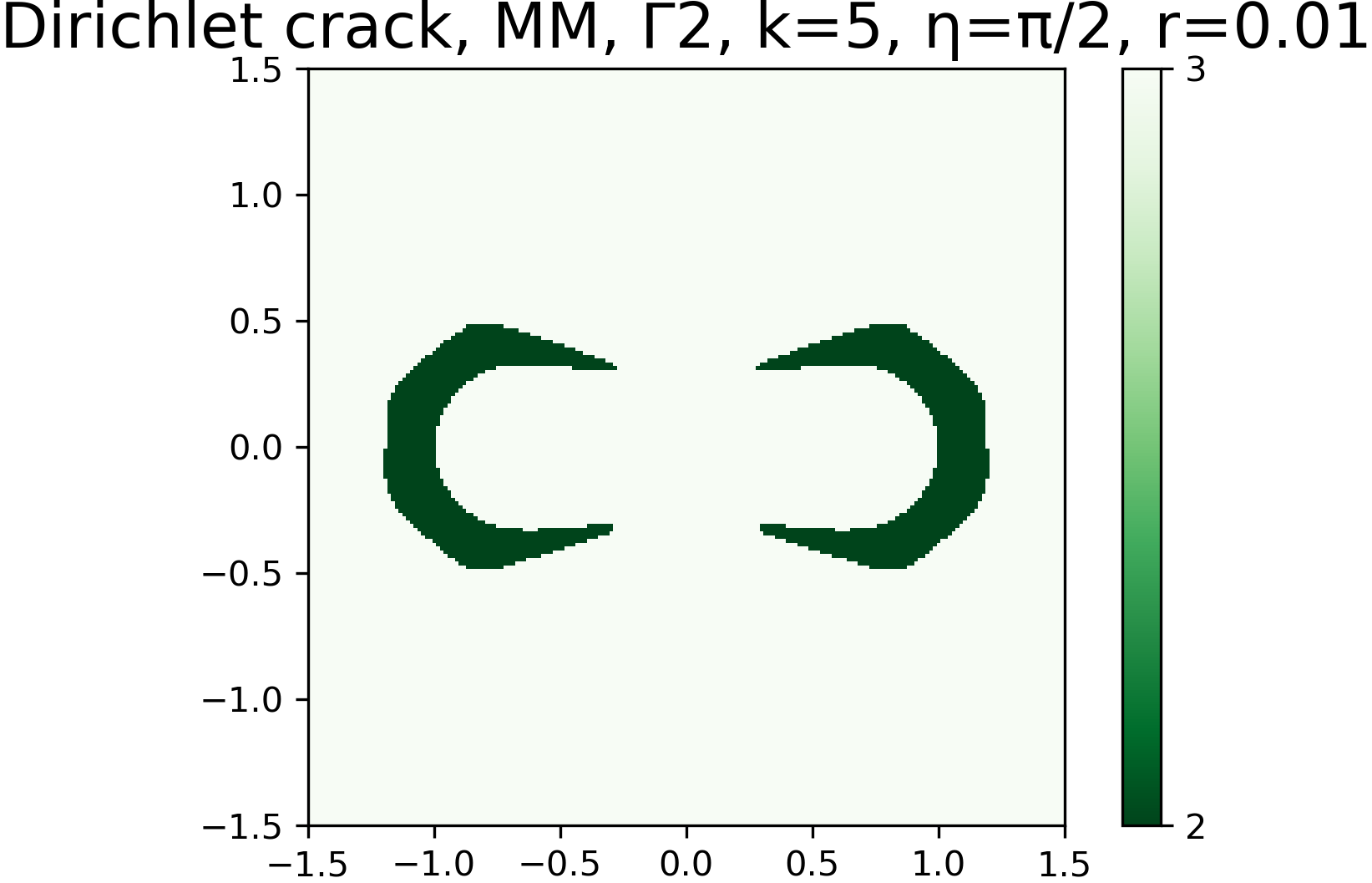}
 \end{center}
 \end{minipage}
 \vspace{1cm} \\ 
\begin{minipage}{0.5\hsize}
  \begin{center}
   \includegraphics[scale=0.5]{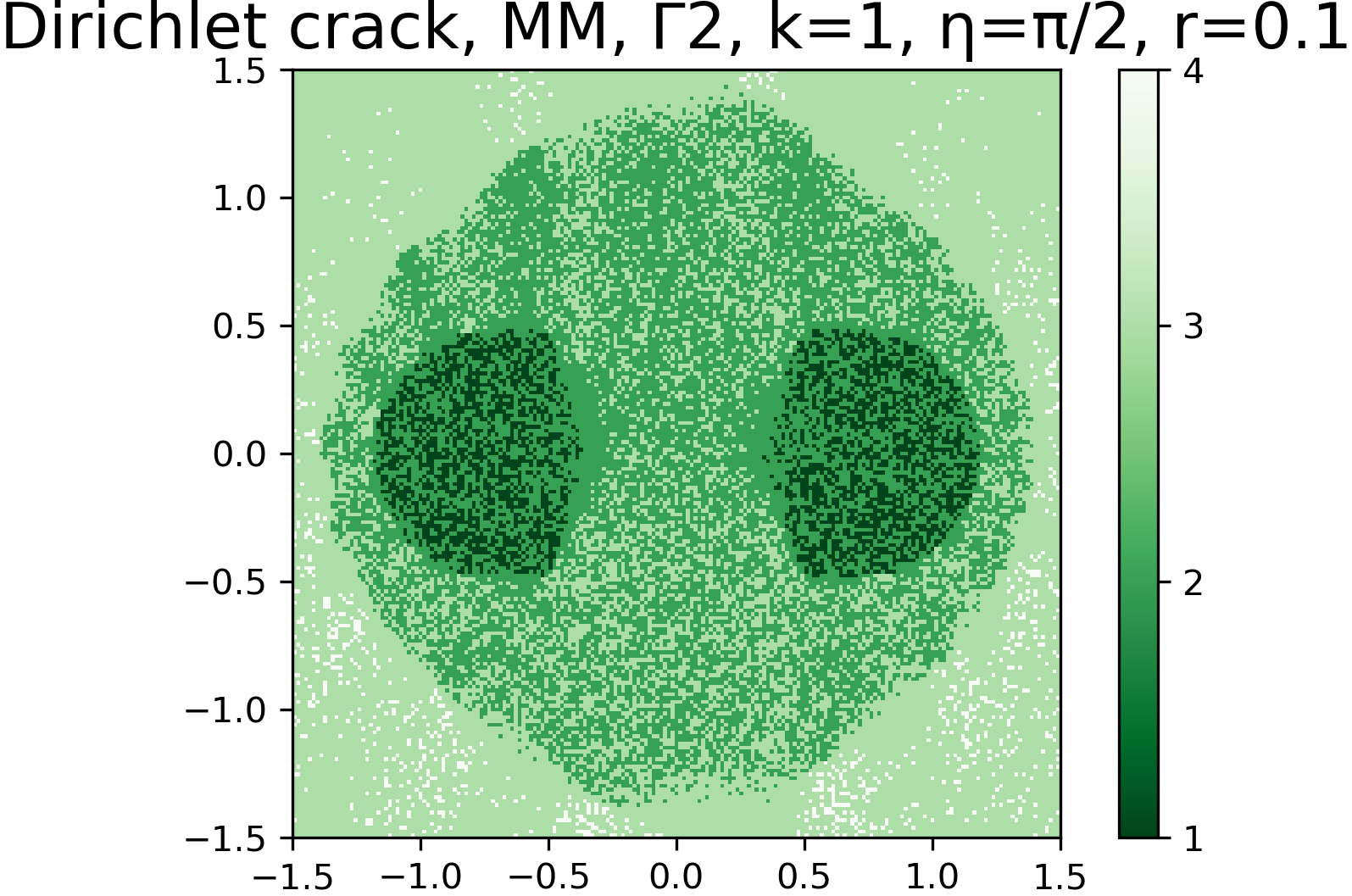}
  \end{center}
 \end{minipage}
 \begin{minipage}{0.5\hsize}
 \begin{center}
  \includegraphics[scale=0.5]{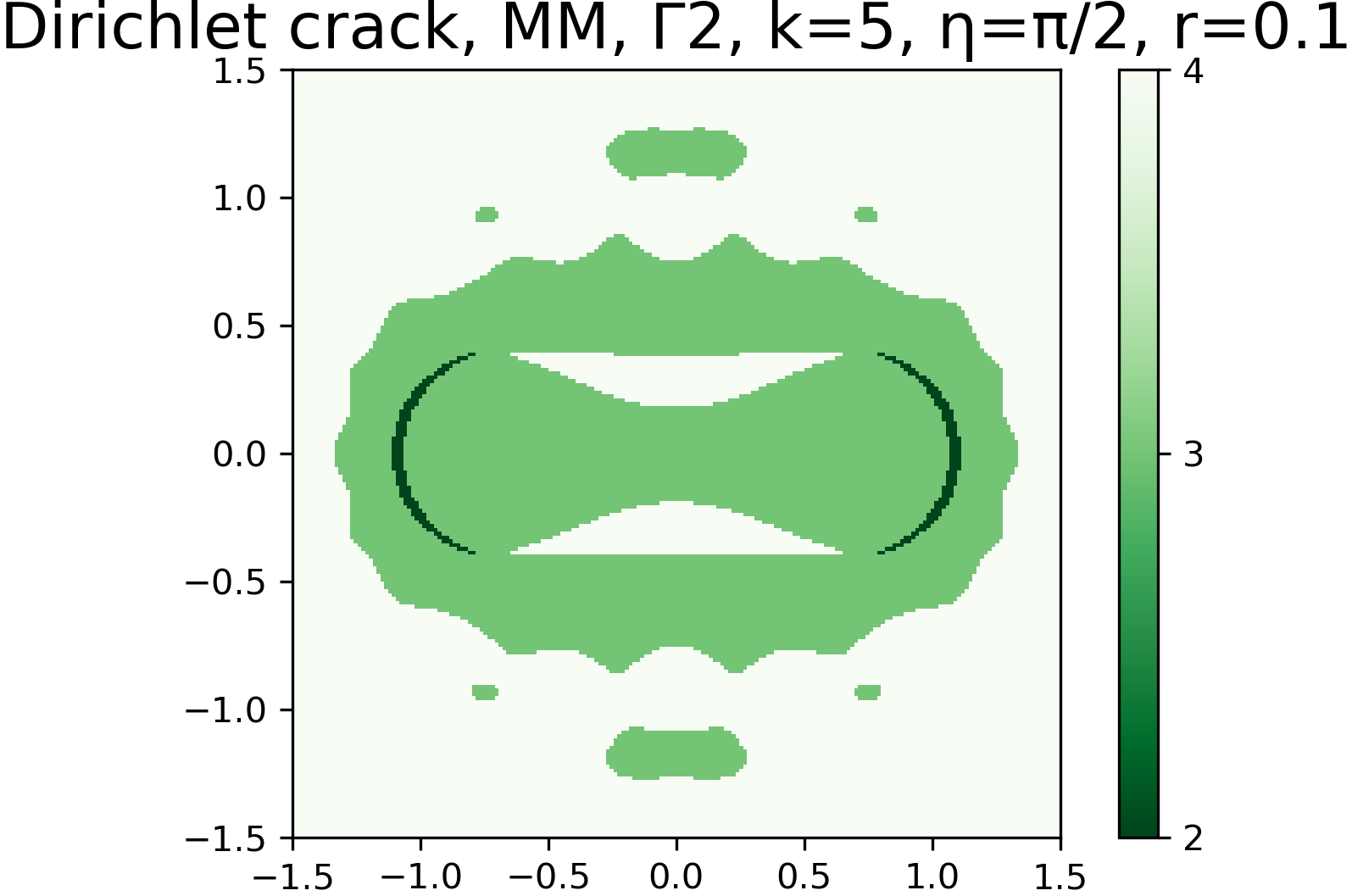}
 \end{center}
 \end{minipage}
 \vspace{1cm} \\ 
\end{tabular}
\caption{Reconstruction for the Dirichlet crack $\Gamma_2$ by the monotonicity method for different angles $\eta=0, \pi/2$, lengths $r=0.01, 0.1$, and wavenumbers $k=1,5$.}\label{MM crack-2}
\end{figure}
\begin{figure}[htbp]
\begin{tabular}{c}
\begin{minipage}{0.5\hsize}
  \begin{center}
   \includegraphics[scale=0.5]{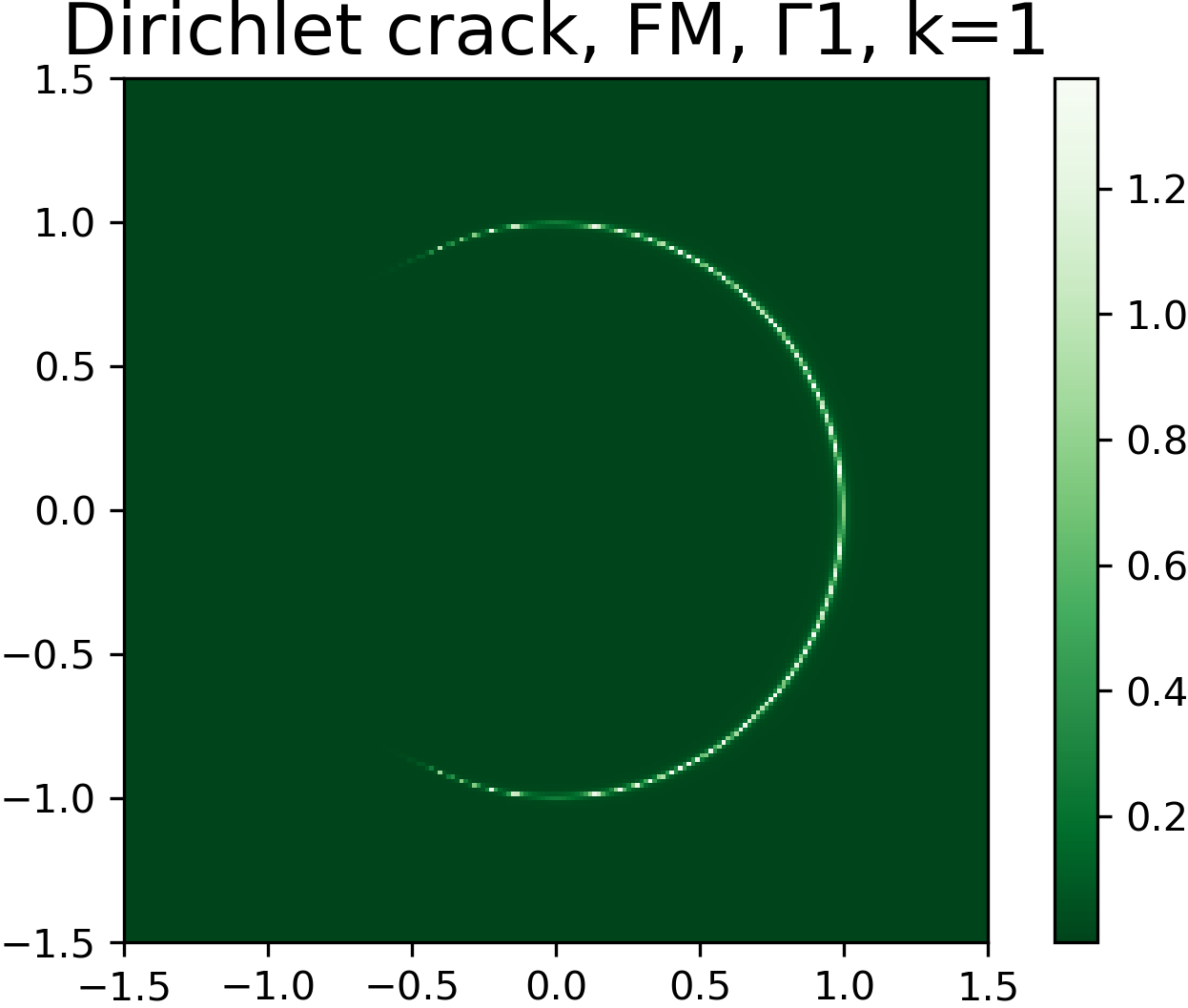}
  \end{center}
 \end{minipage}
 \begin{minipage}{0.5\hsize}
 \begin{center}
  \includegraphics[scale=0.5]{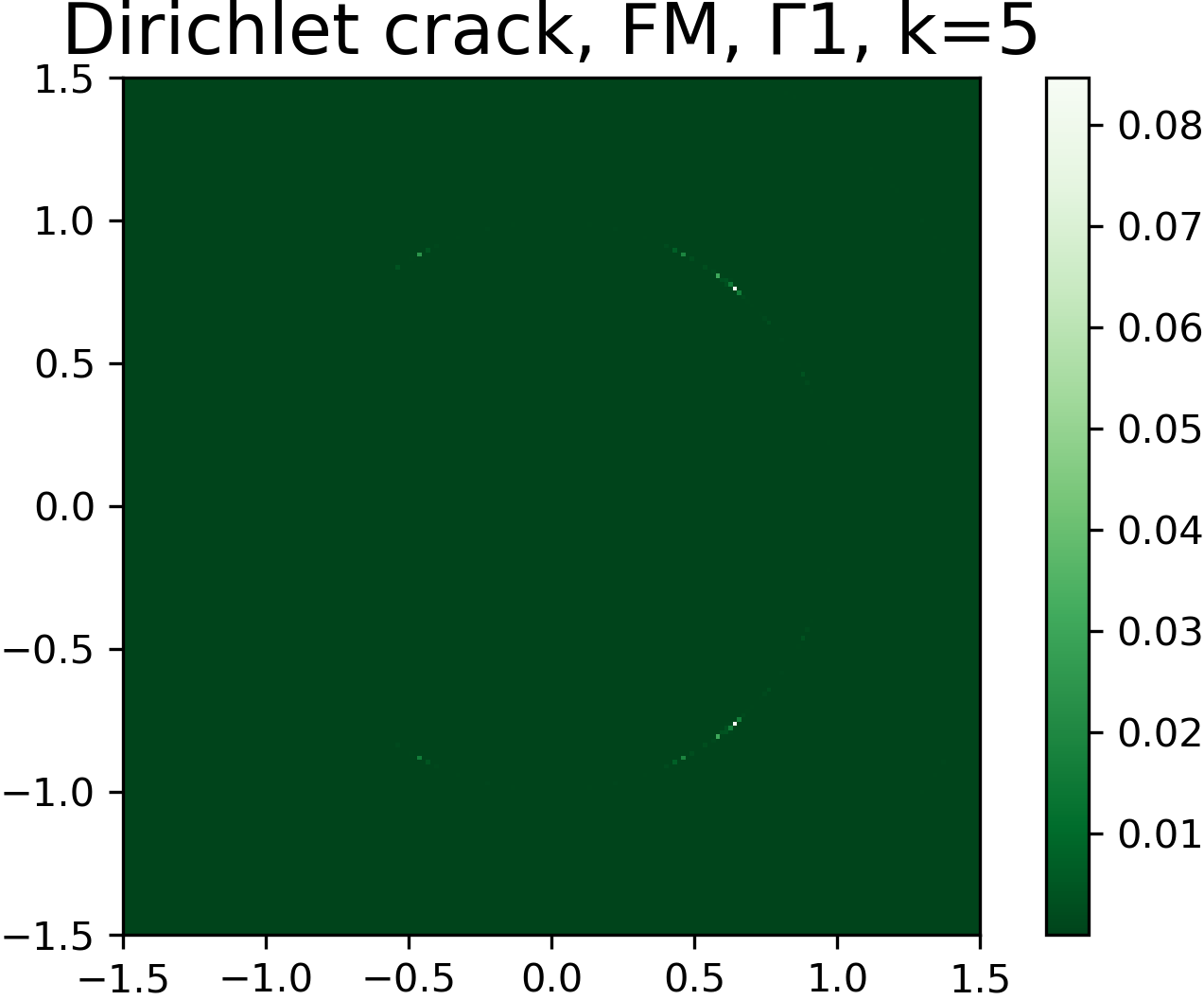}
 \end{center}
 \end{minipage}
 \vspace{1cm} \\ 
\begin{minipage}{0.5\hsize}
  \begin{center}
   \includegraphics[scale=0.5]{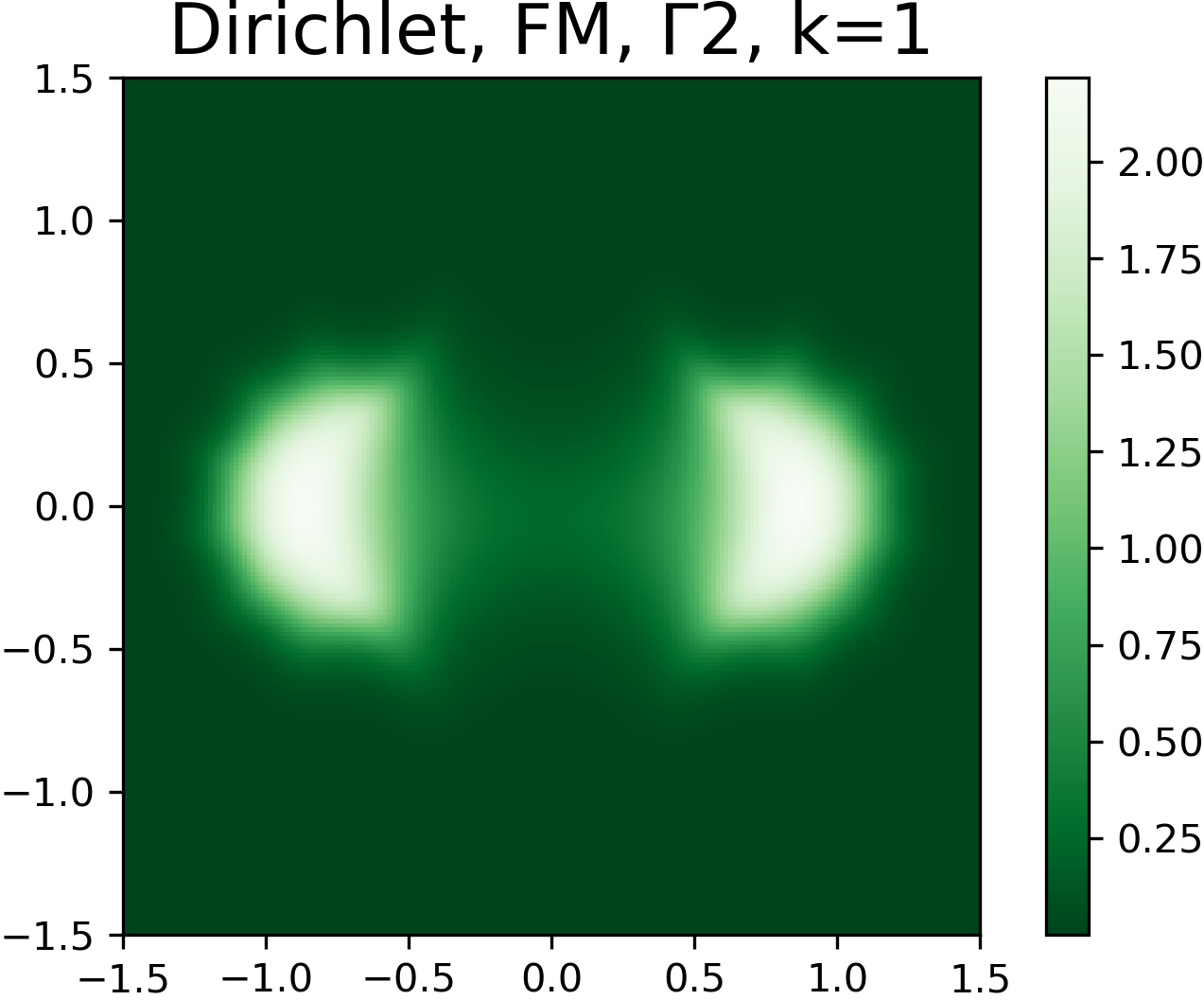}
  \end{center}
 \end{minipage}
 \begin{minipage}{0.5\hsize}
 \begin{center}
  \includegraphics[scale=0.5]{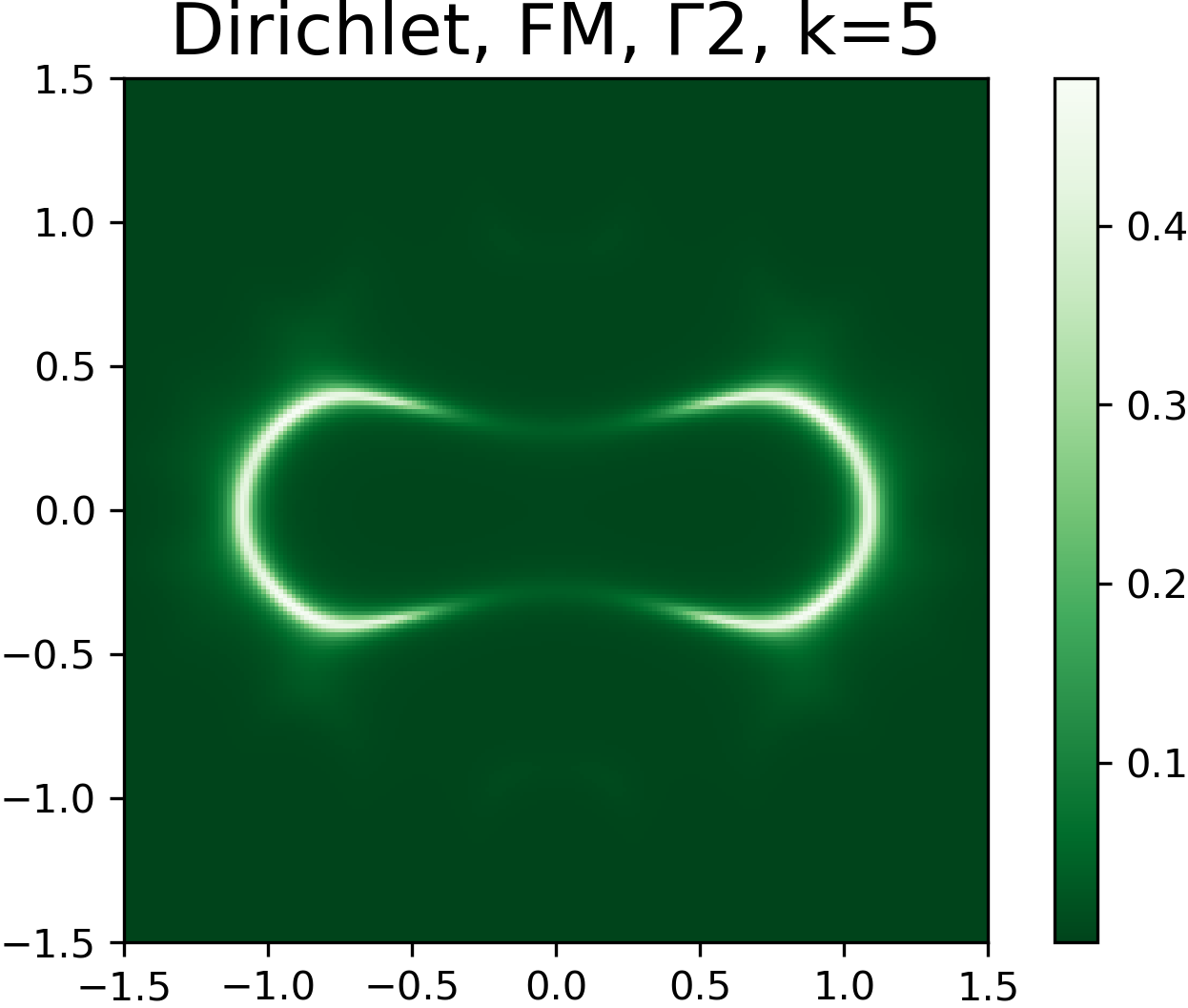}
 \end{center}
 \end{minipage}
\end{tabular}
\caption{Reconstruction for the Dirichlet crack by the factorization method for different wavenumbers $k=1,5$ and shapes $\Gamma_1$, $\Gamma_2$.}\label{FM crack}
\end{figure}

We also provided several numerical examples to compare the factorization method with the monotonicity method. The factorization method is a {\it point test} which checks whether a point $z$ is included in the unknown target or not, while the monotonicity is a {\it domain test} which checks whether a domain $B$ is included in the unknown target or not. For the domain test, we have to find the appropriate choice of $B$. By testing the monotonicity for many different shapes and sizes of $B$, we obtained as accurate reconstructions as the factorization (see Figures \ref{MM obstacle}, \ref{MM medium}, \ref{FM obstacle}, \ref{FM medium}, \ref{MM crack-1}, \ref{MM crack-2}, and \ref{FM crack}). 
\par
\begin{figure}[htbp]
\vspace{-3cm}
\begin{tabular}{c}
\begin{minipage}{0.5\hsize}
  \begin{center}
   \includegraphics[scale=0.5]{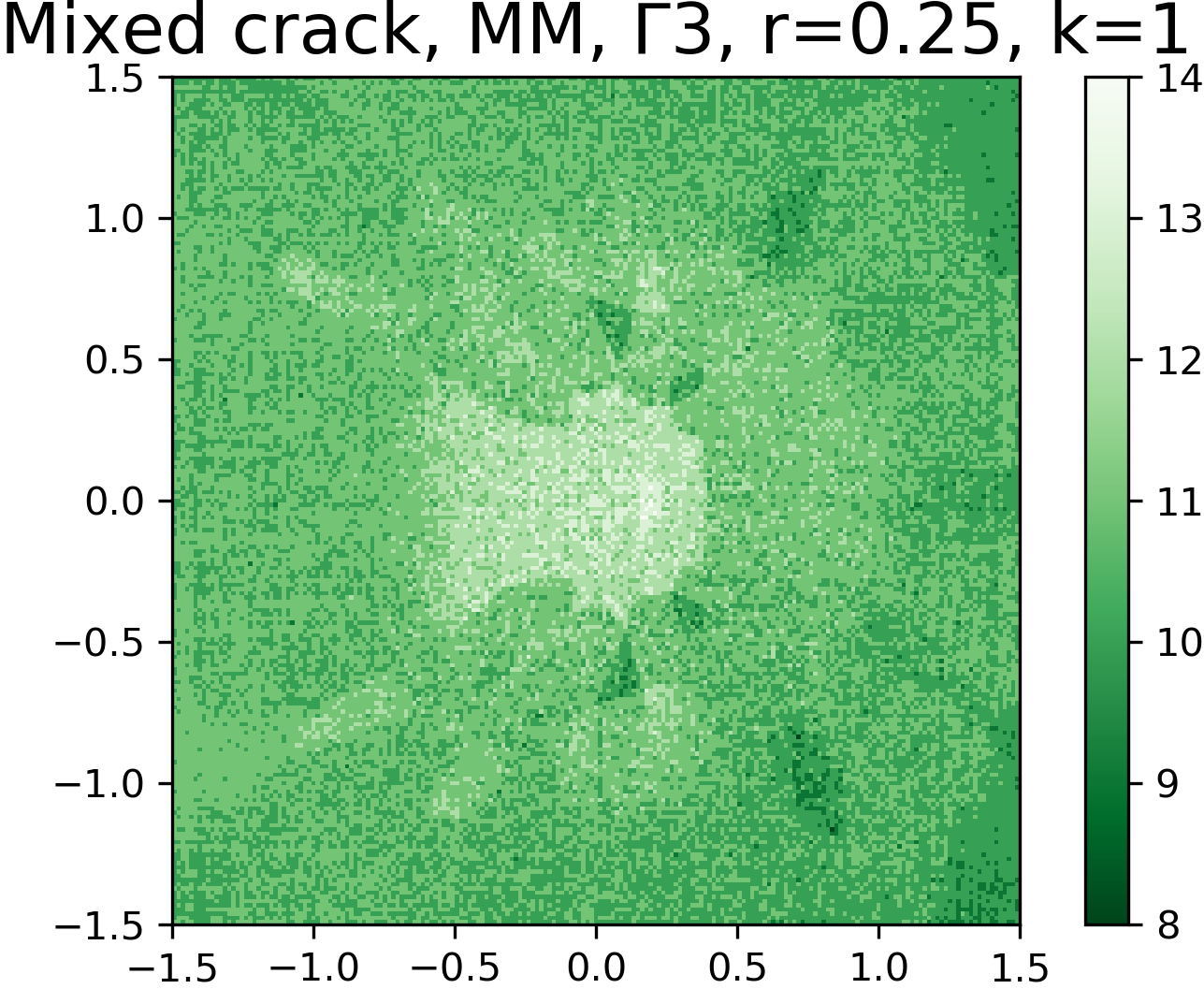}
  \end{center}
 \end{minipage}
 \begin{minipage}{0.5\hsize}
 \begin{center}
  \includegraphics[scale=0.5]{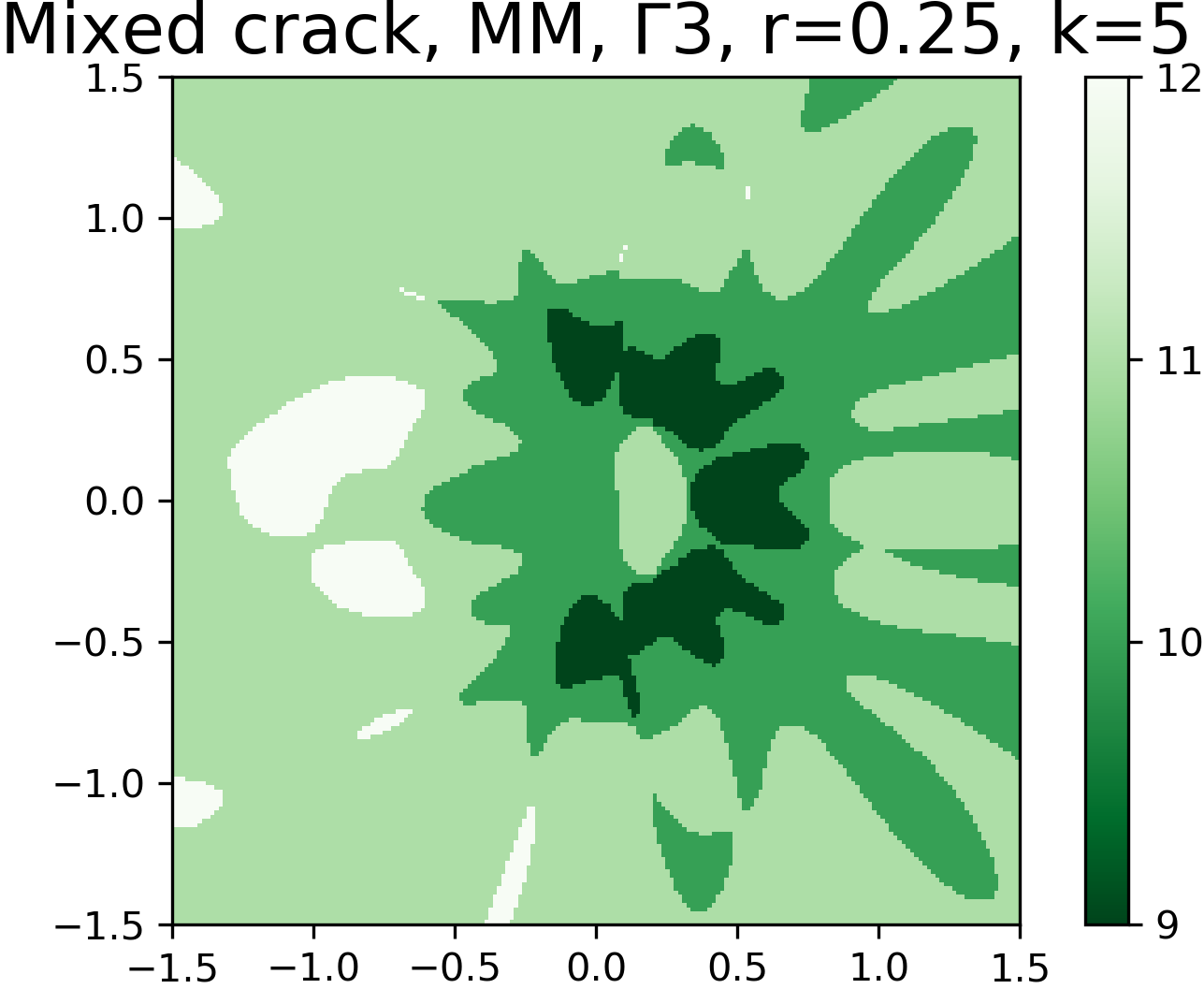}
 \end{center}
 \end{minipage}
 \vspace{1cm} \\ 
\begin{minipage}{0.5\hsize}
  \begin{center}
   \includegraphics[scale=0.5]{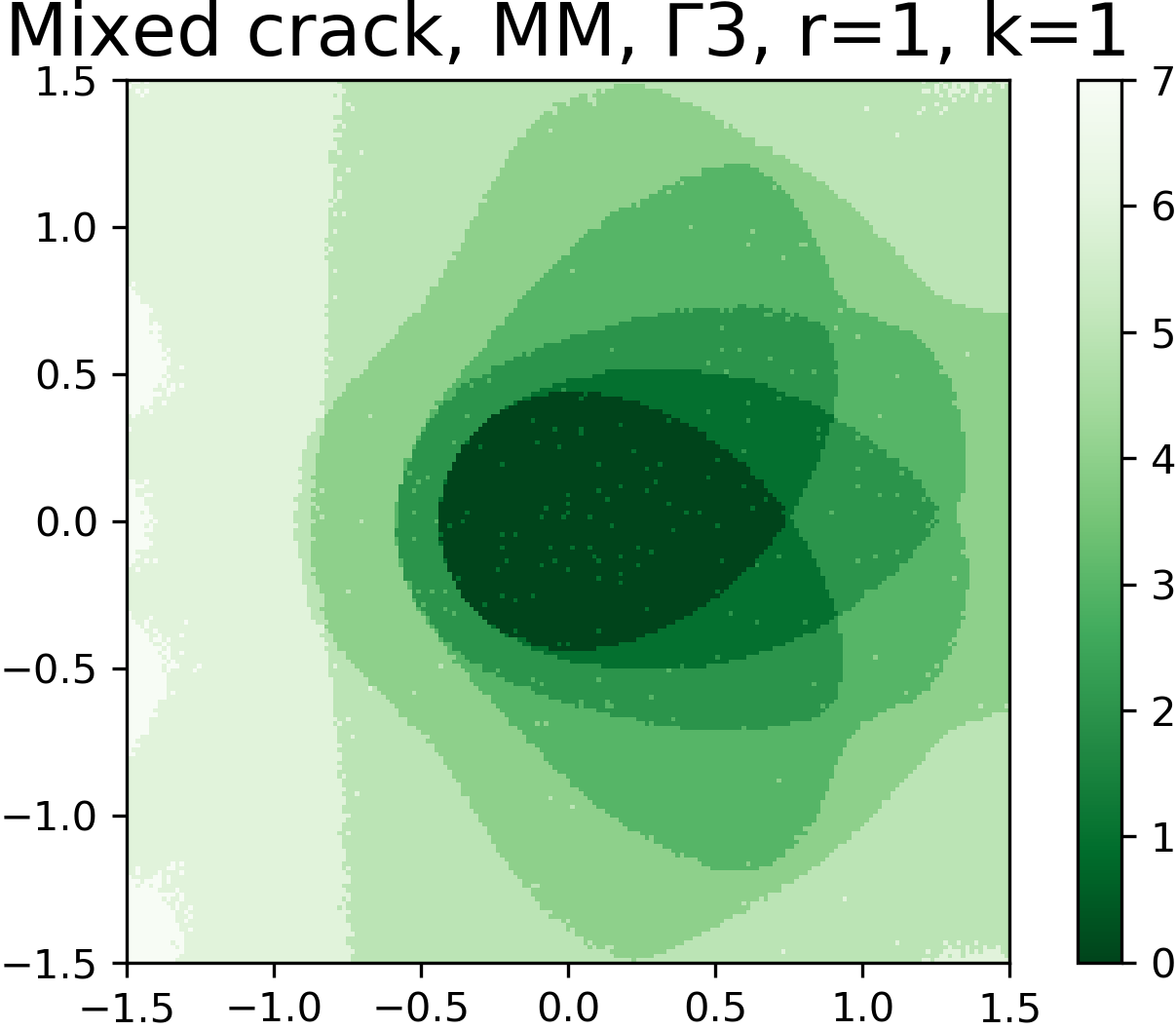}
  \end{center}
 \end{minipage}
 \begin{minipage}{0.5\hsize}
 \begin{center}
  \includegraphics[scale=0.5]{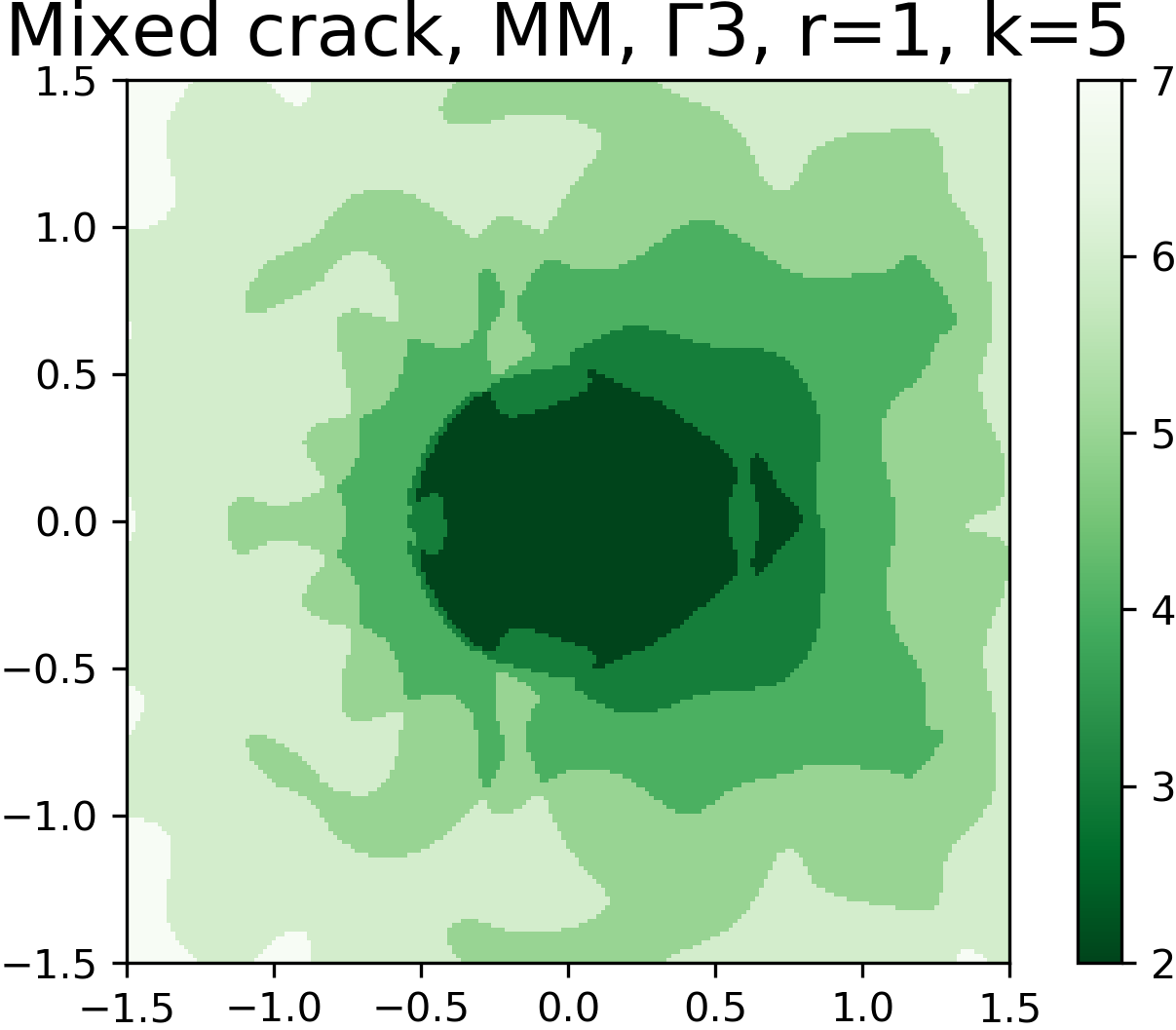}
 \end{center}
 \end{minipage}
 \vspace{1cm} \\ 
\end{tabular}
\caption{Reconstruction for the mixed crack $\Gamma_3$ by the shifting circle test of monotonicity method for different radiuses $r=0.25, 1$ and wavenumbers $k=1, 5$.}\label{MM mixed crack-shifting}
\end{figure}
\begin{figure}[htbp]
\vspace{-3cm}
\begin{tabular}{c}
\begin{minipage}{0.5\hsize}
  \begin{center}
   \includegraphics[scale=0.5]{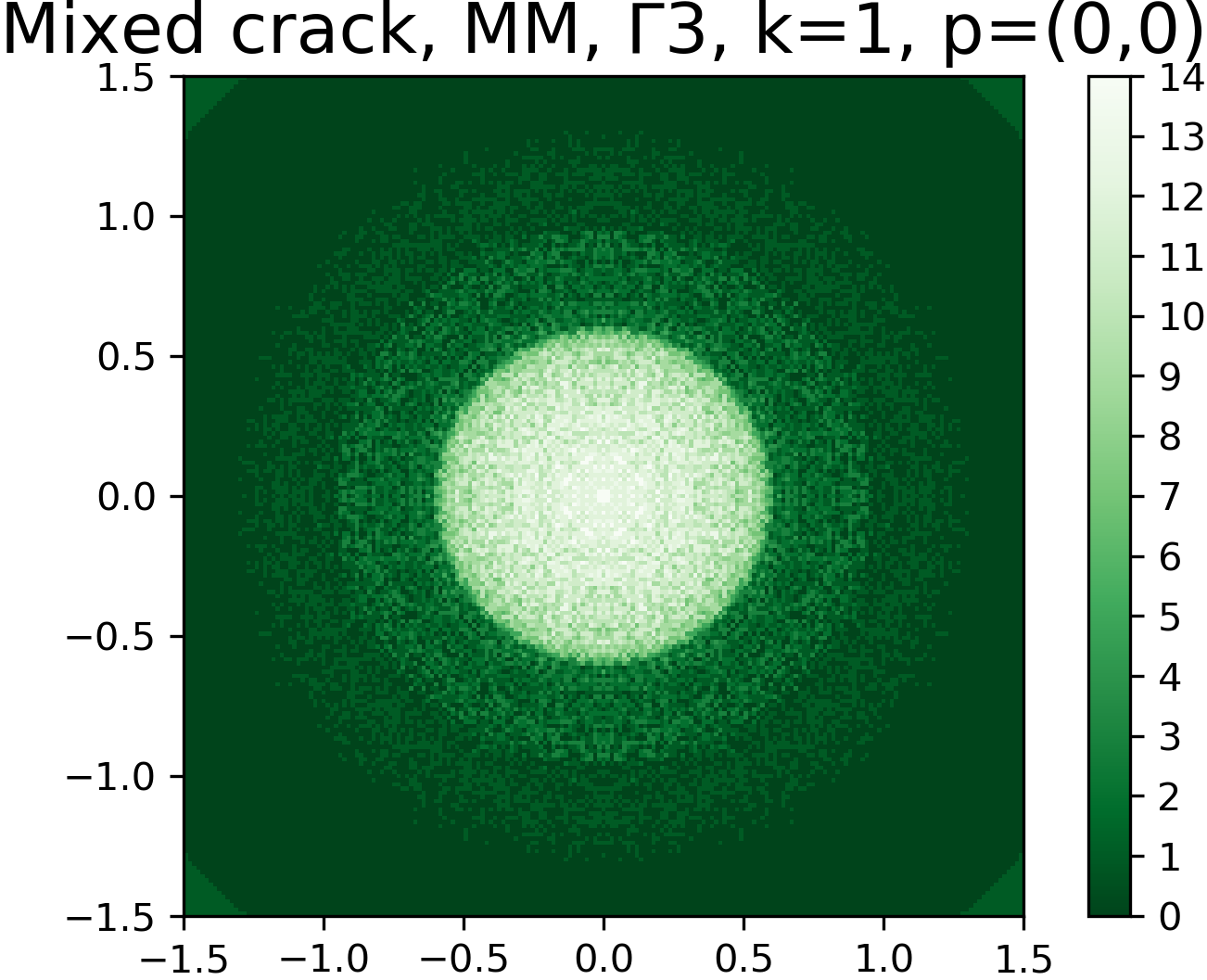}
  \end{center}
 \end{minipage}
 \begin{minipage}{0.5\hsize}
 \begin{center}
  \includegraphics[scale=0.5]{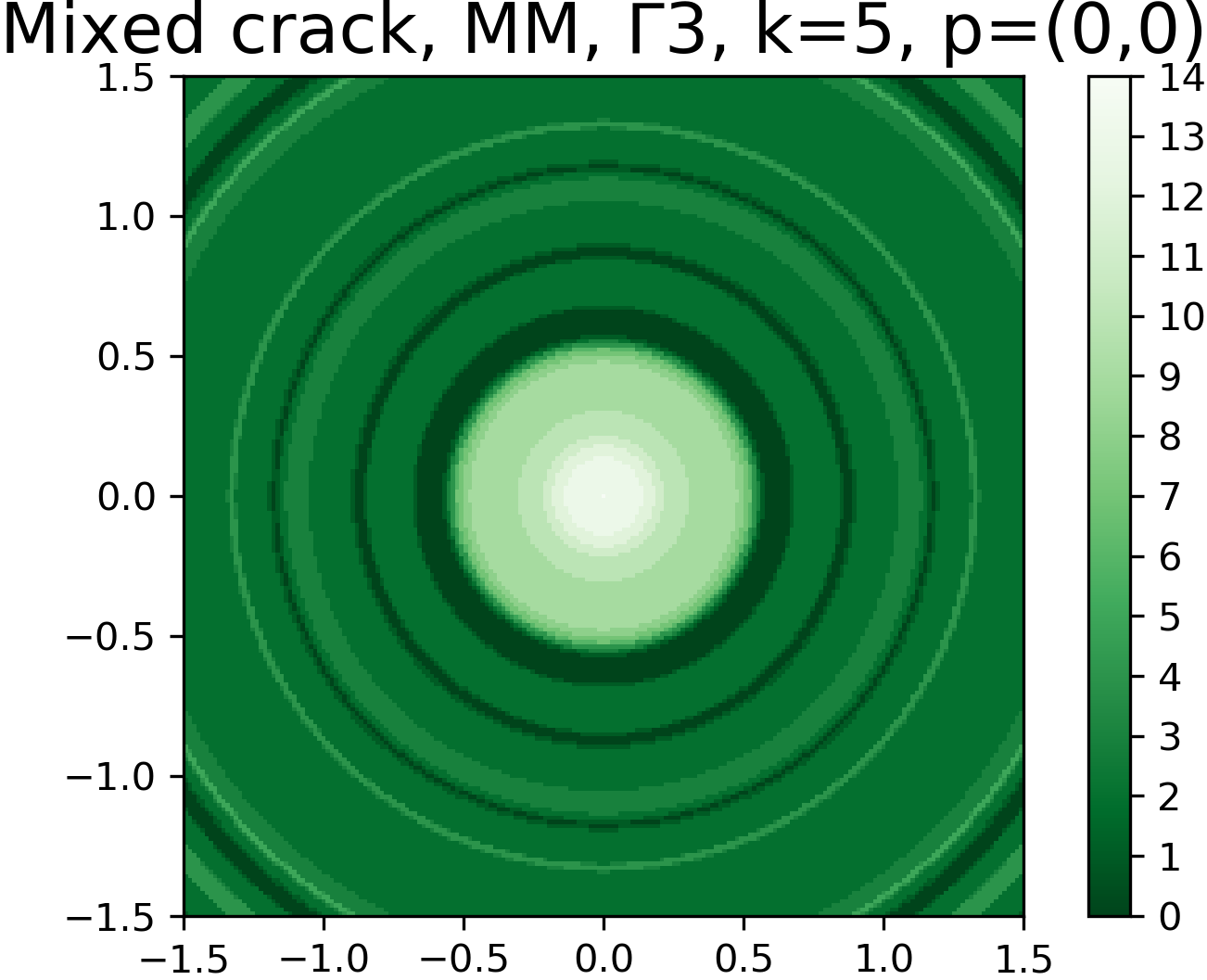}
 \end{center}
 \end{minipage}
 \vspace{1cm} \\ 
\begin{minipage}{0.5\hsize}
  \begin{center}
   \includegraphics[scale=0.5]{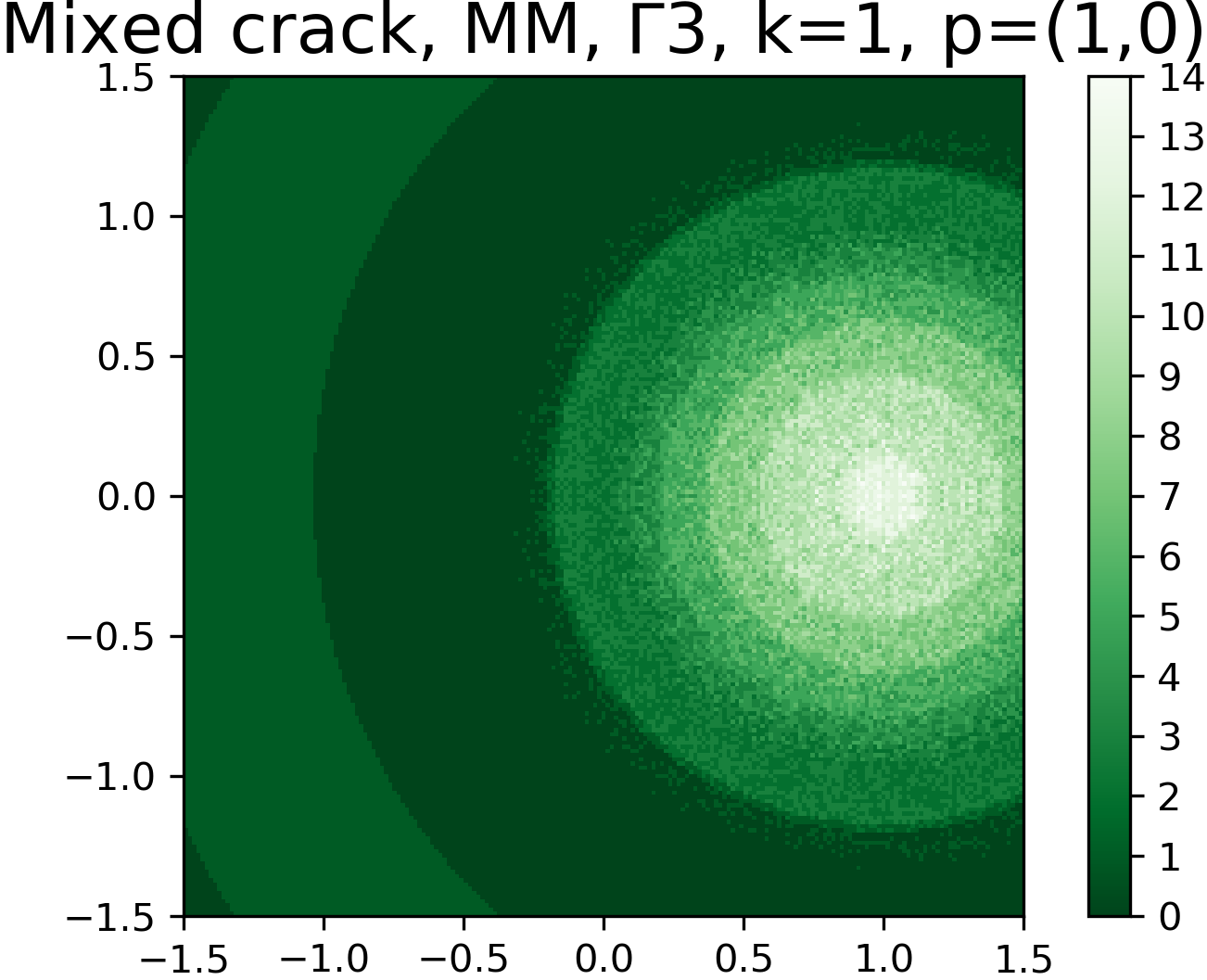}
  \end{center}
 \end{minipage}
 \begin{minipage}{0.5\hsize}
 \begin{center}
  \includegraphics[scale=0.5]{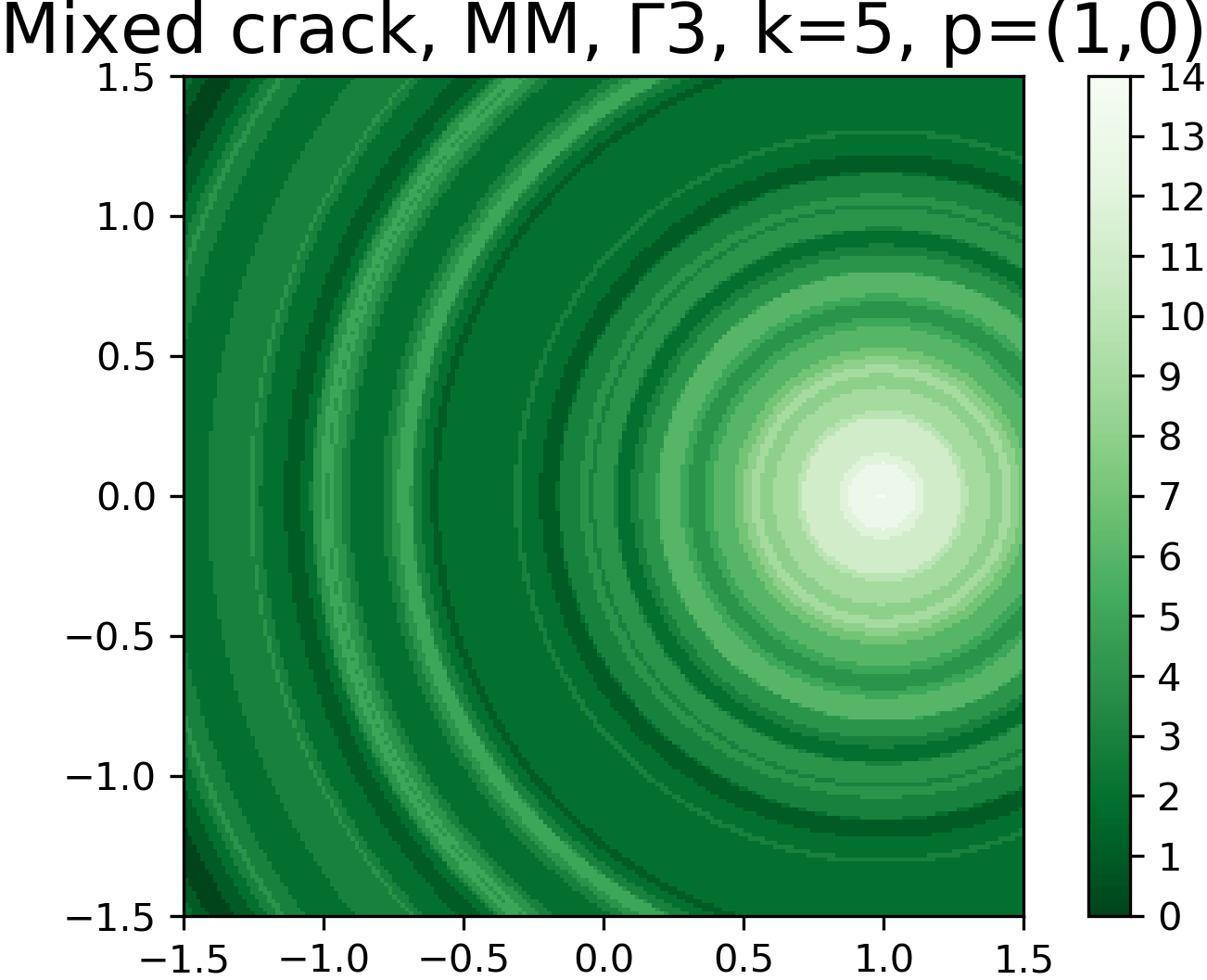}
 \end{center}
 \end{minipage}
 \vspace{1cm} \\ 
\end{tabular}
\caption{Reconstruction for the mixed crack $\Gamma_3$ by the shrinking circle test of monotonicity method for different points $p=(0,0), (1,0)$ and wavenumbers $k=1, 5$.}\label{MM mixed crack-shrinking}
\end{figure}
\begin{figure}[htbp]
\begin{tabular}{c}
\begin{minipage}{0.5\hsize}
  \begin{center}
   \includegraphics[scale=0.5]{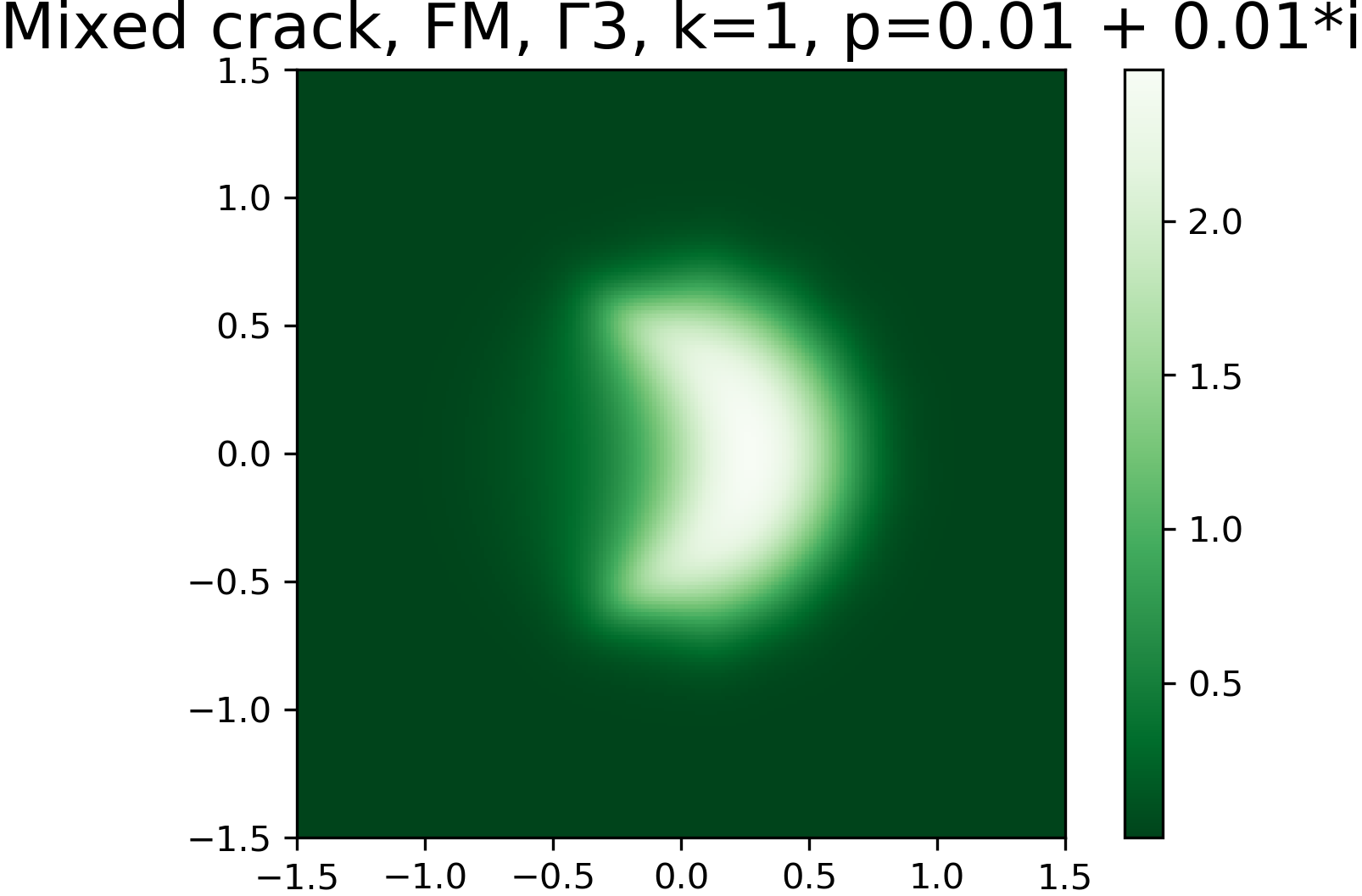}
  \end{center}
 \end{minipage}
 \begin{minipage}{0.5\hsize}
 \begin{center}
  \includegraphics[scale=0.5]{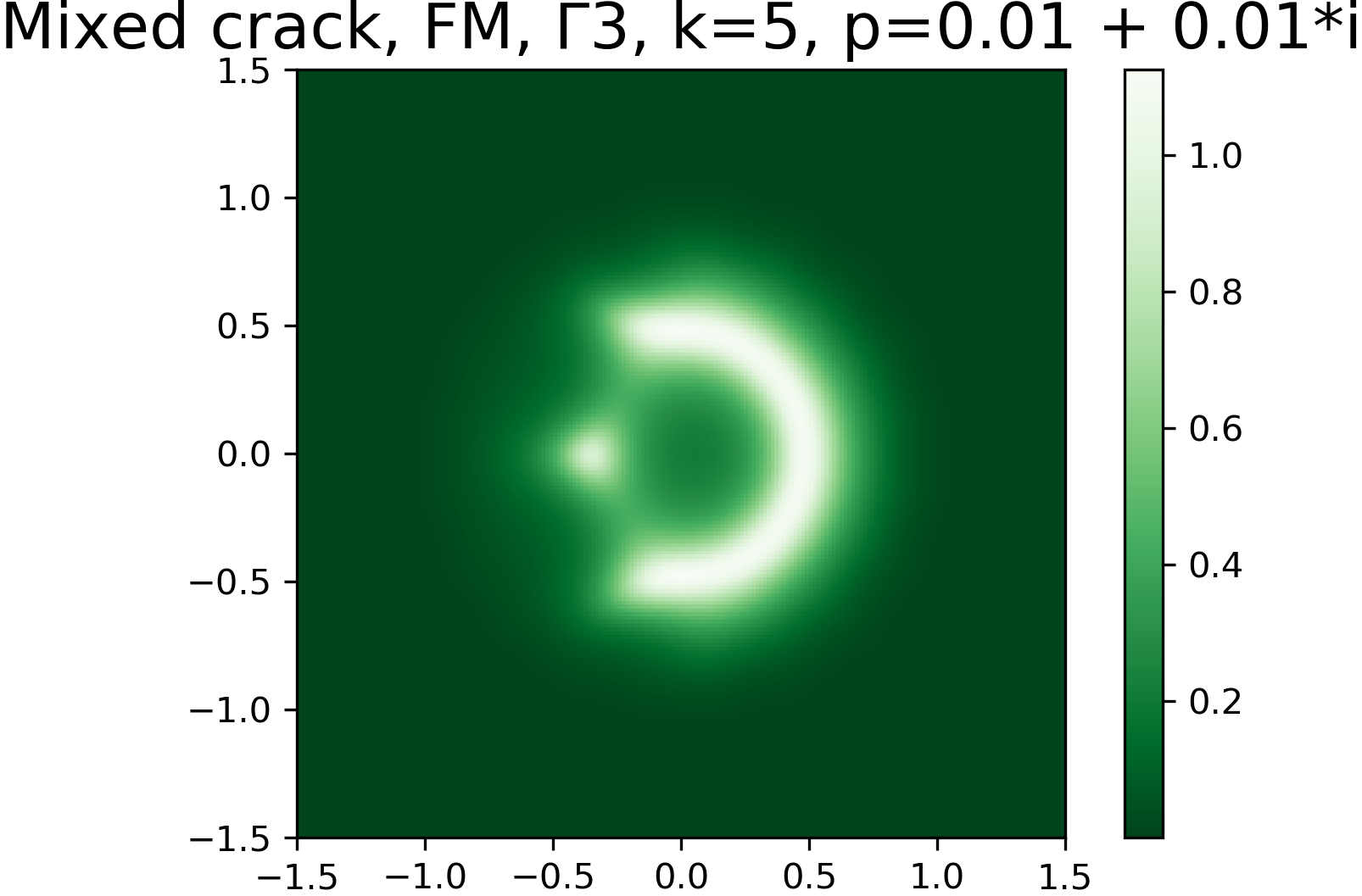}
 \end{center}
 \end{minipage}
 \vspace{1cm} \\ 
\begin{minipage}{0.5\hsize}
  \begin{center}
   \includegraphics[scale=0.5]{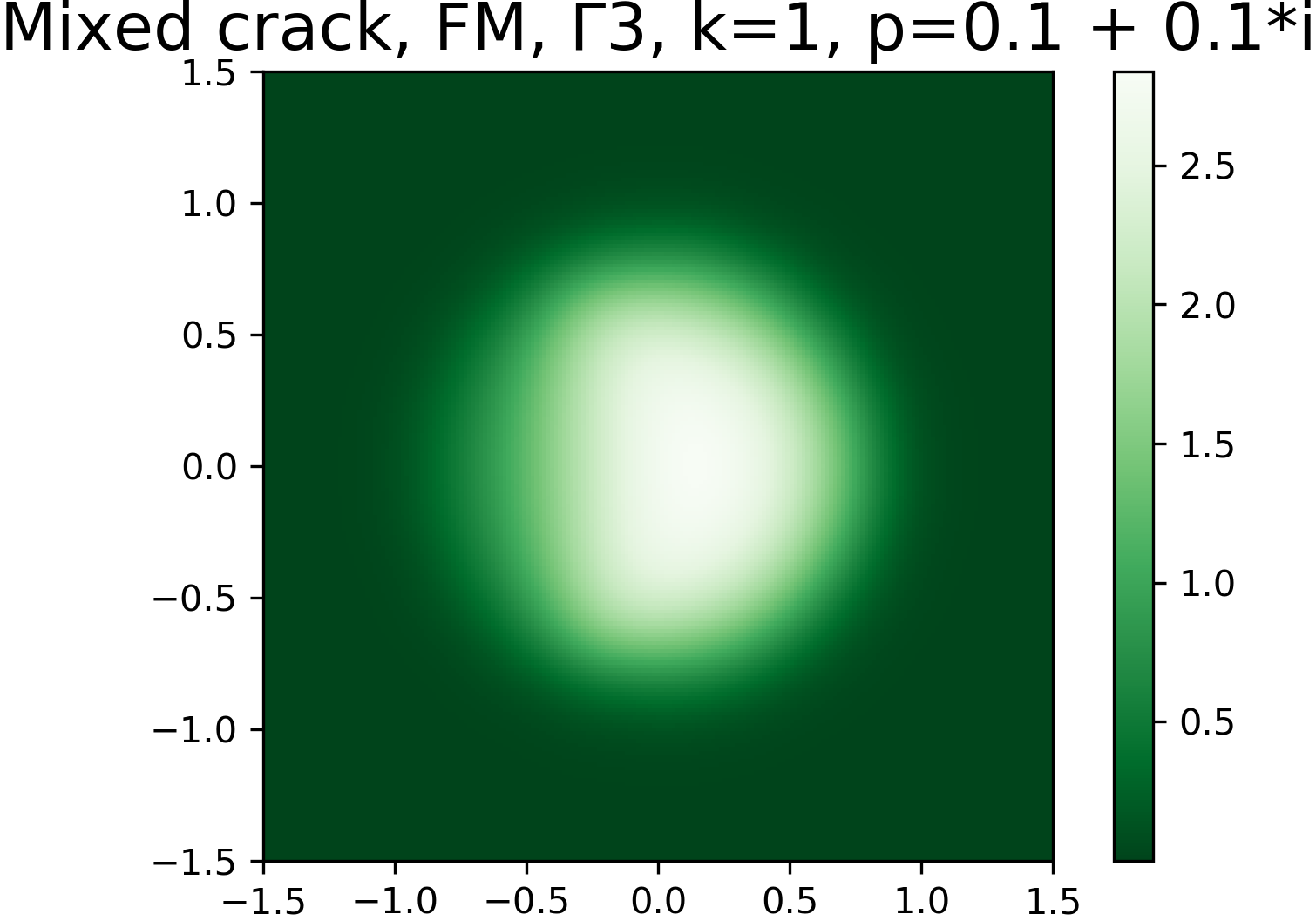}
  \end{center}
 \end{minipage}
 \begin{minipage}{0.5\hsize}
 \begin{center}
  \includegraphics[scale=0.5]{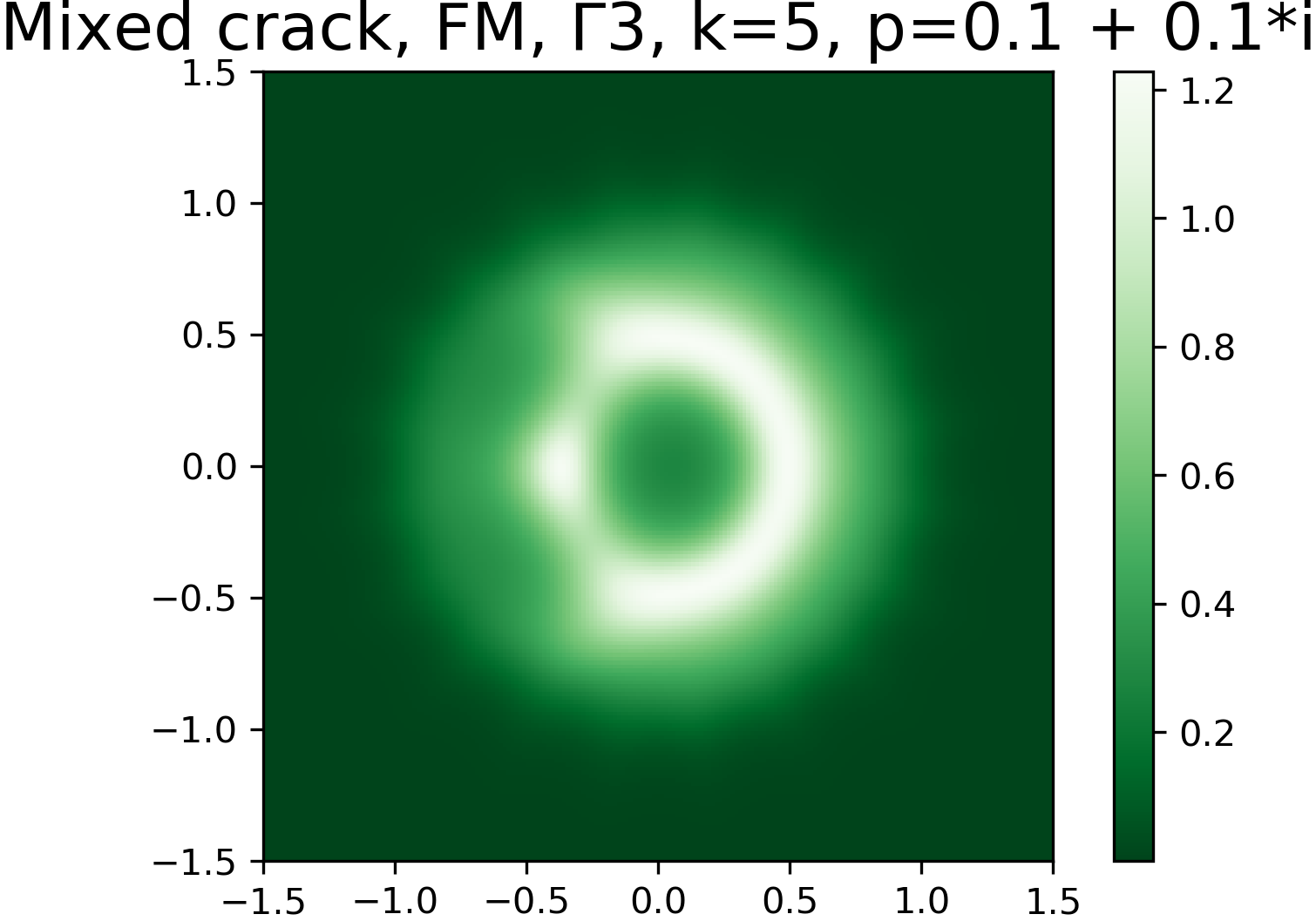}
 \end{center}
 \end{minipage}
\end{tabular}
\caption{Reconstruction for the mixed crack $\Gamma_3$ by the factorization method for different wavenumbers $k=1,5$ and complex numbers $p=0.01+0.01i, 0.1+0.1i$.}\label{FM mixed crack}
\end{figure}

In our numerical examples for the mixture crack, we observed that the factorization method is more accurate than the monotonicity (see Figures \ref{MM mixed crack-shifting}, \ref{MM mixed crack-shrinking}, \ref{FM mixed crack}). This comes from that the factorization method used the Herglotz operator corresponding to the auxiliary closed curve that is a extension of the unknown crack (see (\ref{auxiliary})), which is a very restrictive assumption. On the other hand, the monotonicity can give the reconstruction scheme without using the auxiliary closed curve. Its scheme is the {\it outside domain test} that checks whether a domain $B$ contains $\Gamma$ or not. However in our numerical examples, the outside domain test can not obtain the exact shape of the unknown crack $\Gamma$ although the location and size can be done (see Figures \ref{MM mixed crack-shifting}, \ref{MM mixed crack-shrinking}). As we see numerical examples for the mixed obstacle discussed in \cite{Albicker and Griesmaier}, it also seems to be difficult to obtain the shape of unknown obstacle from outside domain tests. The numerical tests for the mixed problem in the monotonicity has to be developed in the future.

\section*{Acknowledgments}
This work was supported by JSPS KAKENHI Grant Number JP19J10238.

\vspace{5mm}
e-mail: takashi.furuya0101@gmail.com
\vspace{2cm}

\end{document}